\documentclass [11pt,a4paper,titlepage,twoside,openright]{report}
\usepackage{amsfonts,amssymb,amsmath,amscd,latexsym,epsfig}
\usepackage{fancyhdr}
\author{Luca M.A. Martinazzi}

\newtheorem{trm}{Theorem}[chapter]
\newtheorem{prop}[trm]{Proposition}
\newtheorem{cor}[trm]{Corollary}
\newtheorem{lemma}[trm]{Lemma}
\newtheorem{defin}[trm]{Definition}

\newcommand{\R}[1]{\mathbb{R}^{#1}}
\newcommand{\Cl}{\overline}
\newcommand{\M}[1]{\mathcal{#1}}
\newcommand{\norm}[1]{\left\Arrowvert{#1}\right\Arrowvert}
\newcommand{\abs}[1]{\left\arrowvert{#1}\right\arrowvert}
\newcommand{\bra}[1]{\left({#1}\right)}
\newcommand{\bbra}[1]{\left\{{#1}\right\}}
\newcommand{\dual}[1]{\langle {#1} \rangle}
\newcommand{\res}{\mathop{\hbox{\vrule height 7pt width.5pt depth 0pt\vrule height .5pt width 6pt depth 0pt}}\nolimits}
\newcommand{\V}{\mathbf{v}}

\newenvironment{proof}{\noindent\textbf{Proof}}{\hfill$\square$\medskip}
\newenvironment{cenni}{\noindent\textbf{Sketch of the proof}}{\hfill$\square$\medskip}
\newenvironment{prooftrm}{\noindent\textbf{Proof of the theorem}}{\hfill$\square$\medskip}
\newenvironment{oss}{\noindent\textbf{Remark}}{\hfill $\bullet$ \medskip}
\newenvironment{esempio}{\noindent\textbf{Example}}{\hfill $\bullet$ \medskip}

\DeclareMathOperator*{\osc}{osc}
\DeclareMathOperator{\diam}{diam}
\DeclareMathOperator{\diver}{div}
\DeclareMathOperator{\Lip}{Lip}
\DeclareMathOperator*{\dist}{dist}

\DeclareMathOperator{\spt}{spt}
\DeclareMathOperator{\diag}{diag}
\DeclareMathOperator{\loc}{loc}
\DeclareMathOperator{\Min}{Min}
\DeclareMathOperator{\rank}{rank}

\begin{document}

\begin{titlepage}
\begin{center}
\textsc{\Large University of Pisa\\}
\textsc{\large Faculty of Mathematical, Physical and Natural Sciences}\\
\vspace{0.5cm}

\textsc{\Large{Degree Course in Mathematics}}\\
\vspace{3.0cm}

\textsc{\Large Master's Thesis}\\
\Large{May 13, 2004}\\

\vspace{2.5cm}

\LARGE{\textbf{The non-parametric problem of Plateau\\in arbitrary codimension}}\\
\vspace{2.5cm}

\vspace{2.0cm}

\begin{minipage}{0.45\linewidth}
\centering
{\large Author}\\
\vspace{-0.2cm}
\textbf{{\large Luca M.A. Martinazzi}} \\
\vspace{-0.2cm}
{\small martinazzi@math.stanford.edu}\\
\end{minipage}
\hspace{.05\linewidth}
\begin{minipage}{.45\linewidth}
\centering
{\large Advisor}\\
\vspace{-0.2cm}
\textbf{{\large Prof. Mariano Giaquinta}}\\
\vspace{-0.2cm}
{\large Scuola Normale Superiore}
\end{minipage}
\end{center}

\vspace{0.7cm}
\begin{center}
\textsc{\Large Academic Year 2003/2004}
\end{center}
\newpage
\thispagestyle{empty}
\phantom{}
\end{titlepage}

\newpage

\tableofcontents

\newpage

\chapter*{Introduction}
\addcontentsline{toc}{chapter}{Introduction}

The object of this thesis is the study of the non-parametric Plateau problem: given a function $\psi:\partial\Omega\subset\R{n}\rightarrow\R{m}$, may we find a graph $\M{G}_u$, $u:\Cl\Omega\rightarrow\R{m},$ with boundary $\M{G}_{\psi}$ and whose area $\M{H}^n(\M{G}_u)$ is least among the submanifolds of $\R{n+m}$ having the same boundary? The problem is connected to the theory of PDE: if $u$ is a solution of the problem, then the first variation of the area of $\M{G}_u$ is zero, and this is equivalent to an elliptic equation, known as the \emph{minimal surface equation} \eqref{msediv}, if $m=1$, and an elliptic system, the \emph{minimal surface system} \eqref{mssnonpar}, if $m>1$. We say that $\M{G}_u$ is minimal if its first variation is zero.
We always assume the domain $\Omega$ and the boundary data $\psi$ to be $C^\infty$ and the functions $u$ considered to be at least Lipschitz.

In codimension 1 ($m=1$) the non-parametric Plateau problem has been widely studied until the late sixties. In 1968  H. Jenkins and J. Serrin show that the problem is solvable for arbitrary boundary data if and only if $\partial\Omega$ has everywhere non-negative mean curvature. This latter hypothesis gives an \textbf{a priori boundary gradient estimate}. The solution in codimension 1 is unique and, if $\Omega$ is convex, minimizes the area being \textbf{the area functional},  associating to a function $u$ the area of its graph $\M{A}(u)$, \textbf{strictly convex}. Moreover a Lipschitz solution of the minimal surface equation is $C^\infty$ thanks to the celebrated \textbf{theorem of De Giorgi} about the regularity of weak solutions of elliptic equations.

The methods used in codimension 1 don't apply to higher codimension: the a priori gradient estimates don't generalize, the area functional is no longer convex and the regularity theorem of the Giorgi holds only for scalar equations, not for systems.

In 1977 H. Lawson and R. Osserman prove that in codimension greater than 1 the problem of the existence of minimal graphs with prescribed boundary data isn't solvable in general even if the domain $\Omega$ is an $n$-dimensional ball. Also uniqueness and stability fail, due to the non-convexity of the area: they prove the existence of a boundary data $\psi$ for which the minimal surface system admits at least 3 solutions one of which unstable. Lawson and Osserman show, lastly, a Lipschitz but non-$C^1$ graph of least area, in contrast with the regularity theory in codimension 1.

In 2002 Mu-Tao Wang proves some positive results in arbitrary codimension. He shows that the mean curvature flow (the minus gradient flow of the area functional) of the initial graph $\M{G}_\psi$ (now $\psi$ is intended defined on all of $\Cl\Omega$) exists and converges to a minimal graph if the $C^2$ norm of $\psi$ is sufficiently small. The result is based on an a priori boundary gradient estimate and a recent theorem of Brian White giving local estimates for higher order derivatives.

Mu-Tao Wang also describes a region in the Grassmannian of $n$-planes $G(n,m)$ on which the logarithm of the inverse of the area functional is convex; this region contains the tangent planes of the \textbf{area-decreasing} graphs. Applying this result and a regularity theorem of Allard for minimal varifolds yields a Bernstein type theorem: the minimal graph of an area-decreasing function defined on all of $\R{n}$ is an $n$-dimensional plane. This theorem and a theorem of Allard imply that an area-decreasing minimal graph is $C^\infty$.

The exposition of the topics underlines the differences from the geometric and PDE point of view between the Plateau problem in codimension 1 and in higher codimension. The material of chapters \ref{capitolosottov}, \ref{capitolocodim1} and \ref{capitololo} is suitably covered by the literature of the last decades. The proofs of chapters \ref{capitolocodarb} and \ref{capitoloreg} are, on the contrary, very recent. The most important theorems are due to M-T Wang, while various propositions and explications have been added in orded to make the material easily understandable to an undergraduate major and to compare these results with the earlier approaches to the problem.

The ideas presented in this thesis could be further used: the convexity notions for the area among the area-decreasing maps could be useful to prove prove a uniqueness or stability theorem or in a variational approach. I had the possibility to discuss personally of these developments with prof. Mu-Tao Wang at the Columbia University, economically supported by the Scuola Normale Superiore and the research funds of prof. Wang; I very gladly thank them both. In several occasions I discussed the problems connected with my thesis with, apart from my advisor, prof. Luigi Ambrosio and prof. Giovanni Alberti, whom I thank for the interest and the advices.

I want, lastly, to thank sincerely my advisor, prof. Mariano Giaquinta, and the Scuola Normale Superiore. The former for the willingness and cordiality shown during this work, began in september 2002, when I asked him a topic for my third year colloquio at the Scuola Normale. The latter for providing me with a serene, stimulating and productive environment which, together with the University of Pisa, is a fertile breeding ground for a young student willing to enter the research world.

\chapter{Geometry of the submanifolds of $\R{n+m}$}\label{capitolosottov}

\section{Riemannian structures and Levi-Civita connections}

Given a Riemannian manifold $(M,g)$, a Levi-Civita connection on $M$ is an application
$$\nabla:\M{T}(M)\times \M{T}(M)\rightarrow\M{T}(M)$$ ($\M{T}(M)$ is the space of tangent vector fields on $M$) such that

\begin{enumerate}
\item $\nabla_X Y$ is $C^\infty$-linear in $X$:
$$\nabla_{fX_1+gX_2}Y=f\nabla_{X_1}Y+g\nabla_{X_2}Y,\quad\forall f,g\in C^\infty(M);$$
\item $\nabla_X Y$ is $\mathbb{R}$-linear in $Y$:
$$\nabla_X(aY_1+bY_2)=a\nabla_X Y_1+b\nabla_X Y_2,\quad\forall a,b,\in\mathbb{R};$$
\item it satisfies the Leibniz rule for the product:
$$\nabla_X(fY)=f\nabla_X(Y)+D_X fY, \quad\forall f\in C^\infty(M),$$
where $D_X f=X(f)$, being seen as a derivation;
\item it is torsion free: if $[X,Y]:=XY-YX$, then
$$\nabla_X Y-\nabla_YX=[X,Y];$$
\item it is compatible with the metric:
$$D_Xg(Y,Z)=g(\nabla_XY,Z)+g(Y,\nabla_X Z).$$
\end{enumerate}

\begin{trm}\label{unicalc} Every Riemannian manifolds has admits exactly one Levi-Civita connection.
\end{trm}
For a proof see \cite{lee}, theorem 5.4.

\medskip

In what follows we consider $\R{n+m}$ endowed with the usual Riemannian structure, in which the scalar product of two vectors $u,v\in\R{n+m}$ is denoted by $u\cdot v$ or $\dual{u,v}$. We identify $\R{n+m}$ with its tangent space in any of its points. $\R{n+m}$ has an unique Levi-Civita connection: it's the flat connection and we denote it by $\nabla$. Let $\{e_1,\ldots,e_{n+m}\}$ be an ortonormal basis of $\R{m+n}$, globally defined and fixed from now on; then
$$\Gamma_{ij}^k:=(\nabla_{e_i}e_j )^k=0, \quad\forall\, i,j,k.$$

\medskip
An $n$-dimensional submanifold $\Sigma\subset\R{n+m}$ of class $C^r$, $r\geq 2$, will be always endowed with the Riemannian structure provided by the ambient space: it's the only Riemannian structure such the immersion $$\Sigma\hookrightarrow \R{n+m}$$ is an isometry. Thus the metric $g$ on $\Sigma$ is simply the restriction of the metric of $\R{n+m}$. 

We denote by $T\Sigma$ its tangent bundle, of class $C^{r-1}$, and, for each $p\in\Sigma$,  $T_p \Sigma$ will be the tangent space to $\Sigma$ in $p$. Similarly $N\Sigma$ and $N_p\Sigma$ will denote the normal bundle and the normal space in $p$. An arbitrary orthonormal basis of $T_p\Sigma$ will be denoted by $\{\tau_1,\ldots\tau_n\}$ and an orthonormal basis of $N_p\Sigma$ by $\{\nu_1,\ldots,\nu_m\}$.

\medskip

The Levi-Civita connection of $\Sigma$ can be expressed in term of the flat connection $\nabla$ of $\R{n+m}$: $\nabla^\Sigma=\nabla^T$. More precisely, let $X,Y\in\M{T}(\Sigma)$ be tangent vector fields on $\Sigma$; given $\widetilde X$ and $\widetilde Y$, arbitrary extensions to a neighborhood of $\Sigma$ in $\R{n+m}$ of the fields $X$ and $Y$, we have 
\begin{equation}\label{levicivita}
\nabla^\Sigma_X Y=(\nabla_{\widetilde X} \widetilde Y)^T,
\end{equation}
where $(\nabla_{\widetilde X} \widetilde Y)^T$ is the orthogonal projection of $\nabla_{\widetilde X}\widetilde Y$ onto the tangent bundle $T\Sigma$. It may be verified that $\nabla^\Sigma$ doesn't depend on the choice of the extensions $\widetilde X$ and $\widetilde Y$. This is consequence of $\nabla_{\widetilde X}\widetilde Y(p)$ depending only on $X(p)$ and the value of $Y$ on the image of any curve $\gamma:(-\varepsilon,\varepsilon)\rightarrow \R{n+m}$ with $\gamma(0)=p$, $\dot\gamma(0)=X$. From now on, when necessary, the vector fields on $\Sigma$ will be intended as extended, at least locally.

To prove \eqref{levicivita}, we use theorem \ref{unicalc}, that is, thanks to the uniqueness of the Levi-Civita connection, it's enough to prove that $(X,Y)\rightarrow(\nabla_X Y)^T$ is a Levi-Civita connection. The $C^\infty$-linearity in $X$ and the $\mathbb{R}$-linearity in $Y$ are trivial, as well as the Leibniz rule. Let's show that there is no torsion (property 4 in the definition):
$$(\nabla_X Y)^T-(\nabla_Y  X)^T=
(\nabla_X Y- \nabla_Y X)^T =[X,Y]^T=[X,Y].$$

Let's verify the compatibility with the metric:
$$D_X g(Y,Z)=g(\nabla_X Y,Z)+g(X,\nabla_Y Z)=g((\nabla_X Y)^T ,Z)+g(X,(\nabla_Y Z)^T).$$

\subsection{The gradient, divergence and Laplacean operators}

Given a $C^1$ function $f:\Sigma\rightarrow\mathbb{R}$ and $X\in T_p\Sigma$, we define
$$D_X f(p)=\frac{d}{dt}\Big{|}_{t=0}f(\gamma(t)),$$
for any curve $\gamma:(-\varepsilon,\varepsilon)\rightarrow \Sigma$ such that $\gamma(0)=p$ and $\dot\gamma(0)=X$.

The \emph{gradient} on $\Sigma$ of $f$ in $p$ is defined by
$$\nabla^\Sigma f(p)=\sum_{j=1}^n (D_{\tau_j}f(p))\tau_j.$$

It's not hard to prove that if $f$ is defined in a neighborhood of $p$ in $\R{n+m}$, then we have
$$\nabla^\Sigma f(p)=(\nabla f(p))^T,$$
where $\nabla f(p)=\sum_{j=1}^{n+m}\frac{\partial f}{\partial x^j}(p)e_j$. 

In a local frame, that is given a chart $(V,\varphi),$ with $\varphi:V\rightarrow\R{n}$, and given the corresponding local parametrization $F=\varphi^{-1}$ the following holds:
\begin{equation}\label{locgrad}
\nabla^{\Sigma}f=g^{ij}\frac{\partial f}{\partial x^i}\frac{\partial F}{\partial x^j},
\end{equation}
where $\frac{\partial}{\partial x^i} f(p)=\frac{\partial (f\circ\varphi^{-1})}{\partial x^i}(\varphi(p)),$
$g_{ij}=\frac{\partial F}{\partial x^i}\cdot\frac{\partial F}{\partial x^j}$ and
$(g^{ij})$ is the invers matrix of $(g_{ij})$.

\medskip

The \emph{divergence} of a vector field (not necessarily tangent) $\sum_{j=1}^{n+m} X^j e_j$ on $\Sigma$ is defined by

$$\diver^\Sigma X=\sum_{j=1}^{n+m}e_j\cdot(\nabla^\Sigma X^j)=\sum_{i=1}^n(D_{\tau_i}X)\cdot \tau_i.$$ 

In local coordinates, with the same notation as in \eqref{locgrad} and writing $g=\det (g_{ij})$

\begin{equation}\label{locdiv}
\diver^\Sigma X=\frac{1}{\sqrt{g}}\frac{\partial}{\partial x^i}(\sqrt{g}X^i).
\end{equation}

Finally the \emph{Laplacean} on $\Sigma$ of a function in $C^2(\Sigma)$ is defined as
$$\Delta_\Sigma f=\diver^\Sigma\nabla^\Sigma f,$$
which may be written in local coordinates plugging \eqref{locgrad} into \eqref{locdiv}:
\begin{equation}\label{loclapl}
\Delta_\Sigma f=\frac{1}{\sqrt{g}}\frac{\partial}{\partial x^i}\Big(\sqrt{g}g^{ij}\frac{\partial f}{\partial x^j}\Big).
\end{equation}

\section{The second fundamental form and the mean curvature}

\begin{defin}[Second fundamental form]
We define the second fundamental form $h$ to be the normal part of the connection of
$\R{n+m}:$ given $X,Y\in\M{T}(M)$
$$h(X,Y)=(\nabla_X Y)^N.$$
\end{defin}
As before, $X$ and $Y$ are extended.
\begin{prop} The second fundamental form $h$:
\begin{enumerate}
\item is symmetric: $h(X,Y)=h(Y,X);$
\item is $C^\infty$-linear in both variables;
\item $h(X,Y)(p)$ depends only on $X(p)$ and $Y(p).$
\end{enumerate}
In particular $h$ is well defined as a family of bilinear applications
$$h_p:T_p\Sigma\times T_p\Sigma\rightarrow N_p\Sigma.$$
\end{prop}

\begin{proof} Due to the symmetry of $\nabla$ and since for $X,Y\in\M{T}(\Sigma)$ we have $[X,Y]\in\M{T}(\Sigma),$ the following holds true
$$h(X,Y)-h(Y,X)=(\nabla_XY-\nabla_YX)^N=[X,Y]^N=0.$$
To prove 2, we observe that $h$ is the difference of two connections:
$$h(X,Y)=\nabla_X Y-\nabla^\Sigma_X Y$$
it is, thus, $C^\infty$-linear in $X$. Being $h$ symmetric it is also $C^\infty$-linear in $Y$.

Finally, both $\nabla_XY(p)$ and $\nabla_X^\Sigma Y(p)$ depend only on $Y$ and $X(p)$. By symmetry, it's enough to know $Y(p)$ and $X$ and, thus, it's also enough to know only $X(p)$ and $Y(p)$.
\end{proof}

\begin{defin}[Mean curvature]
For each $p\in\Sigma,$ we define the mean curvature $H$ of $\Sigma$ in $p$ to be the trace of the second fundamental form, that is
$$H(p)=\sum_{i=1}^{n}h_p(\tau_i,\tau_i).$$
If $\{v_1,\ldots,v_n\}$ is an arbitrary basis of $T_p\Sigma$ and $g_{ij}:=g(v_i,v_j),$ then
\begin{equation}\label{curvmedia}
H(p)=\sum_{i,j=1}^n g^{ij}h_p(v_i,v_j).
\end{equation}
\end{defin}

We use \eqref{curvmedia} to compute the mean curvature of $\Sigma$: let a local pa\-ra\-me\-tri\-za\-tion $F:\Omega\rightarrow\Sigma$ be given in $p$, that is a diffeomorphism of $\Omega$ with a neighbourhood of $p$. We assume $F(0)=p$. $F$ induces a basis of $T_p\Sigma$, given by $\bbra{\frac{\partial F}{\partial x^i}}_{i=1,\ldots,n}$.
$$\nabla_{\frac{\partial F}{\partial x^i}}\frac{\partial F}{\partial x^j}=\frac{\partial^2 F}{\partial x^i \partial x^j}.$$
Using \eqref{curvmedia} yields
\begin{equation}\label{curvmedia2}
H(p)=\Big( \sum_{i,j=1}^n g^{ij}\frac{\partial^2 F}{\partial x^i \partial x^j}(F^{-1}(p)) \Big)^N,
\end{equation}
$g_{ij}=\frac{\partial F}{\partial x^i}\cdot\frac{\partial F}{\partial x^j}.$

\begin{lemma}[Derivative of a determinant]
Let $g(s)=\det (g_{ij}(s))$, $g_{ij}$ being differentiable in $s$. Then
\begin{equation}\label{derivadet}
\frac{\partial g}{\partial s}=gg^{ij}\frac{\partial g_{ij}}{\partial s}
\end{equation}
\end{lemma}

\begin{prop} Let $F:\Omega\rightarrow\Sigma$ be a local parametrization in $p$. Then $\Delta_\Sigma F(p)\in N_p\Sigma$ and
\begin{equation}\label{curvmedialapl}
H(p)=\Delta_\Sigma F(p).
\end{equation}
The Laplacean of $F$ is defined componentwise.
\end{prop}

\begin{proof}
We prove that $\Delta_\Sigma F(p)$ is orthogonal to $T_p\Sigma$. From now on we suppress $p$. Thanks to \eqref{loclapl} we may write
$$\Delta_\Sigma F\cdot\frac{\partial F}{\partial x^k}=\frac{1}{\sqrt{g}}\frac{\partial}{\partial x^i}\big(
\sqrt{g}g^{ij} \frac{\partial F}{\partial x^j}\cdot \frac{\partial F}{\partial x^k}\big)
-g^{ij}\frac{\partial F}{\partial x^j}\cdot\frac{\partial^2 F}{\partial x^i \partial x^k}.$$
We observe that
$$g^{ij}\frac{\partial F}{\partial x^j}\cdot\frac{\partial^2 F}{\partial x^i \partial x^k}=
\frac{1}{\sqrt{g}}\frac{\partial\sqrt{g}}{\partial x^k}-\frac{1}{2}g^{ij}\frac{\partial g_{ij}}{\partial x^k}.$$
Subsituting the formula for the derivative of a determinant \eqref{derivadet} we obtain
$$\Delta_\Sigma F\cdot\frac{\partial F}{\partial x^k}=0.$$
Since $k$ is arbitrary we conclude that $\Delta_\Sigma F$ is orthogonal to $\Sigma$.

Let us prove \eqref{curvmedialapl}: writing the Laplacean in a local frame and differentiating yields
$$\Delta_\Sigma F= \frac{1}{\sqrt{g}}\frac{\partial}{\partial x^i}\Big(\sqrt{g}g^{ij}\frac{\partial F}{\partial x^j}\Big)=
\frac{1}{\sqrt{g}}\frac{\partial}{\partial x^j}\Big(\sqrt{g}g^{ij}\Big)\partial_j F+g^{ij}\frac{\partial^2F}
{\partial x^i\partial x^j}.$$
Observing that the first term in the right hand side is tangent and using \eqref{curvmedia2} give
$$\Delta_\Sigma F=(\Delta_\Sigma F)^N= \Big(g^{ij}\frac{\partial^2F}{\partial x^i\partial x^j}\Big)^N=H.$$
\end{proof}

\section{The area formula: first variation}\label{sezformula}

We shall call area of $\Sigma$ the $n$-dimensional Hausdorff measure of $\Sigma$, i.e.
$\M{A}(\Sigma):=\M{H}^n(\Sigma)$. It may be computed by means of the area formula.

\begin{trm}[Area formula]
Let $F:\Omega\rightarrow\R{n+m}$ be a \emph{locally Lipschitz} and injective map of an open set $\Omega\subset\R{n}$ into $\R{n+m}$. Let $\Sigma$ be the image of $F$; then
\begin{equation}\label{formulaarea}
\M{H}^n(\Sigma)=\int_\Omega\sqrt{\det{dF^*(x)dF(x)}}dx,
\end{equation}
where $dF^*:\R{n+m}\rightarrow\R{n}$ is the transposed of $dF$.
\end{trm}

For a proof of this theorem see \cite{federer} or \cite{giaquinta}.

If $g_{ij}=\frac{\partial F}{\partial x^i}\cdot\frac{\partial F}{\partial x^j}$, we observe that
$(dF^*dF)_{ij}=\sum_{\alpha=1}^{n+m}\frac{\partial F^\alpha}{\partial x^i}\frac{\partial F^\alpha}{\partial x^j}=
g_{ij}$ thus, being $g=\det g_{ij}$,
\begin{equation}\label{areag}
\M{A}(\Sigma)=\int_\Omega\sqrt{g(x)}dx.
\end{equation}
In particular $\sqrt{g}dx$ is the area element of $\Sigma$ expressed through the pa\-ra\-me\-tri\-za\-tion $F$ so that, given an $\M{H}^n\res\Sigma$-integrable function $f$, we have
$$\int_\Sigma f d\M{H}^n=\int_\Omega f\circ F\sqrt{g}dx.$$

\paragraph{First variation of the area}

\begin{defin}\label{defvar}
Given $\Sigma\subset\R{n+m}$ at least $C^1$, we consider a family of diffeomorphisms $\varphi_t:\R{n+m}\rightarrow\R{n+m}$ such that
\begin{enumerate}
\item $\varphi(t,x):=\varphi_t(x)$ is $C^2$ in $(-1,1)\times \R{n+m};$
\item there exists a compact $K$ non intersecting $\partial\Sigma$ (possibly empty) such that
$\varphi_t(x)=x$ for each $x\notin K$ and $t\in(-1,1)$;
\item $\varphi_0(x)=x$ for each $x\in\R{n+m}$.
\end{enumerate}
\end{defin}

\begin{prop}\label{propvarprima} Set $\Sigma_t=\varphi_t(\Sigma)$ and $X=\frac{\partial \varphi_t}{\partial t}\big|_{t=0}$.  
Assume $\Sigma$ to be at least $C^1$ and admitting a global parametrization $F:\Omega\rightarrow\R{n+m}$.
Then
\begin{equation}\label{varprima}
\frac{d}{dt}\M{A}(\Sigma_t)\Big{|}_{t=0}=-\int_\Sigma\Delta_\Sigma F\cdot X d\M{H}^n,
\end{equation}
where the Laplacien has to be read in the weak sense, that is

\begin{multline}
-\int_\Sigma\Delta_\Sigma F\cdot X =-\int_\Omega 
\frac{1}{\sqrt{g}}\frac{\partial}{\partial x^j}\Big(\sqrt{g}g^{ij}\frac{\partial F^\alpha}{\partial x^j} \Big) X^\alpha\sqrt{g}dx=\\
=\int_\Omega \sqrt{g}g^{ij}\frac{\partial F^\alpha}{\partial x^j}\frac{\partial X^\alpha}{\partial x^i}
dx=\int_\Sigma g^{ij}\frac{\partial F^\alpha}{\partial x^j}\frac{\partial X^\alpha}{\partial x^i}d\M{H}^n.
\end{multline}
\end{prop}
Here and subsequently, integration by parts, even if only formal, is justified by being $\varphi_t=Id$ outside a compact $K$ non intersecting $\partial\Sigma$.

\medskip

\begin{proof}
Being $\varphi$ differentiable and $\varphi_0(y)=y$, we have
\begin{equation}\label{taylor}
\varphi_t(y)=y+tX(y)+o(t), \quad X(y):=\frac{\partial\varphi_t}{\partial t}\Big{|}_{t=0}.
\end{equation}
We differentiate the area formula \eqref{areag} under the integral sign and use the formula for the derivative of a determinant \eqref{derivadet}: set $F_t(x)=\varphi_t(F(x))$ and $g^t_{ij}=\frac{\partial F_t}{\partial x^i}
\cdot\frac{\partial F_t}{\partial x^i}.$ All the derivatives with respect to $t$ are computed for $t=0$ and clearly $g=g^0$.
\begin{equation}\label{contovarprima}
\frac{d}{dt}\int_\Omega\sqrt{g^t}dx=\int_\Omega \frac{\partial \sqrt{g^t}}{\partial t}dx=
\int_\Omega\frac{1}{2\sqrt{g}}\Big(gg^{ij}\frac{\partial g^t_{ij}}{\partial t}\Big)dx.
\end{equation}
To compute $\frac{\partial g^t_{ij}}{\partial t}$ we observe that, thanks to \eqref{taylor}, $\frac{\partial \varphi^t(F(x))}{\partial x^i}
=\frac{\partial F}{\partial x^i}+t\frac{\partial X}{\partial x^i}+o(t)$ and substituting into \eqref{contovarprima} yields
\begin{multline}\label{contivarprima}
\int_\Omega\frac{1}{2\sqrt{g}}\Big(gg^{ij}\frac{\partial g^t_{ij}}{\partial t}\Big)dx=\\
=\int_\Omega\frac{1}{2}\sqrt{g}g^{ij}\frac{\partial}{\partial t}\Big(\big(\frac{\partial F}{\partial x^i}+t\frac{\partial X}{\partial x^i}+o(t) \big)\cdot\big(\frac{\partial F}{\partial x^j}+t\frac{\partial X}{\partial x^j}+o(t) \big)\Big)dx=\\
\int_\Omega\frac{1}{2}\sqrt{g}g^{ij}\Big(\frac{\partial X(F(x))}{\partial x^i}\cdot\frac{\partial F}{\partial x^j}
+\frac{\partial F}{\partial x^i}\cdot\frac{\partial (X(F(x))}{\partial x^j} \Big)dx.
\end{multline}
Due to the symmetry of $g^{ij}$ the last term becomes
\begin{equation}\label{finalevarprima}
\int_\Omega \sqrt{g}g^{ij}\frac{\partial X(F(x))}{\partial x^i}\cdot\frac{\partial F}{\partial x^j}dx
=-\int_\Sigma\Delta_\Sigma F\cdot X.
\end{equation}
\end{proof}

\begin{prop}\label{vardiv}
Let $\Sigma$, $\varphi_t$ and $X$ be as in proposition \ref{propvarprima}. Then
$$\frac{d}{dt}\M{A}(\Sigma_t)\Big{|}_{t=0}=\int_\Sigma \diver^\Sigma X.$$
\end{prop}

\begin{oss}
Differently from proposition \ref{propvarprima}, this proposition doesn't require the existence of a global parametrization, thus it may be considered more intrinsic.
\end{oss}

\begin{proof} Let an arbitrary basis of $T_p\Sigma$ be given, say $\{v_1,\ldots,v_n\}$, and set $g_{ij}=v_i\cdot v_j$. Then, by linearity
$$\diver^\Sigma X=g^{ij}\nabla_{v_i}X\cdot v_j.$$
Consequently, choosing a local parametrization $F$ in $p$, setting $v_i:=\frac{\partial F}{\partial x^i}$ and using $\nabla_{\frac{\partial F}{\partial x^i}}X
=\frac{\partial X(F)}{\partial x^i}$ we obtain
$$\diver^\Sigma X =g^{ij}\frac{\partial X(F(x))}{\partial x^i}\cdot\frac{\partial F}{\partial x^j}.$$
We conclude by comparison with \eqref{contivarprima}.
\end{proof}

\begin{oss} Propositions \ref{propvarprima} and \ref{vardiv} characterize the first variation of the area of a submanifold in the only hypothesis that the submanifold is $C^1$. Actually less is needed: both propositions may be repeated \emph{verbatim} for Lipschitz submanifolds using Rademacher's theorem, see the appendix.
\end{oss}

Now we see how the mean curvature gets involved in the definition of minimal surface and first variation.

\begin{prop}\label{propvarmedia} Let $\Sigma$ be a $C^2$ submanifold and let be given a variation $\varphi_t$ with variation field $X$. Then the first variation of the area of $\Sigma$ with respect to $\varphi$ is
\begin{equation}\label{varcurvmedia}
\frac{d}{dt}\M{A}(\Sigma_t)\Big{|}_{t=0}=-\int_\Sigma H\cdot X.
\end{equation}
\end{prop}

\begin{proof}
Plug \eqref{curvmedialapl}, holding for $C^2$ submanifolds, into proposition \ref{propvarprima}.
\end{proof}

\section{Minimal surfaces}

By minimal surface we mean a submanifold $\Sigma$ whose area is stationary with respect to compactly supported variations keeping its boundary fixed:
\begin{defin}[Minimal surface]
Let $\Sigma$ be a Lipschitz $n$-submanifold of $\R{n+m}$. We shall say that $\Sigma$ is minimal if for every variation $\varphi_t$, definition \ref{defvar}, we have
$$\frac{d}{dt}\Big{|}_{t=0}\M{A}(\Sigma_t)=0.$$
\end{defin}

Thanks to propositions \ref{propvarprima}, \ref{vardiv} and \ref{propvarmedia}, we have the following proposition characterizing minimal surfaces.

\begin{prop}\label{caratmin} Given a Lipschitz submanifold $\Sigma$ of $\R{n+m}$, the following are equivalent:
\begin{enumerate}
\item $\Sigma$ is minimal;
\item for every vector field $X\in C^1_0(\R{n+m};\R{n+m})$ such that $X=0$ in a neighbourhood of $\partial\Omega$
$$\int_\Sigma \diver^\Sigma X=0;$$
\item for each local parametrization $F:\Omega\rightarrow \Sigma$ we have $\Delta_\Sigma F=0$ weakly.
\end{enumerate}
Moreover, if $\Sigma\in C^2$, the preceding statements are equivalent to $H=0$.
\end{prop}

\begin{proof} We have proved that $2\Rightarrow 1$ and $3\Rightarrow 2$. $2\Rightarrow 3$ is also true because
$$0=\int_\Sigma \diver^\Sigma X=-\int_\Sigma \Delta_\Sigma F\cdot X.$$
Being $X$ arbitrary we conclude that $\Delta_\Sigma F=0$.

In order to prove that $1\Rightarrow 2$, it's enough to prove that for every vector field $X\in C^1_0(\R{n+m};\R{n+m})$ vanishing in a neighborhood of $\partial\Sigma$ we may find a family of diffeomorphisms $\varphi_t$ as in definition \ref{defvar} satisfying $\frac{\partial \varphi_t}{\partial t}\big{|}_{t=0}=X$. This may be easily obtained locally: we define a family of variations $\varphi_t^{(i)}$ which may be glued together by means of a partition of unity.

The last claim is an immediate consequence of proposition \ref{propvarmedia}.
\end{proof}

\subsection{The minimal surface system}

Consider a parametrizatione $F:\Omega\rightarrow\R{n+m}$ of a Lipschitz submanifold $\Sigma\subset\R{n+m}$. Thanks to proposition \ref{caratmin}, $\Sigma$ is miniml if and only if $F$ satisfies the following system, called \emph{minimal surface system}:

\begin{equation}\label{ssm}
\sum_{i,j=1}^n\frac{\partial}{\partial x^i}\Big(\sqrt{g}g^{ij}\frac{\partial F^\alpha}{\partial x^j}\Big)=0,\quad \alpha=1,\ldots,n+m,
\end{equation}
where $g=\det (g_{ij})$, $g_{ij}=\frac{\partial F}{\partial x^i}\cdot\frac{\partial F}{\partial x^j}$ and
$(g^{ij})=(g_{ij})^{-1}$.

The definition is well-posed and is intended in the weak sense, i.e., for each $\varphi\in C^\infty_0(\Omega)$
$$\sum_{i,j=1}^n\int_\Omega \sqrt{g}g^{ij}\frac{\partial F^\alpha}{\partial x^j}\frac{\partial \varphi}{\partial x^i}=0.$$

\subsection{Non parametric minimal surfaces} A non parametric surface $\Sigma$ is the graph $\M{G}_u$ of a Lipschitz function $u:\Omega\rightarrow\R{m}$. $\M{G}_u$ is clearly parametrized by the immersion
$$F:=I\times u:\Omega\rightarrow \R{n+m},$$
that is $F(x)=(x,u(x)).$ In this case, the minimal surface system becomes

\begin{equation}\label{mssnonpar}
\left\{ \begin{array}{ll}
\displaystyle
\sum_{i=1}^{n}\frac{\partial}{\partial x^i} \Big(\sqrt{g}g^{ij} \Big)=0 &j=1,\ldots,n\\
\displaystyle
\sum_{i,j=1}^n\frac{\partial}{\partial x^i}\Big(\sqrt{g}g^{ij}\frac{\partial u^\alpha}{\partial x^j}\Big)=0&\alpha=1,\ldots,m.
\end{array}
\right.
\end{equation}
Also the equazions of this system are to be read in the weak sense. It's an elliptic system in divergence form.

Actually, at least in the case of the graph of a $C^2$ function, the system \eqref{mssnonpar} reduces to a quasilinear elliptic system in nondivergence form, as the following proposition shows.

\begin{prop}\label{equivale} Let be $u\in C^2(\Omega)$. Then the system \eqref{mssnonpar} is equivalent to
\begin{equation}\label{eqequivale}
\sum_{i,j=1}^n g^{ij}\frac{\partial^2 u^\alpha}{\partial x^i\partial x^j}=0,\quad\alpha=1,\ldots,m.
\end{equation}
\end{prop}

\begin{proof} Let \eqref{eqequivale} hold true and set $F(x)=(x,u(x))$.
$$\Delta_\Sigma F=\frac{1}{\sqrt{g}}\frac{\partial}{\partial x^i}(\sqrt{g}g^{ij})\frac{\partial F}{\partial x^j}+
g^{ij}\frac{\partial^2 F}{\partial x^i x^j}.$$
The last term vanishes because $g^{ij}\frac{\partial^2 x^k}{\partial x^i x^j}=0$ for every $k=1,\ldots,n$ and
$g^{ij}\frac{\partial^2 u^\alpha}{\partial x^i x^j}=0$ by hypothesis. Since $\Delta_\Sigma F\in N\Sigma$ and $\frac{1}{\sqrt{g}}\frac{\partial}{\partial x^i}(\sqrt{g}g^{ij})\frac{\partial F}{\partial x^j}$ is tangent, we conclude that it has to vanish as well and, thus, $\Delta_\Sigma F=0$.

Conversly, if \eqref{mssnonpar} holds true and $u$ is $C^2$, then, thanks to proposition \ref{caratmin}, $H=0$. We conclude using \eqref{curvmedia2}.
\end{proof}

\section{Singular values: the area-decreasing maps}

From the area formula \eqref{formulaarea}, we know that the area of the graph of a Lipschitz function $u:\Omega\rightarrow\R{m}$ is
\begin{equation}\label{eqareagraf}
\M{A}(\M{G}_u)=\int_\Omega\sqrt{\det(DF^*DF)}dx=\int_\Omega \sqrt{\det\big(I+Du^*Du\big)}dx;
\end{equation}
there always exists a local frame in which $Du$ and the area element have a particularly simple form. 
\begin{prop}[Singular value decomposition]
Let be given $A$, an $m\times n$ matrix. Then there exist $U$ and $V$ orthogonal matrices on $\R{m}$ and $\R{n}$ respectively, such that
$B=UAV$ is a diagonal matrix: if $B=\{\lambda_{\alpha i}\}_{\alpha=1,\ldots,m}^{i=1,\ldots,n}$, then
$\lambda_{\alpha i}=0$ whenever $\alpha\neq i$.
\end{prop}

\begin{proof}
For a proof see \cite{leon}, theorem 7.7.1.
\end{proof}

\begin{oss} We may and do assume that $\lambda_{\alpha i}\geq 0$: indeed changing the sign of the basis vectors is an orthogonal transformation. It's obvious that $\max \lambda_i=|Du|.$
\end{oss}

An application of the singular value decomposition to the differential $Du$ yields
$Du^*Du=\diag \{\lambda_1^2,\ldots,\lambda_n^2\}$, where $\lambda_i:=\lambda_{ii}$ if $i\leq m$ and $\lambda_i=0$ otherwise.
Thus
\begin{equation}\label{areasing}
\M{A}(\M{G}_u)=\int_\Omega \sqrt{\textstyle\prod_{i=1}^n(1+\lambda_i^2(x))}dx.
\end{equation}
\begin{defin}
Let $u:\Omega\rightarrow\R{m}$ be a Lipschitz map. Let $\{\lambda_i(x)\}_{i=1,\ldots, n}$ be the singular values of $Du(x)$. We shall say that $u$ is area-decreasing if there exists $\varepsilon>0$ such that for a.e. $x\in\Omega$ we have
\begin{equation}\label{adcond}
\lambda_i(x)\lambda_j(x)\leq 1-\varepsilon,\quad 1\leq i<j\leq n.
\end{equation}
\end{defin}

The geometric meaning of the \emph{area-decreasing} condition is the following: consider $Du(x)$ restricted to a 2 dimensional subspace $V$ of $\R{n}$. Then for each $A\subset V$ with $\M{H}^2(A)<\infty$ we have $\M{H}^2(Du(x)(A))<\M{H}^2(A).$ Equivalently the Jacobian of $Du(x)\big|_V$ is less than 1.

\medskip

\begin{oss} If $m=1$, that is $u:\Omega\rightarrow\mathbb{R}$, then $u$ is area-decreasing. This follows immediatly from the definition because the nonzero singular values of $Du(x)$ correspond to a basis of the image of $Du(x)$ and therefore  in \eqref{adcond} $\lambda_i$ and $\lambda_j$ cannot be both nonzero.
Actually it's natural for a scalar function to be area-decreasing, since $\M{H}^2(\mathbb{R})=0.$
\end{oss} 

\medskip

As we shall see, in the area-decreasing category it's possible to prove existence, regularity and rigidity theorems for minimal graphs in arbitrary codimension which are natural generalizations of the corresponding theorems in codimension 1.

\section{Statement of the non parametric problem of Plateau}

The non parametric, or Cartesian, Plateau problem requires to find graphs of least area with prescribed boundary. The boundary $\Gamma$ is given as the graph of a smooth given map $$\psi:\partial\Omega\rightarrow\R{m},\quad \Gamma=\M{G}_\psi$$ where $\Omega$ is $C^\infty$ domain in $\R{n}$.\footnote{that is, every point $x\in\partial\Omega$ has a neighborhood diffeomorphic to a half space of dimension $n$.} Consider the set of Lipschitz $n$-submanifolds
$$A=\big\{\Sigma\subset\R{n+m} \; : \; \partial\Sigma=\Gamma\}.$$ 
Such a set is nonempty being the homology of $\R{n+m}$ trivial.

We shall discuss the following problems:

\medskip
\noindent\textbf{Problem 1: Existence of minimizers}\\ Is it possible to find an application $u\in \Lip(\Cl\Omega;\R{m})$ such that $\M{G}_u\in A$ and
$$\M{A}(\M{G}_u)\leq\M{A}(\Sigma),\quad \forall\Sigma\in A?$$ 

\medskip

Weakening the problem.

\medskip

\noindent\textbf{Problem 2: Existence of critical points}\\ Does it exist $u\in\Lip(\Cl\Omega;\R{m})$ such that $\M{G}_u\in A$ and whose graph has vanishing first variation?

\medskip

\noindent\textbf{Problem 3: Stability}\\ Does a solution of problem 2 also solve problem 1? Is it at least stable, that is, small variations don't diminish the area?

\medskip

\noindent\textbf{Problem 4: Uniqueness}\\ Is a solution of problem 1 or problem 2 unique?

\medskip

\noindent\textbf{Problem 5: Regularity}\\ Is a solution of problem 1 or problem 2 regular?

\medskip

We shall not discuss the analogous problems arising when searching for minimizing or stationary graphs in the class

$$B=\{\M{G}_v\;:\;v\in\Lip(\Omega;\R{m}),\;\partial\M{G}_u=\Gamma \}.$$

Nontheless we shall see that in codimension 1 the two classes of problems are often very close together.

\subsection{The Dirichlet problem}
For every map $\psi:\Cl\Omega\rightarrow\R{m}$, we shall call Dirichlet problem for the minimal surface system the following system:
\begin{equation}\label{dirnonpar}
\left\{ \begin{array}{ll}
\displaystyle
\sum_{i=1}^{n}\frac{\partial}{\partial x^i} \Big(\sqrt{g}g^{ij} \Big)=0 &j=1,\ldots,n\\
\displaystyle
\sum_{i,j=1}^n\frac{\partial}{\partial x^i}\Big(\sqrt{g}g^{ij}\frac{\partial u^\alpha}{\partial x^j}\Big)=0&\alpha=1,\ldots,m.\rule{0cm}{.7cm}\\
\displaystyle
u^\alpha\big|_{\partial\Omega}=\psi^\alpha\big|_{\partial\Omega}   &\alpha=1,\ldots,m.\rule{0cm}{.7cm}
\end{array}
\right.
\end{equation}

The Dirichlet problem is equivalent to problem 2 thanks to proposition \ref{caratmin} and the solutions of problem 1 also solve the Dirichlet problem because the first variation of a minimizing surface vanishes.

The regularity problem is intimately connected to the nature of the minimal surface system and, in codimension 1, to the minimal surface equation \eqref{msediv}.

We observe that if $u\in C^2(\Omega)\cap C^0(\Cl\Omega)$, due to proposition \ref{equivale}, the system \eqref{dirnonpar} is equivalent to 

\begin{equation}\label{dirnondiv}
\left\{ \begin{array}{ll}
\displaystyle
\sum_{i,j=1}^{n}g^{ij}\frac{\partial^2 u^\alpha}{\partial x^i\partial x^j}=0 &\alpha=1,\ldots,m\\
\displaystyle
u^\alpha\big|_{\partial\Omega}=\psi^\alpha\big|_{\partial\Omega}   &\alpha=1,\ldots,m.\rule{0cm}{.7cm}
\end{array}
\right.
\end{equation}

\chapter{Codimension 1}\label{capitolocodim1}

\section{Convexity of the area}

Let $\Omega$ be a convex smooth domain in $\R{n}$.
In codimension 1 the area functional, defined on the space $\Lip(\Cl\Omega)$, may be easily rewritten as
\begin{equation}\label{areacodim1}
\M{A}(u)=\int_\Omega\sqrt{1+|Du|^2}dx.
\end{equation}

\begin{prop}[Convexity]\label{propconv}
The area functional
$$\M{A}:\Lip(\Cl\Omega)\rightarrow\mathbb{R} $$
in codimension 1 is stricly convex, that is
$$\M{A}(\lambda u+(1-\lambda) v)\leq\lambda\M{A}(u)+(1-\lambda)\M{A}(v),$$
for every $u,v\in\Lip(\Cl\Omega)$ and $\lambda\in(0,1)$ and equality holds if and only if
$u=v+c$ for some $c\in\R{}$.
\end{prop}

\begin{proof}
Observe that $f(x)=\sqrt{1+x^2}$ is a stricly convex function, being its second derivative
$$f''(x)=\frac{1}{(1+x^2)^\frac{3}{2}}>0.$$
Then the area functional is composition of a linear map $(u\rightarrow Du)$, a convex function $(p\rightarrow |p|)$, another convex function $(x\rightarrow \sqrt{1+x^2})$ and a linear functional (the integral on $\Omega$). Since composition of convex functions is convex, we have the convexity of the area. To verify that this convexity is strict, let $u,v$ be such that $u\neq v+c$. Then
\begin{multline}
\int_\Omega\sqrt{1+|D(\lambda u+(1-\lambda )v|^2}dx\leq
\int_\Omega \sqrt{1+(\lambda|Du|+(1-\lambda) |Dv|)^2}<\\<\lambda\M{A}(u)+(1-\lambda)\M{A}(v).
\end{multline}
The last inequality is strict because $Du\neq Dv$.
\end{proof}

\begin{oss} Convexity is one of the most important properties of the area functional in codimension 1 and it marks a  decided difference between the Plateau problem in codimension 1 and higher. Indeed, as we shall see, existence and uniqueness of graphs of least area in codimension 1 is linked to convexity. Such results are false in higher codimension as we shall see in the counterexamples of Lawson and Osserman.
\end{oss}

\section{Uniqueness and stability}

\begin{trm}\label{mingraf} In codimension 1 the graph of a Lipschitz solution $u:\Cl\Omega\rightarrow\mathbb{R}$ of the minimal surface system \eqref{mssnonpar} minimizes the area among the graphs of Lipschitz functions $v$ such that $u=v$ on $\partial\Omega$.
Moreover $u$ satisfies the \emph{minimal surface equation} in divergence form which is equivalent, in codimension 1, to the minimal surface system:
\begin{equation}\label{msediv}
D_i\frac{D_i u}{\sqrt{1+|Du|^2}}=0.
\end{equation}
This solution is unique.
\end{trm}

Equation \eqref{msediv} is meant to be read in the weak sense.

\begin{proof} \textbf{1.} The minimal surface system implies that the first variation of the area of the graph $\M{G}_u$ vanishes. In particular, for a given function $\varphi\in C^1_c(\Omega)$ we have
\begin{equation}\label{equnic}
0=\frac{d}{dt}\Big|_{t=0}\M{A}(u+t\varphi)=\int_\Omega \frac{\partial}{\partial t}
\sqrt{1+|Du+tD\varphi|^2}dx=-\int_\Omega \frac{D_iuD_i\varphi}{\sqrt{1+|Du|^2}}dx,
\end{equation}
which is the minimal surface equation in divergence form \eqref{msediv}.

\textbf{2.} Equation \eqref{equnic} says that $u$ is  a {critical point} for the area functional. On the other hand convexity implies
$$\M{A}(v)\geq\M{A}(u)+\frac{d}{dt}\Big|_{t=0}\M{A}(u+tv)=\M{A}(u).$$

\textbf{3.} Uniqueness follows by strict convexity of $\M{A}$, which implies that, given $u\neq v$ solutions to the minimal surface equation, we have
$$\M{A}\bra{\frac{u+v}{2}}<\frac{1}{2}(\M{A}(u)+\M{A}(v))=\M{A}(u).$$
The equality follows from $u$ and $v$ being minimizers and this contradicts the inequality.
\end{proof}

\begin{oss} In the class of Lipschitz functions, the minimal surface equation in divergence form \eqref{msediv} is equivalent to the minimal surface system \eqref{mssnonpar}. This means that in order to verify the vanishing of the first variation, it's enough to consider deformations of the form $u+t\varphi$, called non-parametric deformations. We may reduce ourselves to consider such a kind of variations because, being $Du$ bounded, a parametric variation $\varphi_t$, for $t$ small enough preserves the property of being a graph.
\end{oss}

\subsection{Stability under parametric deformations}

We have shown that, given a solution to the minimal surface equation in $\Omega$, its graph minimizes the area among all graphs on $\Omega$ having the same boundary (theorem \ref{mingraf}). Actually, more is true, as the following theorem shows.

\begin{trm}\label{minnonpar}
Let $u:\Cl\Omega\rightarrow \R{}$ be a given Lipschits solution to the minimal surface equation \eqref{msediv} in $\Omega$. Then:
\begin{enumerate}
\item if $\Omega$ is homotopically trivial (for instance, $\Omega$ convex, star-shaper or contractile), then the graph of $u$ minimizes the area among every Lipschitz submanifold $\Sigma\subset\Cl\Omega\times\R{}$ having the same boundary;
\item if $\Omega$ is convex, then the graph of $u$ minimizes the area among all the Lipschitz submanifolds  $\Sigma\subset\R{n+1}$ having the same boundary.
\end{enumerate}
\end{trm}

The proof is based on the existence of a calibration, that is an exact $n$-form $\omega$ of absolute value at most 1, whose restriction to $\M{G}_u$ is the area form.

\begin{prop}[Calibration]\label{calibrazione}
Let $\omega$ be an exact $n$-form in $\Omega\times\R{}$, such that $|\omega|\leq 1$, that is
$$\sum_{1\leq i_1<\ldots<i_n\leq n+1}\omega_{i_1\cdots i_n}^2\leq 1.$$
Let a Lipschitz submanifold $\Sigma_0\subset\Cl\Omega\times\R{}$ with regular boundary be given and assume that
$\omega\big|_{\Sigma_0}$ is the volume form of $\Sigma_0$. Then the area of $\Sigma_0$ is least among the Lipschitz submanifolds $\Sigma\subset\Cl\Omega\times\R{}$ such that
$\partial\Sigma=\partial\Sigma_0$.
\end{prop}
\begin{proof} Being $\omega$ exaxt, we may find an $(n-1)$-form $\eta$ such that $d\eta=\omega$. Let $\Sigma$ be as in the statement of the proposition; then, by Stokes' theorem and since the two submanifolds have the same boundary,
$$\int_{\Sigma-\Sigma_0}\omega=\int_{\partial\Sigma-\partial\Sigma_0}\eta=0.$$
On the other hand, since $|\omega|\leq 1,$
$$\M{A}(\Sigma)\geq\int_{\Sigma}\omega=\int_{\Sigma_0}\omega=\M{A}(\Sigma_0).$$
\end{proof}

\begin{prooftrm} We prove separately the two claims.

\textbf{1.} We consider in $\Omega\times\R{}$ the calibration form
$$\omega(x,y):=\frac{\Big( \sum_{i=1}^n (-1)^{n-i+1}D_iu(x)\widehat{dx^i}dy\Big)+dx^1\cdots dx^n}{\sqrt{1+|Du|^2}}.$$
The minimal surface equation \eqref{msediv} implies $d\omega=0$; being $\Cl\Omega\times\R{}$ homotopically trivial, its de Rham cohomology is zero, thus $\omega$ is exact. Moreover $|\omega|=1$ and the restriction of $\omega$ to $\M{G}_u$ is the volume form of $\M{G}_u$, thus $\omega$ is a calibration for $\M{G}_u$ and proposition \ref{calibrazione} applies because $\partial\M{G}_u$ is smooth.

\textbf{2.} The second claims follows from the first one: let $\Sigma\subset\R{n+1}$ whose boundary be $\partial\M{G}_u$, hence contained in $\Cl\Omega\times\R{}$. The projection of $\Sigma$ onto $\Omega\times\R{}$ is well defined in the following way: for each $x\in\R{n}$ let $\pi_1(x)$ be the point of $\Omega$ of least distance from $x$. Such a point exists by convexity of $\Omega$. Then
$$\pi(x,y):=(\pi_1(x),y),\quad \forall x\in\R{n},\;\forall y\in\R{m}.$$
This projection doesn't increase the area; to the Lipschitz surface (possibly with multiplicity) obtained we apply step 1 and get $$\M{A}(\Sigma)\leq\M{A}(\pi(\Sigma))\leq\M{A}(\M{G}_u).$$
\end{prooftrm}

\begin{oss}
The hypothesis on $\Omega$ are necessary: in \cite{hardt} R. Hardt, C. P. Lau and Fang-Hua Lin proved the existence of a solution of the minimal surface equation whose graph doesn't minimize the area among the $n$-submanifolds of  $\R{n+1}$ having the same boundary.
\end{oss}

\section{Existence}

Let $\Omega \subset \R{n}$ be an open, smooth, connected and bounded domain. As we have seen, the solvability of the Cartesian problem of Plateau in codimension 1 is deeply tied to the solution of the minimal surface equation. In the following theorem we show that, under suitable hypothesis on $\partial\Omega$, it is possible to find a smooth solution to the minimal surface equation with arbitrarily prescribed boundary value. This, thanks to proposition \ref{equivale}, is equivalent to the solution of the problem of Dirichlet.

\begin{equation}\label{mse}
\left\{  \begin{array}{ll}
\displaystyle \sum_{i,j=1}^{n}g^{ij}(Du)\frac{\partial^2 u}{\partial x^i\partial x^j}=0 & \mathrm{in}\, \Omega \\
u=\psi &  \mathrm{su}\,\partial\Omega\rule{0cm}{.6cm}
\end{array}
\right. 
\end{equation}
with $u,\,\psi\in C^\infty(\Omega)\cap C^0(\Cl\Omega),$ $g_{ij}=\delta_{ij}+D_i u D_j u$ and
$(g^{ij})=(g_{ij})^{-1}.$ Explicitily
$$g^{ij}(Du)=\delta_{ij} -\frac{D_iu D_ju}{1+\abs{Du}^2}.$$

\begin{oss} Equation \eqref{mse} is quasilinear and elliptic. On the other hand, it is not uniformly elliptic, that is we cannot find $\lambda>0$ such that
$$g^{ij}(p)\xi_i\xi_j\geq\lambda|\xi|^2, \quad \forall\xi,p\in\R{n}.$$
In general we prove that the ellipticity constant $\lambda$ depends (only) on $|p|$, lemma \ref{ellittlim}, in a suitable way. In particular we know that, for $|p|\rightarrow +\infty$, $\lambda\rightarrow 0$ slowly enough, see \eqref{stimaellitt}.
\end{oss}

To prove the existence theorem \ref{esistcod1}, we use a method based on the fixed point theorem of Caccioppoli-Schauder.

\begin{trm}\label{caccioppoli} Let $T:K\rightarrow K$ a completely continuous operator\footnote{continuous and sending bounded sets into relatively compact sets; we do not assume $T$ to be linear.} which sends a convex, closed, bounded subset $K$ of a Banach space $B$ into itself. Then $T$ has a fixed point, meaning that there exists $\Cl x\in K$ such that $T(\Cl x)=\Cl x.$
\end{trm}

\begin{prop}\label{puntofisso}
Let a Banach space $B$ be given and consider a completely continuous operator $T:B\rightarrow B$ and $M>0$ such that for each pair $(\sigma, u)\in [0,1]\times B$ satisfing
$ u=\sigma Tu$ we get  $\norm{u}<M$. Then $T$ has a fixed point.
\end{prop}
\begin{proof} Let $K=\left\{ u\in B\,| \, \norm{u}\leq M\right\}$ and define the operator

\begin{displaymath}
\Cl{T}(u):=\left\{
\begin{array}{ll}
\displaystyle T(u) & \textrm{if}\, T(u)\in K \\
\displaystyle M\frac{T(u)}{\norm{T(u)}} & \textrm{if}\, T(u)\in B \backslash K\rule{0cm}{.7cm}
\end{array}
\right.
\end{displaymath}

$\Cl{T}$ send $K$ into itself, so that the fixed-point theorem of  Caccioppoli-Schauder, theorem \ref{caccioppoli}, implies that $\Cl{T}$ has a fixed point $u\in K$. Were $\norm{T(u)}\geq M$, we'd have
\begin{equation}\label{punto}
u=\frac{M}{\norm{T(u)}}T(u), \quad \frac{M}{\norm{T(u)}}\in[0,1],
\end{equation}
thus $\norm{u}<M$ by hypothesis, absurd because \eqref{punto} implies that $\norm{u}=M$. So $\norm{u}<M$ and $T(u)=\Cl{T}(u)=u$.
\end{proof}

We shall apply this theorem to the Banach space of functions with H\"older continuous first derivatives $B=C^{1,\alpha}(\Cl{\Omega})$:

\begin{defin}[H\"older functions]\label{defholder} A function $u:=\Cl\Omega\rightarrow\R{d}$ is said to be H\"older continuous with H\"older constant $\alpha\in(0,1]$ if
$$[u]_\alpha:=\sup_{x,y\in\Cl\Omega,\atop x\neq y}\frac{|u(x)-u(y)|}{|x-y|^\alpha}<+\infty.$$
Clearly every H\"older continuous function is continuous; the vector space of such functions is denoted by $C^{0,\alpha}(\Cl\Omega)$ and, endowed with the following norm, is a Banach space:
$$\norm{u}_{0,\alpha}:=\norm{u}_{C^0}+[u]_\alpha.$$
At the same time we define the spaces $C^{r,\alpha}$ with $r\in\mathbb{N}$ as the spaces of functions whose derivatives up to the $r$-th order are H\"older continuous; the corresponding norm is
$$\norm{u}_{r,\alpha}:=\norm{u}_{C^r}+[D^ru]_\alpha$$
where we consider $D^r u$ as a function $\R{d}$-valued for some $d$.
\end{defin}

We consider on $B=C^{1,\alpha}(\Cl\Omega)$ the operator $\widetilde T$ associating to a function $u\in C^{1,\alpha}$ the only solution $v$ to the following Dirichlet problem, whose existence is granted by the theorem which follows.

\begin{equation}\label{mse2}
\left\{  \begin{array}{ll}
\displaystyle\sum_{i,j=1}^{n} g^{ij}(Du)\frac{\partial^2 v}{\partial x^i\partial x^j}=0 & \mathrm{in}\, \Omega \\
\displaystyle v=\psi &  \mathrm{on}\,\partial\Omega.
\end{array}
\right. 
\end{equation}

Such a solution exists in $C^{2,\alpha}(\Cl{\Omega}).$ Indeed $Du\in C^{0,\alpha}(\Cl\Omega),$ thus the coefficients $g^{ij}(Du)$ are H\"older continuous theorem \ref{esistellitt} applies. The inclusion operator $\pi:C^{2,\alpha}(\Cl{\Omega})\rightarrow C^{1,\alpha}(\Cl{\Omega})$ is compact thanks to the corollary \ref{corascoli} to Ascoli-Arzel\`a's theorem. We want to show that, under suitable hypothesis on $\Omega$, the operator
\begin{equation}\label{eqoperatore}
T:=\pi\circ\widetilde{T}:C^{1,\alpha}(\Cl{\Omega})\rightarrow C^{1,\alpha}(\Cl{\Omega})
\end{equation}
is completely continuous, verifies the \emph{a priori estimate} of proposition \ref{puntofisso} and, consequently, has a fixed point, solution of (\ref{mse}).

\begin{trm}\label{esistellitt} Let $a^{ij}\in C^{0,\alpha}(\Cl\Omega)$ be given and assume that they are elliptic and bounded, that is such that we may choose $\lambda,\Lambda>0$ giving
\begin{equation}
\lambda\abs{\xi}^2\leq \sum_{i,j=1}^n a^{ij}(x)\xi_i\xi_j\leq\Lambda\abs{\xi}^2,\, \forall x\in\Omega.
\end{equation}
Then the Dirichlet problem
\begin{equation}
\left\{
\begin{array}{ll}
\displaystyle\sum_{i,j=1}^n a^{ij}(x)\frac{\partial^2 u}{\partial x^i\partial x^j} &\mathrm{in}\;\Omega\\
u=\psi  &\mathrm{on}\;\partial\Omega
\end{array}
\right.
\end{equation}
admits one and only one solution in $C^{2,\alpha}(\Cl\Omega)$. Moreover there exists a constant $C=C(\Omega,\lambda,\Lambda)$ such that
\begin{equation}\label{stimaschauder2}
\norm{u}_{2,\alpha}\leq C(\Omega,\lambda,\Lambda)\norm{\psi}_{2,\alpha}.
\end{equation}
\end{trm}

\subsection{The theorem of Ascoli and Arzel\`a}
A sequence of functions
$u_j:\Omega\rightarrow\R{}$
is said to be \emph{equicontinuous} if for each $x_0\in\Omega,$ $\varepsilon>0$ there exists $\delta>0$ such that
$$|u_j(x)-u_j(x_0)|<\varepsilon,\quad\forall x\in B_\delta(x_0),\;\forall j.$$ 
The same sequence $u_j$ is said to be \emph{equibounded} if there exists $M>0$ such that
$$|u_j(x)|\leq M,\quad \forall x\in\Omega,\forall j.$$
\begin{trm}[Ascoli-Arzel\`a]\label{trmascoli} Every equibounded and equicontinuous sequence of functions
$$u_j:\Omega\rightarrow\R{}$$
admits a subsequence converging uniformly on compact subsets.
\end{trm}

\begin{cor}\label{corascoli}
The immersion $C^{r,\alpha}(\Cl\Omega)\rightarrow C^r(\Cl\Omega)$, $0<\alpha\leq 1$, $r\in\mathbb{N}$,  is compact.
\end{cor}

\begin{proof} Let $u_j$ be bounded in $C^{r,\alpha}(\Cl\Omega)$, that is $\norm{u_j}_{r,\alpha}\leq M$ for some $M>0$. Then the derivatives of highest order are equicontinuous thanks to the estimate
$$|D^r_j(x)-D^r_j(y)|\leq K|x-y|^\alpha,\quad \forall j\in\mathbb{N},x,y\in\Omega.$$
Moreover the lower order derivatives are equicontinuous by boundedness of the highest order derivatives. Applying the theorem of Ascoli and Arzel\`a to each derivative we conclude that there exists a subsequence for which all the derivatives of order less than $r$ converge uniformly.
\end{proof}

\subsection{A priori estimates}

We give an a priori estimate in $C^{2,\alpha}(\Cl\Omega)$ of the solutions of $u=\sigma Tu$, $\sigma\in [0,1],$ being $T$ defined in \eqref{eqoperatore}.

Such an estimate may be obtained in four steps:
\begin{enumerate}
\item $\sup_{\Cl{\Omega}} \abs{u}$
\item $\sup_{\partial \Omega} \abs{Du}$
\item $\sup_{\Cl{\Omega}} \abs{Du}$
\item $\norm{u}_{1,\alpha}$
\end{enumerate}

First of all, we observe that $u=\sigma Tu$ is equivalent to
\begin{equation}\label{mse3}
\left\{  \begin{array}{ll}
\displaystyle\sum_{i,j=1}^n g^{ij}(Du)\frac{\partial^2 u}{\partial x^i\partial x^j}=0 & \mathrm{in}\, \Omega \\
\displaystyle u=\sigma\psi &  \mathrm{on}\,\partial\Omega.
\end{array}
\right. 
\end{equation}

\begin{prop}[Elliptic maximum principle]\label{princmax} Let $a^{ij}\frac{\partial^2}{\partial x^i\partial x^j}$ be an \\ elliptic operator, i.e.
$$\sum_{i,j=1}^n a^{ij}(x)\xi_i\xi_j\geq\lambda|\xi|^2,\quad \forall \xi\in\R{n},$$
for some $\lambda>0$; let arbitrary real-valued functions $b_k$ be given. Then a solution $u\in C^2(\Omega)\cap C^0(\Cl\Omega)$ to
\begin{equation}\label{princmaxeq}
\sum_{i,j=1}^n a^{ij}(x)\frac{\partial^2 u}{\partial x^i\partial x^j}(x)
+\sum_{k=1}^n b^k(x)\frac{\partial u}{\partial x^k}(x)\geq 0,
\end{equation}
satisfies
$$\sup_{\Omega}u=\max_{\partial\Omega}u.$$
\end{prop}

\begin{oss} Applying the maximum priciple to $-u$ we obtain that if
$$\sum_{i,j=1}^n a^{ij}(x)\frac{\partial^2 u}{\partial x^i\partial x^j}
+\sum_{k=1}^n b^k(x)\frac{\partial u}{\partial x^k}(x)\leq 0,$$
then
$$\inf_{\Omega}u=\min_{\partial\Omega}u.$$
\end{oss}

Proofs of the maximum principle are available in almost every book of second order PDE, for instance \cite{gilbarg}.

Since for a fixed $u$ equation (\ref{mse2}) is linear and uniformly elliptic\footnote{for a fixed $u$, $a^{ij}(x):=g^{ij}(Du(x))$ is a function depending only on the variable $x$; being $|Du(x)|$ bounded $x$, lemma \ref{ellittlim} implies the uniform ellipticity of $a^{ij}$.}
, the first step is a simple application of the maximum principle to the solution $u$, thus
$$\sup_{\Cl{\Omega}} \abs{u}\leq \sup_{\partial \Omega} \abs{\psi}.$$

\subsection{Boundary gradient estimates}

In order to have boundary gradient estimates we use barriers; they enable us to prove proposition \ref{propgradbordo}.

The gradient estimates are one of the things most distinguishing the various functionals in the calculus of variations. The costruction of barriers, indeed, is closely related to the structure of the equation (or system), i.e. the coefficients $g^{ij}$, particularly to the behaviour of the ellipticity constant.

\begin{lemma}\label{lemmamax}
Let $u,v\in C^2(\Omega)\cap C^0(\Cl\Omega)$ be such that
\begin{equation}
\left\{  \begin{array}{ll}
\displaystyle \sum_{i,j=1}^n g^{ij}(Du)\frac{\partial^2 u}{\partial x^i\partial x^j}=0 & \mathrm{in}\, \Omega \\
\displaystyle \sum_{i,j=1}^n g^{ij}(Dv)\frac{\partial^2 v}{\partial x^i\partial x^j}\leq 0 & \mathrm{in}\, \Omega \rule{0cm}{0.7cm}\\
\displaystyle u\leq v &  \mathrm{on}\,\partial\Omega \rule{0cm}{0.5cm}
\end{array}
\right. 
\end{equation}
Then $u\leq v$ on all of $\Omega$.
\end{lemma}

\begin{proof} By the mean value theorem of Lagrange there exists $\xi\in (0,1)$ such that
$$g^{ij}(Du)=g^{ij}(Dv)+\sum_{k=1}^{n}\frac{\partial g^{ij}}{\partial p^k}(\xi Dv + (1-\xi)Du)\Big(\frac{\partial u}{\partial x^k} - \frac{\partial v}{\partial x^k}\Big).$$
Subtracting in the previous system and setting $w:= v-u$ we get
\begin{equation}
\left\{  \begin{array}{ll}
\displaystyle\sum_{i,j=1}^n g^{ij}(Dv)\frac{\partial^2 w}{\partial x^i\partial x^j} + \sum_{k=1}^nb^k \frac{\partial w}{\partial x^k} \leq 0 & \mathrm{in}\, \Omega \\
w\geq 0 &  \mathrm{on}\,\partial\Omega\rule{0cm}{0.6cm}
\end{array}
\right. 
\end{equation}
to which the maximum principle, proposition \ref{princmax}, applies. This yields $w\geq 0$ in $\Omega.$
\end{proof}

We have the tools to construct barriers. Let $d:\Omega\rightarrow\mathbb{R}$ the function \emph{distance from the boundary}, smooth in a neighborhood of the boundary ($\Omega$ is a smooth domain). We define
$$N_r:=\left\{x\in\Omega\, \big{|}\, d(x)<r\right\},\quad\Gamma_r:=\{x\in \Omega\, | \, d(x)=r\};$$
these domains are smooth for $r$ small enough and, consequently, we shall always consider $r$ small. We consider on $N_r$ a function $v$ of type

$$v(x)=\psi(x)+h(d(x)),$$
where $h:[0,r]\rightarrow\mathbb{R}_+$ is smooth and satisfies
\begin{equation}\label{condizioni}
h(0)=0,\,h'(t)\geq 1,\,h''(t)\leq 0.
\end{equation}
With these choices we get

$$(1+\abs{Dv}^2)\sum_{i,j=1}^ng^{ij}(Dv)v_{ij}\leq h'' + C h'^2 + h'^3 \Delta d.$$

The behaviour of $\Delta d$ is determined by the mean curvature of $\partial \Omega:$ if $\partial\Omega$ has nonnegative mean curvature
$\Delta d\leq 0$,\footnote{see \cite{giusti} or \cite{gilbarg}} thus,
$$(1+\abs{Dv}^2)\sum_{i,j=1}^n g^{ij}(Dv)v_{ij}\leq h'' + C h'^2.$$
Now, setting $h(d)=k \log(1+\rho d),$ we may choose the constants $k$ and $\rho$ in such a way that conditions (\ref{condizioni}) are satisfied, 
$h(r)\geq 2\sup_{\partial \Omega}\abs{\psi}$ and $h''+C(h')^2\leq 0,$ thus
\begin{equation}
\left\{  \begin{array}{ll}
\displaystyle \sum_{i,j=1}^n g^{ij}(Dv)\frac{\partial^2 v}{\partial x^i\partial x^j} \leq 0 & \mathrm{in}\, N_r \\
\displaystyle \rule{0cm}{0.6cm}v\geq u &  \mathrm{on}\,\partial N_r
\end{array}
\right. 
\end{equation}
This, together with lemma \ref{lemmamax}, implies $u\leq v$ in $N_r$.
Being $u=v$ on $\partial \Omega$ we obtain
\begin{equation}\label{rapporto}
\frac{u(x)-u(y)}{\abs{x-y}}\leq\frac{v(x)-v(y)}{\abs{x-y}}, \quad x \in\Omega,\, y\in\partial\Omega .
\end{equation}
The construction of a lower barrier yields the opposite inequality and the a priori boundary gradient estimate:

\begin{prop}\label{propgradbordo}
Let $\Omega$ be such that $\partial\Omega$ has everywhere nonnegative mean curvature. 
Then there exists a constant $c=c(\Omega,\psi)$ such that, for each $\sigma\in[0,1]$, every solution of the minimal surface equation with prescribed boundary data $\sigma\psi$
\begin{equation}\label{msesigma}
\left\{  \begin{array}{ll}
\displaystyle \sum_{i,j=1}^ng^{ij}(Du)\frac{\partial^2 u}{\partial x^i\partial x^j}=0 & \mathrm{in}\, \Omega \\
\displaystyle\rule{0cm}{.6cm}u=\sigma\psi &  \mathrm{on}\,\partial\Omega.
\end{array}
\right.
\end{equation} 
satisfies $$ \sup_{\partial\Omega}\abs{Du}\leq c.$$
\end{prop}
\begin{proof}
The barriers $v^+$ and $v^-$ constructed for $\psi$ work also for $\sigma\psi$. For a fixed $y\in\partial\Omega$ and a suitable orthonormal frame in $y$, we have
$$ Du(y)=\left(D^{\partial\Omega}u(y),D_{\nu}u(y)  \right),$$
being $\nu$ the interion normal to $\partial\Omega$ in $y$. Since $u=\psi$ on $\partial\Omega$ we have
$D^{\partial\Omega}u=D^{\partial\Omega}\psi$, while (\ref{rapporto}) estimates the normal component $D_{\nu}u$:
$$-k\rho=D_\nu v^-\leq D_\nu u\leq D_\nu v^+=k\rho,$$
thus $|D_\nu u|\leq k\rho$.
\end{proof}

\subsection{Interior gradient estimates}
The interior gradient estimates are consequence of the boundary gradient estimates thanks to the following lemma of Rad\`o.

\begin{lemma}
Let $u$ be a solution to the minimal surface equation (\ref{mse}), then
\begin{equation}\label{riduzione}
\sup\left\{\frac{\abs{u(x)-u(y)}}{\abs{x-y}}\,\Big{|}\, x,y\in\Omega   \right\}=
\sup\left\{\frac{\abs{u(x)-u(y)}}{\abs{x-y}}\,\Big{|}\, x\in\Omega,\, y\in\partial\Omega   \right\}
\end{equation}
\end{lemma}
\begin{proof}
Let $x_1,x_2\in\Omega,$ $x_1\neq x_2$ and  $\tau=x_2-x_1$. Define
$$u_\tau(x):= u(x+\tau),$$
$$\Omega_\tau:=\{x\;:\: x+\tau\in\Omega \}.$$
Both $u$ and $u_\tau$ minimize in $\Omega\cap\Omega_\tau$, which is nonempty.
By the comparison principle, there exists $z\in\partial(\Omega\cap\Omega_\tau)$ such that
$$|u(x_1)-u(x_2)|=|u(x_1)-u_\tau(x_1)|\leq|u(z)-u_\tau(z)|=|u(z)-u(z+\tau)|.$$
Observe that $\partial(\Omega\cap\Omega_\tau)\subset(\partial\Omega\cup\partial\Omega_\tau)$ and, thus, at least one of the points $z,$ $z+\tau$ belongs to $\partial\Omega$. Moreover both $z$ and $z+\tau$ belong to $\Cl{\Omega}$. 
\end{proof}

Thanks to this lemma and the maximum principle for $u$ and the barriers built in $N_r$ it's easy to show that there is an a priori estimates of the right hand side of (\ref{riduzione}). The left part gives an obvious estimate of $\sup_{\Omega}\abs{Du},$ and this yields the interior gradient estimates:\\

\begin{prop}\label{propgradint}
There exists a constant $C=C(\Omega,\psi)$ such that, for each $\sigma\in[0,1]$, a solution $u$ of the minimal surface equation with boundary data $\sigma\psi$ \eqref{msesigma} satisfies
$$\sup_{\Omega}|Du|\leq C.$$

\end{prop}

\subsection{The $C^{1,\alpha}(\Cl\Omega)$ a priori estimates}

In proposition \ref{propgradint} and by the maximum principle we have established an a priori $C^1(\Cl\Omega)$ estimate for the solutions of \eqref{msesigma}, that is an estimate of $\sup_\Omega |u|+\sup_\Omega |Du|$ . To apply the fixed point theorem of Caccioppoli-Schauder we need a $C^{1,\alpha}(\Cl\Omega)$ estimate. Such an estimate comes from a theorem of De Giorgi, theorem \ref{trmDG}, and has been obtained in its global version by O. Lady\v zhenskaya and N. Ural'tseva \cite{lady}.

\begin{prop}\label{proplady}
Let $u$ be a $C^2(\Cl\Omega)$ solution of
\begin{displaymath}
\left\{
\begin{array}{ll}
\displaystyle\sum_{i=1}^n D_i\frac{D_iu}{\sqrt{1+|Du|^2}} & \mathrm{in}\;\Omega,\\
u=\psi&\mathrm{on}\;\partial\Omega,
\end{array}
\right.
\end{displaymath}
with $\psi\in \Lip(\Cl\Omega).$ 
Then, for some $\alpha>0$, the following a priori estimate holds:
$$\norm{u}_{1,\alpha}\leq C\left(\Omega,\lambda,\Lambda\right)\norm{\psi}_{1,\alpha}. $$
\end{prop}

\begin{proof} The a priori estimates of $\sup_{\Omega}|u|+\sup_{\Omega}|Du|$ have already been proved. We may differentiate the minimal surface equation as in proposition \ref{propderiva}; this doesn't require the difference quotient method because we are assuming that $u\in C^2(\Omega)$. We obtain
$$D_i(a^{ij}(Du)D_j w)=0,\quad w:=D_su,$$
with elliptic and bounded coefficients
$$a^{ij}(x)=\frac{1}{\sqrt{1+|Du|^2}}\bra{\delta^{ij}-\frac{D_iu D_ju}{1+|Du|^2}},$$
\begin{equation}\label{ellitt}
\lambda|\xi|^2\leq \sum_{i,j=1}^n a^{ij}\xi_i\xi_j\leq\Lambda|\xi|^2.
\end{equation}
The ellipticity and boundedness constants $\lambda$ and $\Lambda$ are estimated a priori because they depend only on $\sup_\Omega|Du|$: using lemma \ref{ellittlim} yields
$$\lambda=\frac{1}{(1+|Du|^2)^{\frac{3}{2}}},$$
while we may always choose $\Lambda=1$.
Applying the estimate \eqref{stimaDG} of the theorem of De Giorgi we obtain an a priori estimate of $\norm{w}_{C^{0,\alpha}(\Cl\Omega)}$ which, being $w=D_su,$ is the thesis.
\end{proof}

\subsection{The existence theorem}

\begin{trm}\label{esistcod1}
Let $\Omega$ be a smooth, bounded, connected domain whose  boundary has nonnegative mean curvature. Then for each $\psi\in C^{2,\alpha}(\Cl\Omega),$ there exists a unique $C^\infty(\Omega)\cap C^{2,\alpha}(\Cl\Omega)$ solution to the Dirichlet problem for the minimal surface equation (\ref{mse}). 
\end{trm}
\begin{proof}
Uniqueness is consequence of the elliptic maximum principle given in lemma \ref{lemmamax}.

Thanks to the Schauder estimates, see theorem \ref{morrey}, we only need to prove the existence of a solution in $C^{2,\alpha}(\Cl\Omega)$.

Thanks to theorem \ref{puntofisso} and the following remarks and propositions, we have to prove only that the operator
$T=\pi\circ\widetilde{T}$ is completely continuous. Compactness is clear: the Schauder estimates \ref{stimaschauder2} imply that $\widetilde{T}$ sends bounded subsets of $C^{1,\alpha}(\Cl\Omega)$ into bounded subsets of $C^{2,\alpha}(\Cl\Omega),$ which are immersed as relatively compact subsets of $C^{1,\alpha}(\Cl\Omega)$ by $\pi$ thanks to the theorem of Ascoli-Arzel\`a, corollary \ref{corascoli}. To prove continuity consider a sequence
$$u^{(k)}\xrightarrow{C^{1,\alpha}} u$$
and the corresponding sequence $v^{(k)}:=\widetilde{T}u^{(k)}.$

Given a subsequence $u^{(k')}$, thanks to the compact immersion $C^{2,\alpha}\hookrightarrow C^2$, there exists a converging sub-subsequence $u^{(k'')}$ such that
$$v^{(k'')}\xrightarrow{C^2}v $$
We easily observe that also $v$ is a solution to (\ref{mse2}):

\[
\begin{CD}
a^{ij}(Du^{(k'')}) && (v^{(k'')})_{ij} && =0\\
@VC^{0,\alpha}VV @VVC^0V\\
a^{ij}(Du)&&v_{ij}&&=0
\end{CD}
\]
(the sum over $i$ and $j$ is understood) and by uniqueness we have $v=Tu.$ The arbitrariness in the choice of the first subsequence implies
$$Tu^{(k)}\xrightarrow{C^{1,\alpha}} Tu $$
proving, thus, the continuity.
\end{proof}

There are existence theorems in $C^{2,\alpha}(\Cl\Omega)$ or in $C^2(\Omega)\cap C^0(\Cl\Omega)$ even when the mean curvature condition on
$\partial\Omega$ isn't satisfied; in this case, anyway, restrictions on the boundary data $\psi$ are necessary. Indeed the mean curvature condition is sharp:

\begin{prop}\label{meanprecisa} Let $x_0\in\partial\Omega$ be a point where $\partial\Omega$ has negative mean curvature. Then for each neighborhood $U$ of $x_0$ in $\Cl\Omega$ and every $\varepsilon>0$ there exists $\psi:\Cl\Omega\rightarrow\R{}$ with $\spt\psi\subset U$ and $|\psi|<\varepsilon$ such that the Dirichlet problem for the minimal surface equation with boundary data  $\psi$ is not solvable among Lipschitz functions.
\end{prop}

For a proof see \cite{giusti}, chapter 12.

\subsection{Another existence theorem}

Due to proposition \ref{meanprecisa}, an existence theorem for domains with somewhere negative mean curvature needs assumptions on the data $\psi$; such assumptions cannot involve only a $C^0$ estimate, thus conditions on the first derivatives are quite natural. In this direction we have a theorem of Graham Williams \cite{williams} which we state without proof.

\begin{trm}\label{williamsesist} Let $\Omega\subset\R{n}$ be smooth and let be $0<K<\frac{1}{\sqrt{n-1}}$. Then there exist $\delta,C>0$ such that, if
$$\frac{|\psi(x)-\psi(y)|}{|x-y|}\leq K\quad \mathrm{and}\quad |\psi(x)-\psi(y)|<\delta\quad\forall x,y\in\Cl\Omega,\;x\neq y,$$
then the Dirichlet problem for the minimal surface equation with boundary data $\psi$ admits a solution in $C^\infty(\Omega)\cap C^{0,\frac{1}{2}}(\Cl\Omega).$
\end{trm}

We, at last, cite theorem \ref{esistpiccoli}, which we shall prove. It guarantees the existence in arbitrary codimension for sufficiently small boundary data in the $C^{2,\alpha}$ norm. Clearly Williams' theorem is stronger because it yields existence for boundary data small in the $C^1$ norm; on the other hand this theorem doesn't generalize to arbitrary codimension.

\section{Regularity}

The regularity of solutions to the Dirichlet problem for the minimal surface equations in the non-negative mean curvature case may be inferred by existence in $C^{2,\alpha}$ (theorem \ref{esistcod1}) and uniqueness (theorem \ref{mingraf}).

We ask, more generally, whether the solutions to the minimal surface equation are smooth inside $\Omega$. The following theorem gives a positive answer.

\begin{trm}\label{regolcod1}
Let $u:\Cl\Omega\rightarrow\R{}$ be a weak Lipschitz solution to the minimal surface equation \eqref{msediv}. Then $u$ is analytic in $\Omega$.
\end{trm}

\begin{proof} That $u$ is $C^{1,\alpha}$ is an immediate consequence of theorems \ref{trmDG} and \ref{propderiva}. The higher regularity is proved via Schauder estimates, as in theorem \ref{morrey}.
\end{proof}

\subsection{De Giorgi-Nash's theorem}

One of the most beautiful and important theorems in the theory of elleptic equations was proved in 1957 by Ennio De Giorgi \cite{degiorgi} and a few months later independently by John Nash \cite{nash}. This theorem solves the $19^{\mathrm{th}}$ problem of Hilbert:
\begin{quote}
Are the solutions of regular problems in the calculus of variations always necessarily analytic?
\end{quote}

\begin{trm}[De Giorgi-Nash]\label{trmDG} Assume the coefficients $a^{ij}\in L^\infty(\Omega)$ are elliptic, that is such that \eqref{ellitt} holds for some $\lambda,\Lambda>0$. Then every weak solution $u\in W^{1,2}_{\loc}$ 
of
\begin{equation}\label{diver}
D_i(a^{ij}D_j u)=0
\end{equation}
is H\"older continuous, i.e. $u\in C^{0,\alpha}_{\loc}(\Omega)$ for some $\alpha>0$.
Moreover if $u=\varphi$ on $\partial \Omega$, with $\varphi\in \Lip(\Cl\Omega)$, then $u\in C^{0,\alpha}(\Cl\Omega)$, $\alpha=\alpha(\Omega,\lambda,\Lambda)$, and there is a constant $C=C(\Omega,\lambda,\Lambda)$ such that
\begin{equation}\label{stimaDG}
\norm{u}_{C^{0,\alpha}(\Cl\Omega)}\leq C\norm{\varphi}_{C^{0,\alpha}(\Cl\Omega)}.
\end{equation}
\end{trm}

For the proof see \cite{lady}, theorem 14.1 page 201. In order to apply this theorem to the regularity theorem \ref{regolcod1} we have to prove that the first derivatives of $u$ satisfy equation \eqref{diver} for a suitable choice of $a^{ij}$. That's the content of the following proposition.

\begin{prop}\label{propderiva} Let $A(p):=\frac{p}{\sqrt{1+|p|^2}}$ for each $p\in\R{n}$. Then the $W^{1,2}_{\loc}$ solutions of
$$\diver A(Du(x))=0,$$
that is, the solutions of the minimal surface equation are $W^{2,2}_{\loc}$ functions. Moreover, if we set $w=D_s u$ ($s=1,\ldots,n$), we have
$$D_i(a^{ij}D_jw)=0,$$
where $a^{ij}(x):=\frac{\partial A_i}{\partial x^j}(Du(x))=\frac{1}{\sqrt{1+|Du|^2}}\bra{\delta^{ij}-\frac{D_iu D_ju}{1+|Du|^2}}.$
\end{prop}
A proof may be found in \cite{giusti}, theorem C.1; it's based on the difference-quotient method. The point is to prove that the functions
$$\tau_{h,s}Du(x):=\frac{Du(x+he_s)-Du(x)}{h},\quad h\in(0,\varepsilon)$$
have $L^2$ norm equibounded with respect to $h$. This implies $D^2u\in L^2(\Omega)$ and the theorem, then is an easy computation.

\chapter{The counterexamples of Lawson and Osserman}\label{capitololo}

The non-parametric theory of minimal surfaces in codimension 1 had been well developed by the early '70s: the problems of existence, uniqueness and regularity had been positively solved. On the other hand, very little was known on the non-parametric theory in higher codimension.

The counterexamples of Blaine Lawson and Robert Osserman, published in 1977 in \cite{lawson} show the reason of this: the results true in codimension 1 are false in higher codimension. In particular, the existence of minimal graphs (solutions of the minimal surface system \eqref{mssnonpar}) with prescribed boundary value is not guaranteed even for very regular domains. Also the uniqueness of the solution is false; moreover the solutions of the minimal surface system, which in codimension 1 minimize the area thanks to the convexity of the area functional, in higher codimension are not necessary stable. Finally, the smoothness of Lipschitz minimal graphs, provided in codimension 1 by the theorem of De Giorgi \cite{degiorgi}, is not true in higher codimension and in \cite{lawson} a Lipschitz but non $C^1$ minimal graph is shown. This is optimal because in \cite{harvey} this cone is proved to minimize the area, while in \cite{morrey} Morrey proved that a $C^1$ solution to the minimal surface system is analytic.

\section{Non existence}

The existence theorem in codimension 1, theorem \ref{esistcod1}, requires a geometric hypothesis on the domain $\Omega$, namely that the mean curvature of $\partial\Omega$ is everywhere non-negative and that $\Omega$ be regular. If this is the case, we have existence for every smooth boundary data. The following counterexample of Lawson and Osserman shows that, in codimension greater than 1, these hypothesis are not sufficient.

\begin{trm}\label{nonesist} Let $\psi:S^{n-1}\rightarrow S^{m-1}$ be of class $C^2$ non homotopic in $S^{m-1}$ to a constant map, $n> m\geq2$. Then there exists $R_0>0$ depending on $\psi$ such that the Dirichlet problem for the minimal surface system in $B^n$
\begin{equation}\label{dirichmss}
\left\{ \begin{array}{ll}
\displaystyle
\sum_{i=1}^{n}\frac{\partial}{\partial x^i} \Big(\sqrt{g}g^{ij} \Big)=0 
& j=1,\ldots,n\\
\displaystyle
\frac{\partial}{\partial x^i}\Big(\sqrt{g}g^{ij}\frac{\partial u^\alpha}{\partial x^j}\Big)=0
&\alpha=1,\ldots,m\rule{0cm}{.7cm}\\
\displaystyle u^\alpha\big|_{\partial\Omega}=R\psi^\alpha\big|_{\partial\Omega}\rule{0cm}{.7cm}
&\alpha=1,\ldots,m,
\end{array}
\right.
\end{equation}
has no solution if $R\geq R_0$.
\end{trm}

\begin{lemma} Let $\Sigma\subset U$ be a Lipschitz manifold of class $C^2$ in a neighborhood $V$ of its boundary. Then, if $\Sigma$ is minimal in the sense of varifolds in $U\subset\R{n+m}$ open and bounded, we have
\begin{equation}\label{formulabordo}
\M{A}(\Sigma)=\frac{1}{n}\int_{\partial\Sigma}x\cdot\nu(x)d\M{H}^{n-1}(x),
\end{equation}
where $\nu(x)$ is the exterior normal to $\partial\Sigma$ in $x$.
\end{lemma}

\begin{proof} Let
$$\Sigma_r:=\{x\in\Sigma:d(x,\partial\Sigma)\geq r\}.$$
By compactness of $\partial\Sigma$ there exists $r_0>0$ such that for $r\leq r_0$ we have $\Sigma\backslash\Sigma_r\subset V$.
Using the function \emph{distance from the boundary}, we build a smooth function $\varphi_r(x)=\varphi_r(d(x,\partial\Sigma))$ such that
$$\varphi_r\Big|_{\Sigma\backslash\Sigma_{\frac{r}{2}}}=0,\quad \varphi_r\Big|_{\Sigma_r}=1,
\quad 0\leq\varphi_r\leq 1.$$
Consider the vector fields $X(x)=x$ and $\widetilde X(x)=\varphi_r(x)X(x).$ Since $\Sigma$ is minimal and $\widetilde X\in C^1_0(\Sigma)$ we obtain
$$0=\int_\Sigma \diver^\Sigma \widetilde Xd\M{H}^n=\int_\Sigma \nabla^\Sigma\varphi_r\cdot X
d\M{H}^n+\int_\Sigma\varphi_r \diver^\Sigma Xd\M{H}^n.$$
The last term, being $\diver^\Sigma X=n$, satisfies
$$\int_\Sigma \varphi_r \diver^\Sigma Xd\M{H}^n\rightarrow n\M{A}(\Sigma).$$
Being $\varphi_r$ a function of the distance, we have
\begin{equation}\label{contoformulabordo}
\int_\Sigma \nabla^\Sigma\varphi_r\cdot X d\M{H}^n=\int_{\Sigma\backslash\Sigma_r} 
\nabla^\Sigma\varphi_r\cdot X d\M{H}^n\rightarrow \int_{\partial\Sigma}X\cdot
\nu d\M{H}^{n-1};
\end{equation}
indeed $\nabla^\Sigma\varphi_r\sim \nu \frac{\partial \varphi_r(\rho)}{\partial \rho}$ and 
$\int_0^r \frac{\partial \varphi_r(\rho)}{\partial\rho}=1$ for every $r$. We obtain \eqref{contoformulabordo} writing the integrals in charts and applying Fubini-Tonelli's theorem.
\end{proof}

\begin{prooftrm} In the first two steps we prove two inequalities giving an absurd.

\textbf{1.} Let $u$ be a solution of \eqref{dirichmss} and let $\Sigma_R$ be its graph. If $x\in\partial\Sigma_R,$ then $|x|=\sqrt{1+R^2}$. This, applied to 
\eqref{formulabordo}, yields
\begin{equation}\label{massagraf}
\M{A}(\Sigma_R)\leq\frac{\sqrt{1+R^2}}{n}\M{H}^{n-1}(\partial\Sigma_R),
\end{equation}
where we consider $\partial\Sigma_R$ a $C^2$ submanifold of $\R{n+m}$ of dimension $n-1$ and $\Sigma_R$ smooth in a neighborhood of the boundary. The regulatiry is given by Allard's theorem \ref{trmallard2}.
On the other hand
\begin{equation}\label{crescitabordo}
\M{H}^{n-1}(\partial\Sigma_R)= c R^{m-1}
\end{equation}
because, being $\pi_2(\Sigma_R)\subset RS^{m-1}$, $\pi_2:\R{n+m}\rightarrow \{0\}\times\R{m}$, we have that the projection $\pi_2T_x\partial\Sigma_R$ is an $(m-1)$-dimensional vector space in $\R{m}$. Thus a homotety of $\R{m}$ rescales the area of $\partial\Sigma_R$ by a factor $R^{m-1}$. Substituting \eqref{crescitabordo} in \eqref{massagraf}, we obtain
\begin{equation}\label{primadisug}
\M{A}(\Sigma_R)\leq\frac{c\sqrt{1+R^2}}{n}R^{m-1}\leq c_1R^m.
\end{equation}

\textbf{2.} Let $u$ be again a solution of \eqref{dirichmss}. We prove that the image of $u$ contains $B^m$. Were it not so, and were $y_0\in B^n$ a point outside the image of $u$, we would have a retraction
$$\phi:B^m\backslash\{y_0\}\rightarrow S^{m-1}.$$
We may, thus, define the homotopy 
$$F:[0,1]\times S^{n-1}\rightarrow S^{m-1}$$
given by
$$F(t,x)=\phi(u((1-t)x)).$$
It's clear that $F(0,x)=\psi(x)$ and $F(1,x)=u(0)$ for every $x\in S^{n-1}$, contradicting the hypothesis on $\psi$.

Let $x_0\in B^n$ be such that $u(x_0)=0$, and $\zeta:=(x_0,0)\in\Sigma_R$. Take $(x,y)\in\partial\Sigma_R$.
Then
$$|\zeta-(x,y)|\geq |y|=R,$$ whence
$$B_R(\zeta)\cap \partial\Sigma=\emptyset.$$
By the monotonicity formula, proposition \ref{propmonot}, we have that for $0<r\leq R$ the function
$$\sigma(r):=\frac{\M{A}(\Sigma_R)\cap B_r(\zeta)}{\omega_n r^n}$$
is monotone increasing and
$\lim_{r\rightarrow 0^+}\sigma(r)\geq 1.$
This implies that $\sigma(R)\geq  R^n$, that is
\begin{equation}\label{secondadisug}
\M{A}(\Sigma_R)\geq\M{A}(\Sigma_R)\cap B_R(\zeta)\geq\omega_n R^n.
\end{equation}
\textbf{3.} Combining \eqref{secondadisug} with \eqref{primadisug} we have
\begin{equation}\label{sigmar}
\omega_nR^n\leq\M{A}(\Sigma_R)\leq c_1 R^m.
\end{equation}
Since the exponent on the left is greater that the one on the right, this last inequality cannot be true for a sequence $R_i\rightarrow +\infty$. Let
$$R_0=\sup \{R:\omega_n R^n\leq c_1 R^m \}<+\infty.$$
It's clear that thanks to \eqref{sigmar}, $\Sigma_R$ cannot exist for $R>R_0$.
\end{prooftrm}

\begin{oss} The hypothesis $n>m\geq 2$ doesn't include the case of dimension $2$. Actually there exists an important theorem due to T. Rad\`o \cite{rado} saying that the Dirichlet problem for the minimal surface system in dimension $2$ and arbitrary codimension is always solvable when continuous boundary data are prescribed. Moreover such a solution $u$ is continuous up to the boundary and analytic in the interior. For a proof see \cite{lawson}.
\end{oss}

\section{Non uniqueness and non stability}

We have remarked that a solution to the minimal surface equation with prescribed boundary data is unique and minimizes the area among the other graphs on $\Omega$ having the same boundary. Both uniqueness and stability are not true in higher codimension.

\begin{trm} There exists an analytic function $\psi:S^1\rightarrow\R{3}$ such that the minimal surface system \eqref{dirnonpar} with boundary data $\psi$ has at least 3 analytic solutions. Moreover one of these is unstable.\footnote{A surface parametrized by a graph $\M{G}_u$ is said to be unstable if for every $\varepsilon>0$ there exists $v$ with $|u-v|<\varepsilon$ and $\M{A}(\M{G}_v)<\M{A}(\M{G}_u)$.} 
\end{trm}

\begin{cenni} Lawson and Osserman define as boundary data a function $\psi$ symmetric: if $\Gamma\subset\R{5}$ is the graph of $\psi$, then $\Gamma=\sigma(\Gamma)$, where
$$\sigma(x^1,x^2,y^1,y^2,y^3):=(-x^2,x^1,-y^1,-y^3,y^2).$$
Rescaling such a function they define $\Gamma_R:=\M{G}_{r\psi}.$

Making use of Morse's theory the prove the following proposition.

\begin{prop}\label{propmorse} Let $F:B^2\rightarrow \R{5}$ be a parametric minimal surface with boundary data $\Gamma_R$ and let $\lambda:\R{5}\rightarrow\R{}$ be defined by
$$\lambda(x^1,x^2,y^1,y^2,y^3)=:y^1.$$
Then $\lambda\circ F$ has exactly one critical point in $B^2$.
\end{prop}

Thanks to a work of Douglas \cite{douglas}, there exists a surface $\Sigma\subset\R{5}$ parametrized by $F:B^2\rightarrow \R{5}$ (not necessarily a graph over $B^2$) whose boundary is $\Gamma_R$ and which minimizes the area among the Lipschitz surfaces with boundary $\Gamma_R$ homeomorphic to $B^2$. Assume that such a surface is unique. Then, by proposition \ref{propmorse}, $\lambda$ has only one critical point in $\Sigma$. Since $\sigma\circ F$ parametrizes a surface having the same area ($\sigma$ is an isometry) and with the same boundary $\Gamma_R$, if we assume that $\Sigma$ is unique, we have $\Sigma=\sigma(\Sigma)$ and $\sigma(p)=p$ by uniqueness of the critical point. It follows that $p=0\in\Sigma$.

The condition $0\in\Sigma$, and the non existece of other critical points, the shape of $\Gamma_R$ and the symmetry of $\Sigma$ imply $\M{A}_\Sigma\geq 4\pi R$; it's not difficult, anyway, to construct a surface $\Sigma'$ homeomorphic to $B^2$ with boundary $\Gamma_R$ and area
$$\M{A}(\Sigma')\leq (2\pi+\varepsilon)R^2+O(R)+\pi.$$
Consequently we have an absurd for $R$ large enough.
\end{cenni}

After proving the existence of two parametric area minimizing surfaces, the existence of a third parametric minimal surface follows from a work of Morse and Tompkins \cite{morse}. Finally the three surfaces found are non-parametric\footnote{namely they are the graph of a function.} thanks to a theorem of Rad\`o \cite{rado}.

\medskip

\begin{oss} Thanks to this counterexample, the study of the minimal surface system is not equivalent to the solution of the problem of Plateau. It's clear the difference with the codimension 1: a solution of the minimal surface equation minimizes the area among graphs and, if $\Omega$ is convex, minimizes the area also among the parametric surfaces (theorem \ref{minnonpar}).
\end{oss}

\medskip

\section{Non regularity: existence of minimal cones}

The regularity problem for non-parametric minimal surfaces is tightly connected to the nature of the minimal surface system. In codimension 1 this system reduces to an elliptic equation in divergence form \eqref{msediv}; solutions to such an equation are smooth thanks to the theorem of De Giorgi \cite{degiorgi}. There is no analogous theorem for elliptic systems and we may see an example of the difference between elliptic systems and elliptic equations in the following example.

\begin{trm}\label{controregol} Let $\eta:S^3\subset \R{4}\rightarrow S^2\subset\R{3}$ be the Hopf's map defined by
$$\eta(z_1,z_2)=(|z_1|^2-|z_2|^2,2z_1\Cl z_2)\in\mathbb{R}\times\mathbb{C}\cong\R{3},$$
where $(z_1,z_2)\in\mathbb{C}^2\cong\R{4}$. Then the Lipschitz, but not $C^1$, map
$$u:\R{4}\rightarrow\R{3}$$
given by
\begin{equation}\label{conohopf}
u(x)=\frac{\sqrt{5}}{2}|x|\eta\bra{\frac{x}{|x|}},\quad x\neq0
\end{equation}
and $u(0)=0$ satisfies the minimal surface system \eqref{mssnonpar}. 
\end{trm}

\noindent\textbf{Sketch of the proof} Consider the family of immersions
$$i_\alpha:S^3\rightarrow S^6$$ given by
$$i_\alpha(x):=\big(\alpha x,\sqrt{1-\alpha^2}\eta(x)\big),\quad 0\leq\alpha\leq 1.$$
Let $SU(2)$ be the group of unitary matrices of determinant 1 in $\mathbb{C}^2$ and let $SO(3)$ be the group of orthogonal matrices of determinant 1 in $\R{3}$.
We know that $SO(3)\cong SU(2)/\mathbb{Z}_2$ whence there is a natural immersion
$$SU(2)\hookrightarrow SU(2)\times SO(3),$$
through which $SU(2)$ acts on $S^6\subset\R{7}$. Thanks to a theorem of Wu-Yi Hsiang \cite{hsiang}, the orbits of highest volume of the action are minimal submanifolds of $S^6$.\footnote{$S^6$ is a Riemannian manifold. Given a submanifold $\Sigma$ the Levi-Civita connections on $S^6$ and $\Sigma$ are well defined. Consequently we define the mean curvature $H$ of $\Sigma$ in $S^6$; this, in general, doens't agree with the mean curvature of $\Sigma$ seen as a submanifold of $\R{7}$. We shall say that $\Sigma$ is minimal in $S^6 $ if $H=0$.} 
Using the symmetry of $i_\alpha$, it's not difficult to show that $i_\alpha(S^3)$ is a principal orbit and that it's enough to maximize the area among such submanifolds of $S^6$:
$$\M{A}(i_\alpha(S^3))=2\pi^2\alpha(4-3\alpha^2),$$
which attains it's maximum in $[0,1]$ when $\alpha=\frac{2}{3}$.
Then $i_{\frac{2}{3}}(S^3)$ is minimal in $S^6$, whence the cone $C$ built on it is minimal in $\R{7}$ because, in general, the cone built on a minimal submanifold of $S^n$ is a minimal submanifold of $\R{n+1}$. Indeed the mean curvature of $C$ in $x$ has no component parallel to $x$, thus is equal, up to rescaling, to the mean curvature of $i_{\frac{2}{3}}(S^3)$ in $S^6$.
Finally we verify that $C$ is the graph of the function $f$ defined in \eqref{conohopf}.

\medskip

\begin{oss}
Thanks to this counterexample the result of Morrey, theorem \ref{morrey}, saying the a $C^1$ solution to the minimal surface system is analytic, is sharp. To prove a regularity result for Lipschitz solutions we need further hypothesis: we will prove smoothness of Lipschitz minimal graphs which are \emph{area-decreasing}.
\end{oss}

\chapter{Existence in arbitrary codimension}\label{capitolocodarb}

The existence theorem in codimension 1 is based on the a priori boundary estimates for the gradient of a solution $u$. These are obtained by using barriers $v$ satisfying
\begin{equation}\label{facile}
\sum_{i,j=1}^ng^{ij}(Dv)\frac{\partial^2 v}{\partial x^i\partial x^j}\leq 0.
\end{equation}
Actually, being the minimal surface system \emph{non linear}, we should search for barriers of the form
\begin{equation}\label{difficile}
\sum_{i,j=1}^ng^{ij}(Du)\frac{\partial^2 v}{\partial x^i\partial x^j}\leq 0
\end{equation}
and this is more difficult because, since we don't know $Du$ (this is what we want to estimate!), we don't know the ellipticity constant of the coefficients $g^{ij}(Du)$.
In codimension 1 lemma \ref{lemmamax} shows that the construction of barriers satisfying \eqref{facile} is enough, but, being this lemma false in codimension greater than 1, we cannot generalize this procedure.

\medskip

In this chapter we will show that, with suitable assumptions on the $C^2$ norm of the boundary data $\psi$, the construction of barriers is possible. In order to do so, we shall use the parabolic system associated to the minimal surface system (the mean curvature system) and the parabolic maximum principle, showing that there are geometric quantities preserved along the mean curvature flow. Because of the counterexample of Lawson and Osserman, theorem \ref{nonesist}, it's natural to introduce hypothesis on the boundary data.

\medskip

We will use the mean curvature flow to prove the existence of a Lipschitz solution to the minimal surface system; in the following chapter we will prove that the solutions found here (not every solution, in general) are $C^\infty$. These results are due to Mu-Tao Wang, and appeared in 2003 and 2004 in \cite{wangbern} and \cite{wangdir}.

In the section below, we prove a result which was well know at the time of the counterexamples of Lawson and Osserman: a local existence theorem based on the inverse function theorem.

\section{Existence for $C^{2,\alpha}$-small data}

We show that when $\norm{\psi}_{2,\alpha}$ is small enough (depending on $\Omega$), then the Dirichlet problem for the minimal surface system \eqref{dirnondiv} has a smooth solution. This section is independent of the other sections and its results won't be used in what follows.

\begin{trm}[Inverse function]\label{trminversa}
Let $E$ and $F$ be Banach spaces and let $\Phi:E\rightarrow F$ be of class $C^r$, $r\geq1$. Assume that $D\Phi_{x_0}$ is an isomorphism of Banach spaces\footnote{this means that $D\Phi_{x_0}$ is invertible and its inverse is continuous. The latter hypothesis is unnecessary thanks to the open mapping theorem: a linear continuous and surjective map between Banach spaces is open.} for some $x_0\in E$. Then there exists $U\subset E$ and $V\subset F$ open neighborhoods of $x_0$ and $\Phi(x_0)$ respectively such that $\Phi(U)=V$,
$$\Phi\big|_{U}:U \rightarrow V$$
is invertible and the inverse is of class $C^r$.
\end{trm}

The proof of this theorem is the same as the proof of the theorem in $\R{n}$. For the details see \cite{lang}.

\begin{trm}\label{esistpiccoli} Given $\Omega\subset\R{n}$ open, connected and smooth, there exists a constant $C=C(\Omega)$ such that if $\norm{\psi}_{2,\alpha}<C$, then the Dirichelt problem for the minimal surface system with boundary data $\psi$ \eqref{dirnondiv} has a smooth solution.
\end{trm}

\begin{proof} We intend the sum over repeated indices and consider the Banach space operator
$$\Phi:C^{2,\alpha}(\Cl\Omega;\R{m})\rightarrow C^{0,\alpha}(\Omega;\R{m})\times C^{2,\alpha}(\partial\Omega;\R{m})$$
defined by
$$\Phi(u):=\Big(g^{ij}(Du)\frac{\partial^2 u}{\partial x^i\partial x^j},
u\big|_{\partial\Omega}\Big).$$

The differential of $\Phi$ in $u$ is
$$d\Phi_u(v)=\Big( g^{ij}(Du)\frac{\partial^2 v}{\partial x^i\partial x^j}
+ \frac{\partial g^{ij}}{\partial p^\beta_k}(Du)\frac{\partial v^\beta}{\partial x^k}
\frac{\partial^2 u}{\partial x^i\partial x^j},v\big|_{\partial\Omega}\Big).$$
It's easily seen that $d\Phi$ is continuous and that for $u=0$ it reduces to
$$d\Phi_0(v)=\big(\Delta v,v\big|_{\partial\Omega}\big);$$
inverting $d\Phi_0$ means: given $f\in C^{0,\alpha}(\Omega)$ and $h\in C^{2,\alpha}(\partial\Omega)$, solve in $C^{2,\alpha}(\Cl\Omega)$ the Dirichlet problem
\begin{displaymath}
\left\{
\begin{array}{ll}
\Delta v=f & \mathrm{in} \; \Omega\\
v=h & \mathrm{on}\;\partial\Omega.
\end{array}
\right.
\end{displaymath}
As well known this problem has always a solution in $C^{2,\alpha}(\Cl\Omega)$ and the solution is unique (maximum principle). Thus the operator $d\Phi_0$ is invertible and its inverse is continuous by the open mapping theorem.

Note that $\Phi(0)=0$; then the inverse function theorem \ref{trminversa} implies the existence of a neighborhood of $0$ $V\subset C^{0,\alpha}(\Omega)\times C^{2,\alpha}(\partial\Omega)$ contained in the image of $\Phi$. In particular there exists $C>0$ such that $\{0\}\times B_C(0)$ is contained in the image of $\Phi$ and this, together with theorem \ref{morrey}, concludes the proof.
\end{proof}

\section{Parabolic linear equations}

We state, without proof, some results from the theory of parabolic linear equations; a classical reference is the book of Lady\v zhenskaya, Ural'tseva and Solonnikov \cite{lady2}, see also Lieberman \cite{lieberman}.

A second order linear parabolic equation is a differential equation which may be written in the form
\begin{equation}\label{eqparab}
\frac{\partial u}{\partial t}(x,t)-\sum_{i,j=1}^na^{ij}(x,t)\frac{\partial^2 u}{\partial x^i\partial x^j}(x,t)=0
\end{equation}
in a domanin $\Omega_T:=\Omega\times (0,T)\subset\R{n+1}$, where the coefficients $a^{ij}$ are elliptic and bounded, that meaning that there exist $\lambda,\Lambda>0$ satisfying
$$\lambda|\xi|^2\leq a^{ij}(x,t)\xi_i\xi_j\leq\Lambda|\xi|^2, \quad \forall \xi\in\R{n}.$$
If we may choose $\lambda$ independently of $x$ and $t$, then equation \eqref{eqparab} is said to be \emph{uniformly parabolic}.
We denote by $\partial^*\Omega$ the parabolic boundary of $\Omega_T$, that is
$$\partial^*\Omega_T:=\Omega\times\{0\}\cup\partial\Omega\times[0,T).$$
The parabolic maximum principle is the parabolic analogue of the elliptic maximum principle \ref{princmax}:

\begin{prop}[Parabolic maximum principle] Given a solution $u:\Cl\Omega_T\rightarrow \R{}$ to the uniformly parabolic inequality
$$\frac{\partial u}{\partial t}(x,t)-\sum_{i,j=1}^na^{ij}(x,t) \frac{\partial^2 u}{\partial x^i\partial x^j}(x,t)\geq 0,$$
we have
$$\inf_{\Omega_T}u=\inf_{\partial^*\Omega_T}u,$$
while, if
$$\frac{\partial u}{\partial t}(x,t)-\sum_{i,j=1}^na^{ij}(x,t) \frac{\partial^2 u}{\partial x^i\partial x^j}(x,t)\leq 0,$$
we get $$\sup_{\Omega_t}u=\sup_{\partial^*\Omega_T}u.$$
\end{prop}

We consider the initial-boundary data problem for a linear parabolic equation with nonconstant coefficients
\begin{equation}\label{parab}
\left\{  \begin{array}{ll}
\displaystyle \frac{\partial u}{\partial t}(x,t)-\sum_{i,j=1}^na^{ij}(x,t)\frac{\partial^2u}{\partial x^i\partial x^j}(x,t)=0 & \mathrm{in}\, \Omega_T \\
\displaystyle u=\psi &  \mathrm{on}\,\partial^*\Omega_T
\end{array}
\right. 
\end{equation}
and study its solvability in suitable weighted H\"older spaces.

\begin{defin} Let $\Omega_T=\Omega\times(0,T)$, with $\Omega\subset\R{n}$ smooth and bounded domain; for $u:\Omega_T\rightarrow \R{d}$,
$\alpha\in (0,1]$ we define the seminorms
$$[u]_\alpha=\sup_{X,Y\in \Omega_T \atop X\neq Y}\frac{\abs{u(X)-u(Y)}}{\abs{X-Y}^\alpha},$$
where $X=(x,t)$, $Y=(y,s)$ and $\abs{X-Y}= \max \{ \abs{x-y}, \abs{t-s}^{1/2}\}$ is the parabolic distance;
$$[u]_{1,\alpha}= [Du]_\alpha $$
$$[u]_{2,\alpha}= [D^2u]_{\alpha}+ [D_tu]_\alpha$$
We define the corresponding norms
$$\norm{u}_{\alpha}= \sup_{\Omega_T} \abs{u} + [u]_\alpha$$
$$\norm{u}_{1,\alpha}= \sup_{\Omega_T} \abs{u}+\sup_{\Omega_T}\abs{Du} + [u]_{1,\alpha}$$
$$\norm{u}_{2,\alpha}= \sup_{\Omega_T} \abs{u}+\sup_{\Omega_T}\abs{Du}+\sup_{\Omega_T}\abs{D^2u}+\sup_{\Omega_T}\abs{u_t} 
+ [u]_{2,\alpha}$$
\end{defin}
We define the corresponding H\"older spaces
$$C^{r,\alpha}(\Omega_T)=\left\{u:\Omega\rightarrow\R{m}\, |\, \norm{u}_{r,\alpha} <\infty  \right\},
\quad r=0,1,2$$
which, endowed with the corresponding norms, are Banach spaces.

Observe that the first derivatives with respect to $t$ are treated as second order space derivatives; something similar happens in the definition of parabolic distance.

\begin{oss}
H\"older continuous functions are uniformly continuous, thus their continuity extends up to the boundary of $\Omega_T$. For instance a function in $C^{2,\alpha}(\Omega_T)$ has second order space derivatives and first order time derivative continuous up to the boundary.

This makes the spaces introduced up to now insufficient in the study of the solvability of an initial-boundary value problem as \eqref{parab}. Indeed, if the data $\psi$ don't satisfy the \emph{compatibility condition}
$$\frac{\partial \psi}{\partial t}(x,t)-\sum_{i,j=1}^na^{ij}(x,t)\frac{\partial^2\psi}{\partial x^i\partial x^j}(x,t)=0,\quad\mathrm{on}\;\partial^*\Omega,$$
then it's clearly impossible that a $C^{2,\alpha}(\Omega_T)$ function satisfy \eqref{parab}.

For this reason we introduce some weighted spaces using a \emph{distance from the boundary} function.
\end{oss}

\begin{defin}[Weighted parabolic H\"older spaces]
We introduce the function \emph{distance from the parabolic boundary}:
$$d(x)=\dist(x,\partial^*\Omega);\;d(X,Y):=\min\{d(X),d(Y)\}.$$
Define
$$[u]_0^*=\osc_\Omega f, \; \abs{f}_0^*=\sup_\Omega \abs{f}$$
\begin{displaymath}
\abs{u}_0^{(\delta)}=\left\{
\begin{array}{ll}
\sup_{\Omega_T} d^\delta\abs{u} & \mathrm{if}\, b\geq 0\\
\sup_{\Omega_T} (\diam \Omega)^\delta\abs{f} & \mathrm{if}\, b< 0\\
\end{array}
\right.
\end{displaymath}
$$[u]_\alpha^{(\delta)}=\sup \left\{ d(X,Y)^{\delta+\alpha}\frac{\abs{u(X)-u(Y)}}{\abs{X-Y}^\alpha}\right\}$$
$$[u]_{1,\alpha}^{(\delta)}=\sup 
\left \{d(X,Y)^{\delta +1+\alpha}\frac{\abs{Du(X)-Du(Y)}}{\abs{X-Y}^\alpha}\right\}$$
$$[u]_{2,\alpha}^{(\delta)}=\sup 
\left\{d(X,Y)^{\delta +2+\alpha}\left( \frac{\abs{D^2u(X)-D^2u(Y)}}{\abs{X-Y}^\alpha}+ \frac{\abs{D_tu(X)-D_tu(Y)}} {\abs{X-Y}^\alpha}\right) \right\}$$
$$\norm{u}_{\alpha}^{(\delta)}= \abs{u}_0^{(\delta)}+[u]_\alpha^{(\delta)}$$
$$\norm{u}_{1,\alpha}^{(\delta)}= \abs{u}_0^{(\delta)}+\abs{Du}_0^{(1+\delta)}+[u]_{1,\alpha}^{(\delta)}$$
$$\norm{u}_{2,\alpha}^{(\delta)}= \abs{u}_0^{(\delta)}+\abs{Du}_0^{(1+\delta)}+\abs{D^2u}_0^{(2+\delta)}+
\abs{u_t}_0^{(2+\delta)}+[u]_{2,\alpha}^{(\delta)}$$
\end{defin}

Consequently we define the spaces $C^{k,\alpha}_{(\delta)}$, $k=0,1,2$. The functions in these spaces are continuous with the proper number of derivatives in the interior, but not in general up to the boundary.

The main existence result from the linear theory we use is contained in the following theorem.

\begin{trm}\label{esistparab}
Set $\alpha\in(0,1)$, $\delta\in(1,2)$. Assume that the coefficients $a^{ij}$ satisfy the uniform ellipticity condition  \eqref{ellitt} in $\Omega_T$, are H\"older continuous, meaning that 
$$\norm{a^{ij}}_\alpha^{(0)}<\infty$$
and assume
\begin{equation}\label{eqmod}
\abs{a^{ij}(X)-a^{ij}(Y)}\leq \zeta(\abs{X-Y})
\end{equation}
for some continuous increasing function $\zeta$ with $\zeta(0)=0$.
Then the probelm \eqref{parab} with initial-boundary data $\psi \in C^\delta$ is uniquely solvable in $C^{2,\alpha}_{(-\delta)}$ and
\begin{equation}\label{stimaesistenza}
\norm{u}_{2,\alpha}^{(-\delta)}\leq C\left( \norm{a^{ij}}_\alpha^{(0)},\lambda,\Lambda,\Omega_T,\zeta  \right)\norm{\psi}_\delta
\end{equation}
\end{trm}

\begin{oss} If we choose $\delta$ and $\alpha$ such that $\delta - \alpha > 1+\theta$, $\theta\in (0,l)$ we obtain
$$[Du]_\theta \leq\norm{u}_{1,\alpha}^{(-\delta)}\leq C \norm{u}_{2,\alpha}^{(-\delta)}<\infty,$$
which implies that the solution has first derivatives H\"older continuous and, thus, continuous up to the boundary.
\end{oss}

If the compatibility condition is satisfied, it's possible to obtain $C^{2,\alpha}$-regularity up to the boundary:

\begin{trm}\label{esistcomp} Set $\alpha\in(0,1)$ and let the coefficients $a^{ij}$ in $\eqref{parab}$ be $\alpha$-H\"older continuous,  elliptic and satisfy \eqref{eqmod}. Also let $\psi\in C^{2,\alpha}(\Omega_T).$
Assume that the compatibility condition for the initial-boundary data $\psi$ holds:
$$\frac{\partial \psi}{\partial t}(x,t)-\sum_{i,j=1}^na^{ij}(x,t)\frac{\partial^2\psi}{\partial x^i\partial x^j}(x,t)=0
,\quad \mathrm{on}\;\partial\Omega\times\{0\}.$$
Then problem \eqref{parab} has a unique solution in $C^{2,\alpha}(\Cl\Omega).$
\end{trm}

\section{The Dirichlet problem for the minimal surface system}

Recall that the non-parametric Dirichlet problem in non-divergence form for the minimal surface system is
\begin{equation}\label{mss}
\left\{  \begin{array}{ll}
\displaystyle \sum_{i,j=1}^n g^{ij}(Du)\frac{\partial^2 u}{\partial x^i\partial x^j}=0 & \mathrm{in}\, \Omega \\
u=\psi &  \mathrm{on}\,\partial\Omega
\end{array}
\right. 
\end{equation}
with $\psi\in C^\infty(\Cl{\Omega};\R{m}),$ $g_{ij}(Du)=\delta_{ij}+\sum_l \frac{\partial u^l}{\partial x^i} \frac{\partial u^l}{\partial x^j}$ and $(g^{ij}(Du))=(g_{ij}(Du))^{-1}$.

\medskip

To solve this quasilinear elliptic system, we study the associated parabolic system, which corresponds to the non-parametric mean curvature flow:
\begin{equation}\label{mcs}
\left\{  \begin{array}{ll}
\displaystyle \frac{\partial u}{\partial t}-\sum_{i,j=1}^n g^{ij}(Du)\frac{\partial^2 u}{\partial x^i\partial x^j}=0 & \mathrm{in}\, \Omega_\infty \\
u=\psi &  \mathrm{on}\,\partial^*\Omega_\infty,
\end{array}
\right. 
\end{equation}
with $\psi\in C^\infty(\Cl\Omega_\infty)$ and $u=(u^1,\dots,u^m)\in C^2(\Omega_\infty)\cap C^0(\Cl\Omega_\infty).$

\medskip

\section{The mean curvature flow}

Let be given an $n$-submanifold $\Sigma$ in $\R{n+m}$ parametrized by
$F:\Cl\Omega\rightarrow \R{n+m}$. The mean curvature flow of $\Sigma$ is a family of embeddings
$$F_t:\Cl\Omega\rightarrow\R{n+m},\quad t\in [0,T)$$ such that, defined $F(x,t)=F_t(x),$
\begin{equation}\label{mcfpar}
\left\{
\begin{array}{ll}
\displaystyle \frac{\partial F}{\partial t}(t,x)=H(F(t,x)) &\mathrm{in}\;\Omega_T\\
\displaystyle F=F_0 &\mathrm{on}\;\partial^*\Omega_T,\rule{0cm}{.6cm}
\end{array}
\right.
\end{equation}
where $H(F(t,x))$ is the mean curvature vector of the submanifold $\Sigma_t:=F_t(\Omega)$ in $F(t,x)$.

The mean curvature flow is the minus gradient flow of the area functional with respect to the $L^2$ scalar product as may be seen in \eqref{varcurvmedia}: we deform a given surface $F_0(\Omega)$ in the direction in which the area decreases most.

The relation between the parametric mean curvature flow and the system in \eqref{mcs} is described in the following proposition.

\begin{prop}\label{propmcf}
Let $F:\Cl\Omega_T\rightarrow\R{n+m}$ be a solution to the parametric mean curvature flow \eqref{mcfpar} and assume that $\Sigma_t:=F_t(\Omega)$ may be written, for each $t\in[0,T)$, as tha graph of a function with gradient bounded on $\Omega$. Then there exists a family of diffeomorphisms leaving the boundary fixed $r_t:\Cl\Omega\rightarrow\Cl\Omega$ such that $\widetilde F_t=F_t\circ r_t$ is of the form
$$\widetilde F(x)=(x,u(x)),\quad x\in\Omega,$$
where
$$u:\Cl\Omega_T\rightarrow \R{m}$$
solves \eqref{mcs} with initial-boundary data $\psi:\Cl\Omega\rightarrow \R{m}$ such that
$$F_0(x)=(x,\psi(x)).$$
Conversly, given $u$ solution of \eqref{mcs}, the family of embeddings
$$\widetilde F_t:=I\times u:\Omega\rightarrow\R{n+m}$$ satisfies
$$\Big(\frac{\partial \widetilde F}{\partial t}(t,x)\Big)^N=\widetilde H(t,x),$$
where $\widetilde H(t,x)$ is the mean curvature of $\Sigma_t:=\widetilde F_t$ in $\widetilde F(t,x)$.
\end{prop}

\begin{proof}
The map
$$\pi_1\circ F_t:\Cl\Omega\rightarrow \Cl\Omega$$
is a bijection because we assumed that $\Sigma_t$ is a graph on $\Cl\Omega$. Define
$$r_t:=(\pi_1\circ F_t)^{-1}.$$
Clearly $r_t$ is the identity on $\partial\Omega$.
Let $\widetilde F(t,x)=F(t,r(t,x))$; then
$$\frac{\partial \widetilde F}{\partial t}(t,x)=\frac{\partial F}{\partial t}(t,r(x,t))+
dF_t\Big(\frac{dr}{dt}\Big)(t,x),$$
whence
\begin{equation}\label{eqmoto}
\bra{\frac{\partial \widetilde F}{\partial t}(t,x)}^N=\bra{\frac{\partial F}{\partial t}(t,r(t,x))}^N
=\widetilde H(t,x).
\end{equation}
We observe that
$$\frac{\partial \widetilde F}{\partial t}\in\{0\}\times\R{m},\quad\sum_{i,j=1}^n g^{ij}(Du)\frac{\partial^2 \widetilde F}{\partial x^i\partial x^j}\in\{0\}\times\R{m};$$
thanks to \eqref{curvmedia2} we have
\begin{equation}\label{curvmediatilde}
\Big(\sum_{i,j=1}^n g^{ij}(Du)\frac{\partial^2 \widetilde F}{\partial x^i\partial x^j}\Big)^N=\widetilde H,
\end{equation}
moreover the projection of $\{0\}\times\R{m}$ onto $N_{\widetilde F(t,x)}\Sigma_t$ is injective, thus we get
$$\frac{\partial \widetilde F}{\partial t}= \sum_{i,j=1}^n g^{ij}(Du)\frac{\partial^2 \widetilde F}{\partial x^i\partial x^j},$$
which is equivalent to \eqref{mcs}.
The converse is similar.
\end{proof}

In what follows we will not make use of the above proposition but in equation \eqref{curvmediatilde} which enables us to compute the variation of the area of the graphs moving by mean curvature, i.e. solving \eqref{mcs}.

\medskip

We prove the existence for all times of the mean curvature flow using a method of continuity, as done by Mu-Tao Wang \cite{wanglongtime}: we show that the set of times for which the solution exists is both open and closed\footnote{We say that a solution exists in $[0,t_0]$ if there exist for $(x,t)\rightarrow(x_0,t_0^-)$ the limits of the second order space derivatives and of the first order time derivative.}. Using the Caccioppoli-Schauder's fixed point theorem we prove the existence for small times, thus proving openess. For the closure we need the a priori estimates of the subsequent sections: we estimate the gradient on the boundary and on the interior and, thanks to a theorem of Brian White, also the higher order derivatives.

\subsection{Existence of the mean curvature flow for small times}

\begin{trm}\label{trmtempipiccoli} Let $\psi\in C^\delta(\Cl\Omega_\infty)$ for some $\delta\in(1,2)$. Then there exists a constant $\varepsilon>0$ such that the problem \eqref{mcs} has a solution $u\in C^{2,\alpha}_{(\delta)}(\Omega_{\varepsilon})$.
\end{trm}

\begin{proof}
Choose $\theta\in(1,\delta)$ and set $M:=1+[\psi]_\theta<+\infty$. For some $\varepsilon>0$ to be fixed we set
$$K=\{v\in C^\theta(\Omega_\varepsilon):[v]_\theta\leq M\}$$
and we define the non-linear operator
$$T:K\rightarrow C^{\theta}(\Omega_\varepsilon)$$
which to $u\in K$ associates the solution $Tu=v$ to the uncoupled linear system
\begin{equation}\label{mcs2}
\left\{  \begin{array}{ll}
\displaystyle \frac{\partial v}{\partial t}-\sum_{i,j=1}^n g^{ij}(Du)\frac{\partial^2 v}{\partial x^i\partial x^j}=0 & \mathrm{in}\, \Omega_\varepsilon \\
v=\psi &  \mathrm{on}\,\partial^*\Omega_\varepsilon,
\end{array}
\right. 
\end{equation}
Thanks to theorem \ref{esistparab}, such a solution exists in $C^{2+\alpha(\theta+1)}_{(\delta)}(\Omega_\varepsilon)$ and 
$$[v]_1\leq[v]_\delta\leq C[v]_{2+\alpha(\theta-1)}^{(-\delta)}\leq C(M).$$
Consequently $|v-\psi|\leq C\varepsilon$ in $\Omega_\varepsilon$ and, by interpolation, $[v-\psi]_\theta\leq C
\varepsilon^{\frac{\delta-\theta}{\delta}}.$ Then $[u]_\theta\leq M$ for some $\varepsilon>0$ small enough which, from now on, we fix. The operator $T$ sends $K$ into itself; observing that, thanks to Ascoli-Arzel\`a's theorem, $K$ is a compact subset of $C^1(\Omega_\varepsilon)$ and that it's convex (seminorms are convex), we may apply Caccioppoli-Schauder's fixed point theorem and we obtaian a fixed point $u\in C^{2+\alpha(\theta-1)}_{(-\delta)}(\Omega_\varepsilon)$ for $T$. It is a solution to problem \eqref{mcs} in $\Omega_\varepsilon$ and, thanks to theorem \ref{esistparab}, $u\in C^{2+\alpha}_{(-\delta)}(\Omega_\varepsilon)$.

\end{proof}

\section{Boundary gradient estimates}

\begin{lemma}[Ellipticity and boundedness of $g^{ij}$]\label{ellittlim}
Let 
$$g_{ij}(Du)=\delta_{ij}+\sum_{\alpha=1}^m\frac{\partial u^\alpha}{\partial x^i}
\frac{\partial u^\alpha}{\partial x^j}$$ and $(g^{ij}(Du))$ the inverse matrix of $(g_{ij}(Du))$. Then we have  
\begin{equation}\label{stimaellitt}
\frac{1}{1+ \eta}\abs{\xi}^2\leq g^{ij}(Du)(x)\xi_i\xi_j\leq \abs{\xi}^2,\quad \eta=\sup_{\Omega_T}\abs{Du}^2,
\end{equation}
for every $\xi\in\R{n}$.
\end{lemma}

\begin{proof} 
It's clear that $(g_{ij})$ is symmetric, thus diagonalizable. Let
$\lambda_1\leq\cdots\leq\lambda_n$ be its eigenvalues (possibly repeated). Then also $(g^{ij})$ is diagonalizable and its eigenvalues are $\frac{1}{\lambda_1},\ldots,\frac{1}{\lambda_n}$. Now it's easy to verify that
$$\sum_{i,j=1}^ng_{ij}(Du)(\xi_i\xi_j)=|\xi|^2+\Big(\sum_{1\leq i\leq n \atop 1\leq \alpha\leq m}
\frac{\partial u^\alpha}{\partial x^i}\xi_i \Big)^2\geq |\xi|^2.$$ 
Let $\xi_*$ be an eigenvector relative to the smallest eigenvalue, say $\lambda_1$. The above equation implies that, denoting by $g$ the linear map associated to the matrix $(g_{ij})$, we have
$$|\xi_*|^2\leq\dual{\xi_*,g\xi_*}=\lambda_1|\xi_*|^2$$
and therefore $\lambda_i\geq 1.$
Similarly we obtain that, if $\lambda_n$ is the largest eigenvalue, then
$$\lambda_n=\sup_{|\xi|=1} |g(\xi)|=1+|Du|^2.$$
Estimate \eqref{stimaellitt} follows immediatly by the estimates on the eigenvalues.
\end{proof}

\begin{trm} Let $\Omega$ be bounded, convex, smooth and let $u\in C^{2,\alpha}_{(-\delta)}(\Omega_T;\R{m})$ be a solution to \eqref{mcs}.
Then we have the following estimate

\begin{equation}\label{stimabordo}
\abs{Du}<4n\diam\Omega(1+\eta)\sup_\Omega\abs{D^2\psi}+\sqrt{2}\sup_{\partial\Omega}\abs{D\psi}\; 
\mathrm{on}\,\partial\Omega\times[0,T),
\end{equation}
where $\eta=\sup_{\Omega_T}\abs{Du}^2$.
\end{trm}

Observe that the hypothesis $u\in C^{2,\alpha}_{(-\delta)}(\Omega_T;\R{m})$ implies $u\in C^2(\Omega_T;\R{m})\cap C^1(\Cl\Omega\times(0,T);\R{m})$.

\begin{proof}
Take $p\in\partial \Omega$ and set $\Gamma\subset\R{n}$ to be the hyperplane tangent to $\partial\Omega$ in $p$; define  $d$ as the function \emph{distance from $\Gamma\times(0,T)$} in $\Omega_T$, i.e.
$$d(x,t)=\dist(x,\Gamma).$$
Since $d$ is linear $\sum g^{ij}D_{ij}d=0.$ For a fixed $1\leq l\leq m$ we define the \emph{barrier}
$$v(x,t)=k\log(1+\rho d)-(u^l -\psi^l).$$
We compute
\begin{equation}\label{contobarriera}
\frac{\partial v}{\partial t} -g^{ij}(Du)\frac{\partial^2 v}{\partial x^i\partial x^j}=\frac{k\rho^2}{(1+\rho d)^2}g^{ij}(Du)\frac{\partial d}{\partial x^i}\frac{\partial d}{\partial x^j}-g^{ij}(Du)\frac{\partial ^2\psi^l}{\partial x^i\partial x^j}.
\end{equation}
We used that fact that $u$ is a solution and $\psi$ doesn't depend on $t$.
Thanks to the ellipticity estimate on $(g^{ij})$, inequality \eqref{stimaellitt}, and to $\abs{Dd}=1$ we have $g^{ij}D_idD_j d\geq\frac{1}{1+\eta}$, thus
$$\frac{\partial v}{\partial t}-g^{ij}\frac{\partial^2 v}{\partial x^i\partial x^j}\geq \frac{k\rho^2}{(1+\rho d)^2} \frac{1}{1+\eta}.$$
$$g^{ij}(Du)\frac{\partial^2\psi^l}{\partial x^i\partial x^j}\leq n\sup_{\Omega}\abs{D^2 \psi}.$$
Therefore, if
\begin{equation}\label{stimabar}
\frac{k\rho^2}{(1+\rho \diam\Omega)^2} \frac{1}{1+\eta}\geq n\sup_{\Omega}\abs{D^2 \psi},
\end{equation}
we get $v_t -\sum g^{ij}D_{ij}v\geq 0$ on $\Omega_T.$ 
Now $v(x,t)\geq 0$ on $\partial^*\Omega,$ thus the strong parabolic maximum principle implies $v > 0$ in $\Omega_T.$
Since $v(p,t)=0$ for every $t\in [0,T)$ we have
$$0 <\frac{\partial v}{\partial n}= k\rho - \frac{\partial(u^l-\psi^l)}{\partial n}, \; \mathrm{i.e.}$$
$$\frac{\partial u^l}{\partial n} < k\rho + \frac{\partial \psi^l}{\partial n}.$$
The construction of a lower barrier yields an analogous estimate for $-\frac{\partial u} {\partial n}$
whence $$\abs{\frac{\partial u}{\partial n}}< k\rho+\abs{\frac{\partial \psi}{\partial n}}$$
and, since $u\big|_{\partial\Omega}=\psi\big|_{\partial\Omega}$
$$ D^{\partial\Omega}u=D^{\partial\Omega}\psi,$$ we have
\begin{equation}\label{eqstimab}
\abs{Du}< \sqrt{\left( k\rho+\abs{ \frac{\partial \psi}{\partial n}} \right)^2 +\abs{D^{\partial\Omega} \psi}^2}
\leq k\rho+\sqrt{2}\abs{D\psi}\; \mathrm{in}\, p.
\end{equation}
To obtain \eqref{stimabar} we set $\rho=(\diam \Omega)^{-1}$ and 
$$ k=4n(\diam\Omega)^2(1+\eta)\sup_\Omega\abs{D^2\psi},$$
whence \eqref{eqstimab} becomes \eqref{stimabordo}.
\end{proof}

\section{Interior gradient estimates: the function $*\omega$}
We introduce the function
\begin{equation}\label{omega}
*\omega=\frac{1}{\sqrt {\det (I + Du^TDu)}}=\frac{1}{\sqrt {\prod_{i=1}^n (1+\lambda_i^2)}},
\end{equation}
where the number $\lambda_i$ are the singular values of $Du$, i.e. the square roots of the eigenvalues of $Du^TDu$.
The following relations are easily verified:
\begin{equation}\label{relaz}
*\omega>\frac{1}{\sqrt{2-\delta}}\Rightarrow \abs{Du}^2<1-\delta;\quad \abs{Du}<\sqrt{(2-\delta)^{1/n}-1} \Rightarrow *\omega>\frac{1}{\sqrt{2-\delta}}.
\end{equation}

Let $\omega$ be the $n$-form on $\R{n+m}$ defined by
$$\omega(e_1,\ldots,e_n)=1$$
$$\omega(e_{i_1},\ldots,e_{i_n})=0\quad\textrm{if}\; i_1<\ldots<i_n, \; i_n>n.$$

The covariant derivatives of a tensor are well defined on a Riemannian manifold:
consider in particular $\omega$ belonging to $\M{T}_n(\Sigma)$, the space of \emph{covariant} $n$-tensors. By definition
$$\nabla_X^\Sigma \omega(Y_1,\ldots,Y_n):=D_X\omega(Y_1,\ldots,Y_n)-\sum_{i=1}^n
\omega(Y_1,\ldots,\nabla^\Sigma_X Y_i,\ldots,Y_m).$$
Moreover the Laplacean of a tensor may be defined as
$$\Delta_\Sigma\omega=\nabla_{\tau_k}^\Sigma\nabla_{\tau_k}^\Sigma \omega.$$

\begin{lemma}\label{lemmalapl}
Let $\omega$ be defined as above for the Riemannian submanifold $\Sigma\subset\R{n+m}.$ Let be given an orthonormal basis $\{\tau_1,\ldots,\tau_n\}$ in a neighborhood of a fixed point $p$. Then in $p$ we have 
$$(\Delta_\Sigma \omega)(\tau_1,\ldots,\tau_n) =\Delta_\Sigma(\omega(\tau_1,\ldots,\tau_n))=\Delta_\Sigma*\omega.$$
\end{lemma}
\begin{proof} Set
$$\omega(\tau_1,\ldots,\tau_n)=\omega_{1\cdots n},\quad (\Delta_\Sigma \omega)(\tau_1,\ldots,\tau_n)=(\Delta_\Sigma \omega)_{1\cdots n}.$$
Then
\begin{multline}
(\Delta_\Sigma\omega)_{1\cdots n}=
D_{\tau_k}\bra{\nabla_{\tau_k}^\Sigma\omega(\tau_1,\ldots,\tau_n)}-\sum_i D_{\tau_k}\omega
(\tau_1,\ldots,\nabla^\Sigma_{\tau_k}\tau_i,\ldots,\tau_n)=\\
=D_{\tau_k}D_{\tau_k}(\omega(\tau_1,\ldots,\tau_n))-2\sum_{i,k} D_{\tau_k}\bra{\omega(\tau_1,\ldots,\nabla_{\tau_k}^\Sigma\tau_i,\ldots,\tau_n)}+\\
+\sum_{i,j,k}\omega\bra{\tau_1,\ldots,\nabla_{\tau_k}^\Sigma\tau_i,\ldots,\nabla_{\tau_k}^\Sigma\tau_j,\ldots,\tau_n)}= a+b+c,
\end{multline}
where $a=\Delta_\Sigma(\omega(\tau_1,\ldots,\tau_n))$ because $\{\tau_k\}$ is an orthonormal basis of the tangent space,
$b=0$ because $\dual{\nabla_{\tau_k}^\Sigma\tau_i,\tau_i}=\frac{1}{2}D_{\tau_k}\dual{\tau_i,\tau_i}=0$ and $\omega$ is alternating. Finally, also $c=0$:
\begin{multline}
\sum_{i,j,k}\omega\bra{\tau_1,\ldots,\nabla_{\tau_k}^\Sigma\tau_i,\ldots,\nabla_{\tau_k}^\Sigma\tau_j,\ldots,\tau_n)}=\\
=-\sum_{i,j,k\atop i\neq j}\dual{\nabla^\Sigma_{\tau_k}\tau_i,\tau_j} \dual{\tau_i,\nabla^\Sigma_{\tau_k}\tau_j}\omega_{1\cdots n}-
\sum_{j,k}\dual{\nabla^\Sigma_{\tau_k}\tau_j,\nabla^\Sigma_{\tau_k}\tau_j}\omega_{1\cdots n}=\\
=\bra{\sum_{i,j,k\atop i\neq j}\dual{\nabla^\Sigma_{\tau_k}\tau_i,\tau_j}^2-\sum_{j,k}\dual{\nabla^\Sigma_{\tau_k}\tau_j,
\nabla^\Sigma_{\tau_k}\tau_j}}\omega_{1\cdots n}=0.
\end{multline}
The last equality is justified by the fact that for $i=j$ fixed
$$\sum_{j\atop i\neq j}\dual{\nabla^\Sigma_{\tau_k}\tau_i,\tau_j}^2=\abs{\nabla^\Sigma_{\tau_k}\tau_i}^2=
\dual{\nabla^\Sigma_{\tau_k}\tau_i,\nabla^\Sigma_{\tau_k}\tau_i}.$$
Summing over $i$ and $k$ we conclude.
\end{proof}

\begin{lemma}[Codazzi's equation]
Let $h^\alpha_{ij}=\big(\nabla_{\tau_i}\tau_j\big)\cdot\nu_\alpha$ and $h^\alpha=H\cdot\nu_\alpha$ be the coefficients in local coordinates of the second fundamental form and of the mean curvature, respectively:
$$h(X,Y)=h^\alpha_{ij}X^i Y^j\nu_\alpha,\quad H=h^\alpha \nu_\alpha.$$
Then
\begin{equation}\label{eqcodazzi}
h^\alpha_{ik,k}=h^\alpha_{,i}
\end{equation}
where commas denote the covariant derivatives.
\end{lemma}
\begin{proof}
The connection of $\R{n+m}$ is flat, meaning that the curvature vanishes, therefore
\begin{multline}\nonumber
h^\alpha_{ik,k}= D_{\tau_k}\dual{\nabla_{\tau_i}\tau_k,\nu_\alpha}=
\dual{\nabla_{\tau_k}(\nabla_{\tau_i}\tau_k) ,\nu_\alpha}+\dual{\nabla_{\tau_i}\tau_k,\nabla_{\tau_k}\nu_\alpha}=\\
=\dual{\nabla_{\tau_i}(\nabla_{\tau_k}\tau_k) ,\nu_\alpha}+\dual{\nabla_{\tau_k}\tau_k,\nabla_{\tau_i}\nu_\alpha}=
D_{\tau_i}\dual{H,\nu_\alpha}.
\end{multline}

\end{proof}

\noindent\textbf{Notation} In what follows we will write
$$\omega_{1\cdots \alpha^i \cdots\beta^j\cdots n}:=\omega_{1\cdots (i-1) \alpha (i+1) \cdots (j-1)\beta (j+1)\cdots n}$$
to denote that $\alpha$ occurs in the $i$-th place and $\beta$ in the $j$-th.
\begin{prop}\label{evolomega} Along the mean curvature flow $\omega$ satisfies the following equation:
\begin{equation}\label{evolomega2}
\bra{\frac{\partial}{\partial t}-\Delta_\Sigma}\omega_{1\cdots n}=\omega_{1\cdots n}
\sum_{\alpha,i,k}(h^\alpha_{i k})^2
-\sum_{i,j, k,\alpha,\beta\atop i\neq j}\omega_{1\cdots \alpha^i \cdots\beta^j\cdots n}h^\alpha_{ik}h^\beta_{jk},
\end{equation}
where in the last sum $\alpha$ occupies the $i$-th place and $\beta$ the $j$-th.
\end{prop}

\begin{proof}
Being constant, $\omega$ is parallel on $\R{n+m}$, that is $\nabla\omega=0$. Thus 
$$\bra{\nabla_{\tau_k}^\Sigma\omega}_{1\cdots n}=\bra{(\nabla_{\tau_k}^\Sigma-\nabla_{\tau_k})\omega}_{1\cdots n}
= \sum_i \omega(\tau_1,\ldots,\nabla_{\tau_k}\tau_i-\nabla_{\tau_k}^\Sigma \tau_i,\ldots,\tau_n).$$
Observing that $\nabla_{\tau_k}^N\tau_i=\sum_\alpha h^\alpha_{ik}\nu_\alpha$ e $\omega_{1\cdots n,k}:=
\nabla_{\tau_k}^\Sigma\omega(\tau_1,\ldots,\tau_n)$ 
\begin{equation}\label{gradomega}
\omega_{1\cdots n,k}=\sum_{i,\alpha} \omega_{1\cdots  \alpha^i \cdots n}h^\alpha_{ik}.
\end{equation}
Similarly
$$\omega_{1\cdots \alpha^i\cdots n,k}=-\sum_{l}\omega_{1\cdots l^i\cdots n}h^\alpha_{lk}+
\sum_{\beta,j\atop i\neq j}\omega_{1\cdots \beta^j \cdots\alpha^i \cdots n}h^\beta_{jk}.$$
\begin{equation}\label{contiomega}
\omega_{1\cdots n,kk}=\sum_{\alpha,i} \omega_{1\cdots  \alpha^i \cdots n,k}h^\alpha_{ik}+
\sum_{\alpha,i} \omega_{1\cdots  \alpha^i\cdots n}h^\alpha_{ik,k}.
\end{equation}
We may apply Codazzi's equation \eqref{eqcodazzi}, obtain $h^\alpha_{i k,k}=h^\alpha_{,i}$ and by
\eqref{contiomega} we get:

\begin{multline}\label{lapl}
\omega_{1\cdots n,kk}=\\=-\sum_{\alpha,i,l,k}\omega_{1\cdots l^i\cdots n}h^\alpha_{lk}h^\alpha_{ik}+
\sum_{i\neq j}\omega_{1\cdots \beta^j \cdots\alpha^i \cdots n}h^\beta_{jk}h^\alpha_{i k}+
\sum_{\alpha,i} \omega_{1\cdots  \alpha^i\cdots n}h^\alpha_{,i}=\\
=-\omega_{1\cdots n}\sum_{i,k,\alpha}(h^\alpha_{ik})^2+
\sum_{i\neq j}\omega_{1\cdots \beta^j \cdots\alpha^i \cdots n}h^\beta_{jk}h^\alpha_{ik}
+ \sum_{\alpha,i}\omega_{1\cdots  \alpha^i\cdots n}h^\alpha_{,i}
\end{multline}
To compute $\frac{\partial}{\partial t}\omega$ we fix a time $t>0$, a point $p\in\Sigma_t$ and consider
$F:\Omega_T\rightarrow \R{n+m}$ a parametrization of the mean curvature flow satisfying
$\frac{\partial F}{\partial t}\in N\Sigma$ and such that
$\{\partial_1,\ldots,\partial_n\}_p$ (we intend $\partial_i:=\frac{\partial F}{\partial x^i}$) is an orthonormal basis of $T_p\Sigma$ which evolves remaining an orthonormal basis, say for all the times in $(t-\varepsilon,t+\varepsilon)$. This may always be done with a local reparametrization of $\Omega$ based on the inverse function theorem.

With this choices and with $g_{ij}:=\frac{\partial F}{\partial x^i}\cdot\frac{\partial F}{\partial x^j}$ we have that in $p$ 
\begin{equation}\label{baseon}
0=\frac{\partial}{\partial t}g_{ij}=\frac{\partial}{\partial t}\dual{(\nabla_{\partial_i}H)^T,\partial_j}
+\frac{\partial}{\partial t}\dual{\partial_i,(\nabla_{\partial_j}H)^T}.
\end{equation}
In $(t-\varepsilon,t+\varepsilon)\times\{p\}$ we have $\frac{\partial}{\partial t}g_{ij}=0$, $g=\sqrt{\det g}=1$;
moreover $\frac{\partial}{\partial t}\partial_i=\frac{\partial^2 F}{\partial t\partial x^i}=
\frac{\partial}{\partial x^i}H=\nabla_{\partial_i}H$ thus
\begin{multline}
\frac{\partial}{\partial t}\omega(\tau_1,\ldots,\tau_n)=\frac{\partial}{\partial t}
\bra{\frac{1}{g}\omega(\partial_1,\ldots,\partial_n)}=\frac{\partial}{\partial t}\omega(\partial_1,\ldots,\partial_n)=\\
=\sum_i\omega\bra{\partial_1,\ldots,\partial_{i-1},(\nabla_{\partial_j}H)^N,\partial_{i+1},\ldots,\partial_n)}+\\
+\sum_i\omega\bra{\partial_1,\ldots,\partial_{i-1},(\nabla_{\partial_j}H)^T,\partial_{i+1},\ldots,\partial_n)}.
\end{multline}
The last sum vanishes because setting $i=j$ in \eqref{baseon} we obtain that $(\nabla_{\partial_j}H)^T$ has no component along $\partial_j$ and the other components are unimportant because $\omega$ is alternating. Therefore we have
\begin{equation}\label{tempo}
\frac{\partial}{\partial t}\omega_{1\cdots n}=\sum_i\omega_{1\cdots  \alpha^i\cdots n}h_{\alpha, i}.
\end{equation}
Combining \eqref{tempo} with \eqref{lapl} and applying lemma \ref{lemmalapl} the proof follows.
\end{proof}

We recall the singular value decomposition.

\begin{lemma}\label{singdec} Given a linear application $L:\R{n}\rightarrow\R{m}$ there exist orthonormal basis $\{v_i\}_{i=1,\ldots,n}$ and $\{w_\alpha\}_{\alpha=1,\ldots,m}$ of $\R{n}$ and $\R{m}$ respectively such that the matrix
$\lambda_{i\alpha}$ associated to $L$ in such basis is diagonal, i.e. $\lambda_{i\alpha}=0$ if $i\neq\alpha$.
\end{lemma}

Apply this to $Du_t(x):\R{n}\rightarrow\R{m}$ fixing $t>0$, $x_0\in\Omega$ to find orthonormal basis
$\{v_i\}_{i=1,\ldots,n}$ and $\{w_\alpha\}_{\alpha=1,\ldots,m}$ as in lemma \ref{singdec}. To such basis we associate a basis of the tangent space and a basis of the normal space to the graph of $u(t,\cdot)$ in $(x_0,u(t,x_0)):$

$$\left\{ \tau_i=\frac{1}{\sqrt{1+\sum_\beta\lambda^2_{i\beta}}}(v_i+\sum_\beta\lambda_{i\beta}w_\beta)
\right\}_{i=1,\ldots,n}$$
$$\left\{ \nu_\alpha=\frac{1}{\sqrt{1+\sum_j\lambda^2_{j\alpha}}}(w_\alpha-\sum_j\lambda_{j\alpha}v_j)
\right\}_{\alpha=1,\ldots,m}$$
Observe that, defined $\pi$ the projection of $\R{n+m}$ on the first $n$ coordinates, we have
\begin{equation}\label{pi}
\pi(\nu_\alpha)=-\sum_j\lambda_{j\alpha}\pi(\tau_j).
\end{equation}
Since $\omega (a_1,\ldots,a_n)=\omega (\pi(a_1),\ldots,\pi(a_n))$, we may use \eqref{pi} to compute
\begin{equation}\label{sempl}
\omega_{1\cdots \beta^j \cdots\alpha^i \cdots n}=\omega_{1\cdots n}(-\lambda_{\beta j}\lambda_{\alpha i}
+\lambda_{\beta i}\lambda_{\alpha j})
\end{equation}
Now proposition \ref{evolomega} may be written in terms of the singular values of $Du$:

\begin{multline}\label{evolomega3}
\bra{\frac{\partial}{\partial t}-\Delta_\Sigma}*\omega=*\omega
\Big(\sum_{\alpha,i,k} (h^\alpha_{ik})^2
+\sum_{i,j, \alpha,\beta k\atop i\neq j}(-\lambda_{\beta j}\lambda_{\alpha i}
+\lambda_{\beta i}\lambda_{\alpha j})h^\alpha_{i k}h^\beta_{j k}\Big)=\\
=*\omega\Big(\sum_{\alpha,i,k} (h^\alpha_{ik})^2
+\sum_{i,j, k\atop i\neq j}(-\lambda_j\lambda_i h^i_{ik}h^j_{jk}+\lambda_j\lambda_i h^i_{jk} h^j_{ik}\Big)
\end{multline}

\begin{prop}\label{stimona} Let $*\omega\geq\frac{1}{\sqrt{2-\delta}}.$
Then
\begin{equation}\label{stimona2}
\bra{\frac{\partial}{\partial t}-\Delta_\Sigma}*\omega\geq\delta|A|^2,
\end{equation} 
where $|A|^2=\sum_{\alpha,i,j}|h^\alpha_{ij}|^2$ is the squared norm of the second fundamental form.
\end{prop}

\begin{proof}
We shall use the fact that the hypothesis on $*\omega$ implies $\sum_i \lambda_i^2\leq 1-\delta$ and $0\leq\lambda_i\lambda_j
\leq 1- \delta,$ for instance following \eqref{relaz}. 

\textbf{1.} Assume $n\leq m$.
$$\sum_{\alpha,i,k}(h^\alpha_{i k})^2=\sum_{\alpha,i,k\atop \alpha>n} (h^\alpha_{i k})^2 + \sum_{i,k}(h^i_{ik})^2
+\sum_{i,j,k\atop i\neq j} (h^i_{jk})^2.$$
We estimate the right hand side of \eqref{evolomega3}, which we simplify because $\lambda_{i\alpha}$ is diagonal:
\begin{multline}
\sum_{\alpha,i,k}(h^\alpha_{i k})^2
-\sum_{i,j,k\atop i\neq j}\lambda_i\lambda_j h^i_{ik}h^j_{jk}
+\sum_{i,j,k\atop i\neq j}\lambda_i\lambda_j h^i_{jk}h^j_{ik}\geq\\
\geq\delta |A|^2+ (1-\delta)\sum_{i,k}(h^i_{ik})^2 +(1-\delta)\sum_{i,j,k\atop i\neq j} (h^i_{jk})^2+\\
-\sum_{i,j,k\atop i\neq j}\lambda_i\lambda_j h^i_{ik}h^j_{jk}-(1-\delta)\sum_{i,j,k\atop i\neq j}|h^i_{jk}h^j_{ik}|\geq \delta|A|^2
+ \sum \lambda_i^2\sum_{i,k}(h^i_{ik})^2 +\\+(1-\delta)\sum_{i,j,k \atop i\neq j} (h^i_{jk})^2
-\sum_{i,j,k\atop i\neq j}\lambda_i\lambda_j|h^i_{ik}h^j_{jk}|-(1-\delta)\sum_{i,j,k\atop i\neq j}|h^i_{jk}h^j_{ik}|\geq\\
\geq\delta|A|^2+(1-\delta) \sum_{i,j,k \atop i\neq j}(h^i_{jk})^2-|h^j_{ik}h^i_{jk}| )
+\Big(\sum_{i,k}h^i_{ik}\lambda_i \Big)^2-\sum\lambda_i\lambda_j h^i_{ik}h^j_{jk}\geq\\
\geq \delta|A|^2 +\sum_{i,j,k\atop i\neq j}(|h^i_{jk}|-|h^j_{ik}|)^2+\Big(\sum_{i,k}\lambda_i h^i_{ik}\Big)^2\geq
\delta|A|^2.
\end{multline}

\textbf{2.} The case $m\leq n$ may be deduced by the case $m=n$ and thus by the above step observing that in
$$\sum_{\alpha,i,k}(h^\alpha_{i k})^2 -\sum_{i,j,k\atop i\neq j}\lambda_i\lambda_j h^i_{ik}h^j_{jk}
+\sum_{i,j,k\atop i\neq j}\lambda_i\lambda_j h^i_{jk}h^j_{ik}$$
the second and third terms are 0 when $i>m$ or $j>m$, while in the first term we may neglect the terms corresponding to  $i>m$ or $j>m$ because they are positive.
\end{proof}

\begin{trm}\label{stimainterna} Assume that the initial data $\psi$ satisfies
\begin{equation}\label{datoiniz}
8n\diam\Omega\sup_\Omega\abs{D^2\psi}+\sqrt{2}\sup_{\partial\Omega}\abs{D\psi}< \sqrt{2^{1/n}-1}.
\end{equation}
Then there exists $\delta\in(0,1)$ such that the solutions to the mean curvature flow \eqref{mcs} satisfy
\begin{equation}\label{stimagrad}
\sup_{\Omega_T}\abs{Du}^2 < 1-\delta.
\end{equation}
\end{trm}

\begin{proof}
The hypothesis implies
$\frac{1}{\sqrt {\det (I + D\psi^TD\psi)}}>\frac{1}{\sqrt{2}},$
that is $*\omega>\frac{1}{\sqrt{2}}$ at the time $t=0$.
By continuity of $*\omega$ and compactness of $\Cl{\Omega}$ we may find $\delta\in(0,1)$ such that $*\omega> \frac{1}{\sqrt{2-\delta}}$
at the time $t=0$ and
$$8n\diam\Omega\sup_\Omega\abs{D^2\psi}+\sqrt{2}\sup_{\partial\Omega}\abs{D\psi}< \sqrt{(2-\delta)^{1/n}-1}.$$
Thus, by proposition \ref{stimona}, as long as the condition $*\omega> \frac{1}{\sqrt{2-\delta}}$ holds true, we have
\begin{equation}\label{evol}
\bra{\frac{\partial}{\partial t}-\Delta_\Sigma}*\omega\geq\delta|A|^2.
\end{equation}
where $A$ is the second fundamental form.
Now assume that there is a first moment $t_0$ such that for some $x_0\in\Omega$ we have $*\omega=\frac{1}{\sqrt{2-\delta}}.$ Then thanks to the boundary gradient estimates \eqref{stimabordo}, and since thanks to
\eqref{relaz} $\sup_{\Omega\times [0,t_0)}\abs{Du}^2\leq 1-\delta$, we get
$$\sup_{\partial\Omega\times [0,t_0)}\abs{Du}< 8n\diam\Omega\sup_\Omega\abs{D^2\psi}+\sqrt{2}\sup_{\partial\Omega}\abs{D\psi}< \sqrt{(2-\delta)^{1/n}-1},$$
which, due to \eqref{relaz}, implies that $*\omega>\frac{1}{\sqrt{2-\delta}}$ on the boundary. Therefore $x_0$ is an interior point where $*\omega$ attains a minimum smaller then the infimum on the boundary. This is absurd because the parabolic maximum principle applies to \eqref{evol}. Thus for every time (in $[0,T)$) it's true that $*\omega>\frac{1}{\sqrt{2-\delta}}$ and, consequently, $\abs{Du}^2<1-\delta$.
\end{proof}

\section{Long time existence of the mean curvature flow}
The a priori estimates of $\sup_{\Omega_T}|u|+\sup_{\Omega_T}|Du|$, in codimension 1, yield also the a priori estimates on the higher order derivatives (proposition \ref{proplady}), thanks to De Giorgi's theorem. In arbitrary codimension the a priori estimates on the higher order derivatives are in general not available: if we could obtain a $C^{1,\alpha}$ estimate depending only on $\sup|Du|$ and, thus, on the $L^{\infty}$-norm of the coefficients $g^{ij}$, we could also prove the smoothness of the Lipschitz solutions of the minimal surface system. On the other hand, this latter result is false because of the minimal cone exhibited by Lawson and Osserman, theorem \ref{controregol}.

To obtain the higher order estimates we will use a theorem of Brian White \cite{white} and we will prove them in the area decreasing case.

\medskip

The work done to prove long time existence may be divided into the following steps:

\begin{enumerate}
\item we use theorem \ref{trmtempipiccoli} to prove existence for small times;
\item assuming that we have proved existence in $[0,t_0)$, we study the possibility that a singularity appears at the time $t_0$ and, thanks to the gradient a priori estimates and the area-decreasing condition, we obtain that the second fundamental forms of the surfaces $\Sigma_t,$ $t<t_0$
vanish on a suitable parabolic blow-up;
\item from the latter step we may apply White's theorem, obtaining the $C^{2,\alpha}$ a priori estimates;
\item the solution converges with its derivatives as $t\rightarrow t_0^-$, therefore we may reapply the short time existence theorem and conclude that there is existence for alla times.
\end{enumerate}

\subsection{The parabolic blow-up, the Gaussian density and White's theorem}
We may consider a mean curvature flow $F$ in $\R{n+m}$ as a subset of $\R{n+m}\times[0,t_0)$, simply associating to
$$F:\Cl\Omega\times[0,t_0)\rightarrow\R{n+m}$$
its trace in the space-time:
\begin{equation}\label{eqtraccia}
\mathcal{M}=\{(F(x,t),t):x\in\Cl\Omega,t\in [0,t_0)\}.
\end{equation}

\begin{defin}[Gaussian density] Let $\mathcal{M}$ be a mean curvature flow as in \eqref{eqtraccia}. Then the Gaussian density of $\mathcal{M}$ in $X=(x,t)$ with radius $r$ is
$$\Theta(\mathcal{M},X,r):=\int\limits_{y\in\mathcal{M}(t-r^2)}\frac{1}{(4\pi r^2)^{\frac{n}{2}}}
e^{\frac{-|y-x|^2}{4r^2}}d\mathcal{H}^n(y),$$
where we intend that $\M{M}(t-r^2)=\M{M}\cap\{X=(y,s):s=t-r^2\}.$
\end{defin}
We may prove, in analogy with the monotonicity formula for minimal surfaces, that the quantity $\Theta(\M{M},X,r)$ is monotonous with respect to $r$, so that its limit exists and we may define
\begin{equation}\label{eqdenspar}
\Theta(\M{M},X):=\lim_{r\rightarrow 0}\Theta(\M{M},X,r).
\end{equation}

\begin{trm}[White]\label{trmwhite} For every $0<\alpha<1$ there exist $\varepsilon=\varepsilon(n,m,\alpha)>0$ and $C=C(n,m,\alpha)>0$ such that if $\M{M}$ is the mean curvature flow of an $n$-dimensional graph in $\R{n+m}$ and if for a certain open set $U\in\R{n+m}\times [0,t)$ and every $X\in U$ and $0<r<\dist(x,U^c)$ the following holds
$$\Theta(\M{M},X,r)\leq 1+\varepsilon,$$
then 
$$\norm{u\big|_U}_{2,\alpha}\leq C,$$
where $u$ is the function whose graph parametrizes the mean curvature flow and the H\"older norm is the parabolic one.
\end{trm}

\begin{oss} This theorem, to be compared with Allard's theorem, which is the elliptic analogue, says that we get \emph{local} estimates on the higher order derivatives if we may control the Gaussian density (which is the analogue of the density defined in \eqref{eqdens}).
\end{oss}

To apply this theorem we need the notion of parabolic blow-up. Similarly to the blow-up defined in proposition \ref{blowup}, the parabolic blow-up is a dilatation of the space-time performed in such a way that the mean curvature flow system is preserved. To obtain this we need to treat the time variable differently from the space variable, as done in the introduction to the parabolic H\"older norms.

\begin{defin} The parabolic blow-up of a space-time $\R{n+m}\times[0,t_0)$ in $(y_0,t_0)$ with parameter $\lambda$ is the bijection
$$D_\lambda:\R{n+m}\times[0,t_0)\rightarrow \R{n+m}\times [-\lambda^2t_0,0)$$
defined by
\begin{equation}\label{eqbupar}
D_\lambda(y,t)=(\lambda(y-y_0),\lambda^2(t-t_0)).
\end{equation}
\end{defin}

Studying the density in a point of a mean curvature flow as defined in \eqref{eqdenspar} is equivalent to studying 
$$\lim_{\lambda\rightarrow+\infty}\Theta(D_\lambda \M{M},0,r).$$

\subsection{The long time existence theorem}

\begin{trm}
Let $\psi\in C^{\infty}(\Cl\Omega)$ be satisfying \eqref{datoiniz}. Then the mean curvature flow, solution to \eqref{mcs}, exists in $C^\infty(\Omega_\infty)\cap C^1(\Cl{\Omega}_\infty)$; moreover there exists $\delta>0$ such that $|Du(x,t)|\leq 1-\delta$ for every $(x,t)\in\Cl{\Omega}_\infty$.
\end{trm}

\begin{proof} We proceed in several steps.

\textbf{1.} Thanks to theorem \ref{trmtempipiccoli} there are $\varepsilon>0$ and a solution $u\in C^{2,\alpha}_{(\delta)}(\Omega_{\varepsilon})$ to the system \eqref{mcs}.

\textbf{2.} Thanks to the remark following theorem \ref{esistparab}, the solution found belongs to $C^1(\Cl\Omega_\varepsilon)$ meaning that for every $t\in(0,\varepsilon)$, $u(t,\cdot)\in C^{1}(\Cl\Omega)$. Consequently we may apply the interior and boundary gradient estimates given by theorem \ref{stimainterna} and conclude that $|Du|\leq 1-\delta$.

\textbf{3.} The set of times for which there is a solution is closed: if $t_0<+\infty$ is the $\sup$ of the times for which a solution $u$ exists, then the limits of the space derivatives and of the time derivative exist as $t\rightarrow t_0^-$.
To prove this we study the possibility that a singularity appears in $(y_0,t_0)$, $y_0\in\Omega$.

We use the backward heat kernel\footnote{it's called \emph{backward} heat kernel because compared with the standard heat kernel we have $t<t_0$ and instead of $t-t_0$ we have $t_0-t$; it is used to study the mean curvature flow \emph{before} a certain time.} in $(y_0,t_0)$, introduced by Huisken in \cite{huisken}:

$$\rho_{y_0,t_0}(y,t):=\frac{1}{(4\pi(t_0-t))^{\frac{n}{2}}}\exp\Big(\frac{-|y-y_0|^2}{4(t_0-t)}\Big);$$

By a monotonicity formula, Huisken proves that $\lim_{t\rightarrow t_0}\int \rho_{y_0,t_0} d\mu_t$ exists, where
$\mu_t=\M{H}^n\res\Sigma_t$ and $\Sigma_t$ is the surface moving by mean curvature. Moreover the backward heat kernel satisfies the following equation, proved by M-T. Wang in \cite{wangfour}
$$\frac{d}{dt}\rho_{y_0,t_0}=-\Delta_{\Sigma_t}\rho_{y_0,t_0}-\rho_{y_0,t_0}\Big(\frac{|F^N|^2}{4(t-t_0)^2}+
\frac{F^N\cdot H}{t-t_0}\Big),$$
where $F^N(x,t)$ is the projection of $F(x,t)$ into $N_{F(x,t)}\Sigma_t$. Recalling that the mean curvature flow satisfies
$$\frac{d}{dt}d\mu_t=-|H|^2d\mu_t$$
and using proposition \ref{stimona} to get
$$\frac{d}{dt}*\omega\geq\Delta_{\Sigma_t}*\omega+\delta|A|^2,$$
we obtain

\begin{multline}\label{eqlungotempo}
\frac{d}{dt}\int (1-*\omega)\rho_{y_0,t_0}d\mu_t\leq \int[\Delta_{\Sigma_t}(1-*\omega)-\delta|A|^2]\rho_{y_0,t_0}d\mu_t+\\
-\int (1-*\omega)\Big[\Delta_{\Sigma_t}\rho_{y_0,t_0}+\rho_{y_0,t_0}\Big(\frac{|F^N|^2}{4(t-t_0)^2}+
\frac{F^N\cdot H}{t-t_0}\Big)\Big]d\mu_t-\\-\int(1-*\omega)|H|^2\rho_{y_0,t_0}d\mu_t.
\end{multline}

Every integral is intended over the space-time and respect to the measure $\mu_t=\M{H}^n\res\Sigma_t$. Reordering the terms on the right hand side of the above inequality we get

\begin{multline}
\int\big[\Delta_{\Sigma_t}(1-*\omega)\rho_{y_0,t_0}-(1-*\omega)\Delta_{\Sigma_t}\rho_{y_0,t_0}\big]d\mu_t
-\delta\int|A|^2\rho_{y_0,t_0}d\mu_t+\\
-\int(1-*\omega)\rho_{y_0,t_0}\Big[\frac{|F^N|^2}{4(t-t_0)^2}+
\frac{F^N\cdot H}{t-t_0}+|H|^2\Big]d\mu_t.
\end{multline}

The first term vanishes integrating by parts and the third one completes to squares: \eqref{eqlungotempo} becomes

\begin{multline}
\frac{d}{dt}\int (1-*\omega)\rho_{y_0,t_0}d\mu_t\leq\\
\leq -\delta\int|A|^2\rho_{y_0,t_0}d\mu_t-\int(1-*\omega)\rho_{y_0,t_0}
\Big|\frac{|F^N|}{2(t_0-t)}+H\Big|^2d\mu_t
\end{multline}

Being $*\omega>0$ and $\rho_{y_0,t_0}>0$ it's clear that
$$\int (1-*\omega)\rho_{y_0,t_0}d\mu_t\leq \int\rho_{y_0,t_0}d\mu_t;$$
being the last integral finite we find that
$$\frac{d}{dt}\int (1-*\omega)\rho_{y_0,t_0}d\mu_t\leq C-\delta\int|A|^2\rho_{y_0,t_0}d\mu_t$$
for some constant $C>0$

For $\lambda>1$ we apply a parabolic dilatation $D_\lambda$ in $(y_0,t_0)$ as defined in \eqref{eqbupar}. If $\M{M}$ is the trace of the mean curvature flow, we now study $\M{M}^\lambda:=D_\lambda(\M{M});$ we denote the new time parameter by $s$, so that $t=t_0-\frac{s}{\lambda^2}$, and the volume form induced after the blow-up is
$d\mu^\lambda_s:$ it is the volume form on the surface
$$\Sigma^\lambda_s=F^\lambda_s(\Cl\Omega):=\lambda F_{t_0+\frac{s}{\lambda^2}}(\Cl\Omega).$$

By means of a change of variable
$$\frac{d}{ds}\int(1-*\omega)\rho_{0,0}d\mu^\lambda_s=\frac{1}{\lambda^2}\frac{d}{dt}\int(1-*\omega)\rho_{y_0,t_0}d\mu_t
\leq \frac{C}{\lambda^2}-\frac{\delta}{\lambda^2}\int\rho_{y_0,t_0}|A|^2d\mu_t.$$
Observe that $*\omega$ is invariant under the parabolic blow-up (which is a homotethy in the space variables and, thus, doesn't alter the differentials). Using the fact that the second fundamental form $A$ rescales as $\frac{1}{\lambda}$ (because it's obtained from the second derivatives) and that $\rho_{y_0,t_0}d\mu_t$ is invariant under parabolic blow-up, we obtain
$$\frac{1}{\lambda^2}\int\rho_{y_0,t_0}|A|^2d\mu_t=\int\rho_{0,0}|A|^2d\mu^\lambda_s.$$
Therefore
$$\frac{d}{ds}\int(1-*\omega)\rho_{0,0}d\mu^\lambda_s\leq \frac{C}{\lambda^2}-\delta\int\rho_{0,0}|A|^2d\mu^\lambda_s.$$
Integrating with respect to $s$ from $-1-\tau$ to $-1$ for some $\tau>0$ we get
\begin{multline}\label{eqcontolambda}
\delta\int_{-1-\tau}^{-1}\int\rho_{0,0}|A|^2d\mu^\lambda_sds\leq\\
\leq\int(1-*\omega)\rho_{y_0,t_0}d\mu_{-1}^\lambda-\int(1-*\omega)\rho_{y_0,t_0}d\mu_{-1-\tau}^\lambda +\frac{C}{\lambda^2}.
\end{multline}
Observe that
$$\int(1-*\omega)\rho_{0,0}d\mu^\lambda_{s}=\int(1-*\omega)\rho_{y_0,t_0}d\mu_{t_0+\frac{s}{\lambda^2}}$$
and, since $\lim_{t\rightarrow t_0^-}\int(1-*\omega)\rho_{y_0,t_0}d\mu_t$ exists, we conclude that the right hand side of \eqref{eqcontolambda} goes to zero as $\lambda\rightarrow+\infty.$ For every $\tau>0$ we may thus choose a sequence $\lambda_j\rightarrow+\infty$ such that
$$\int_{-1-\tau}^{-1}\int\rho_{0,0}|A|^2d\mu^{\lambda_j}_s\leq C(j)$$
with $C(j)\rightarrow 0$. But $\tau$ is arbitrary and, up to choose each $C(j)$ smaller, we may find a sequence  $\tau_j$ such that $\frac{C(j)}{\tau_j}\rightarrow 0$ and a sequence $s_j\in[-1-\tau_j,1]$ such that
\begin{equation}\label{eqinvestiga}
\int_{-1-\tau}^{-1}\int\rho_{0,0}|A|^2d\mu^{\lambda_j}_{s_j}\leq \frac{C(j)}{\tau_j}.
\end{equation}
To study \eqref{eqinvestiga} we observe that
$$\rho_{0,0}(F_{s_j}^{\lambda_j})(F^{\lambda_j}_{s_j})=\frac{1}{4\pi(- s_j)}\exp\Big(\frac{-|F_{s_j}^{\lambda_j}|^2}{4(-s_j)} \Big),$$
where $F_{s_j}^{\lambda_j}=\lambda_j F_{t_0+\frac{s_j}{\lambda_j^2}}.$

If we consider, for every $R>0$, the ball $B_R(0)\subset\R{n+m}$ and, for $j$ large enough, assume $-1<s_j<-\frac{1}{2}$, then
$$\int\rho_{0,0}|A|^2d\mu^{\lambda_j}_{s_j}\geq\frac{1}{2\pi}e^{\frac{-R^2}{2}}
\int_{\Sigma^{\lambda_j}_{s_j}\cap B_R(0)}|A|^2d\mu^{\lambda_j}_{s_j}.$$
From that follows that for every compact $K\subset \R{n+m}$ we have
\begin{equation}\label{eqasvanisce}
\int_{\Sigma^{\lambda_j}_{s_j}\cap K}|A|^2d\mu^{\lambda_j}_{s_j}\rightarrow 0,\quad j\rightarrow +\infty.
\end{equation}

We now want to prove that \eqref{eqasvanisce} and $*\omega$ being bounded from below (a priori gradient estimates) imply
\begin{equation}\label{eqlimiterho}
\lim_{j\rightarrow +\infty}\int\rho_{y_0,t_0}d\mu_{t_0+\frac{s_j}{\lambda_j^2}}=\lim_{j\rightarrow+\infty}\int\rho_{0,0}
d\mu_{s_j}^{\lambda_j}\leq 1.
\end{equation}
Assume that the origin is a limit point for $\Sigma_{s_j}^{\lambda_j},$ otherwise there is nothing to prove. Thanks to the a priori gradient estimates, each $\Sigma_t$ is the graph of a function $u_t$ with gradient equiboundedwith respect to $t$. Set
$$u_j:=u_{t_0+\frac{s_j}{\lambda_j^2}}$$
and perform an elliptic blow-up (see proposition \ref{blowup}) of the graph of $u_j$ with parameter $\lambda_j$. The surface that we obtain is the graph of a surface which we denote by
$$\widetilde u_j:\lambda_j \Cl\Omega\rightarrow\R{m}.$$
Of course also $|D\widetilde u_j|$ is equibounded and we assumed that $$\lim_{j\rightarrow +\infty}\widetilde u_j(0)=0.$$ The hypothesis on the gradients and corollary \ref{corascoli} to Ascoli-Arzel\`a's theorem imply that we may assume $\widetilde u_j\rightarrow \widetilde u_\infty$ in $C^{0,\alpha}$ on every compact set for some $0<\alpha<1$, $u_\infty:\R{n}\rightarrow\R{m}$ being Lipschitz. It may be proven, as in \cite{ilmanen} the following inequality:
$$|A_j|\leq|\nabla^{\Sigma_{s_j}^{\lambda_j}}d\widetilde u_j|\leq(1+|D\widetilde u_j|^2)^{\frac{3}{2}}|A_j|, $$
being $A_j$ the second fundamental form on $\Sigma_{s_j}^{\lambda_j}$.

It follows that
$$\widetilde u_j\rightarrow \widetilde u_\infty,\quad \mathrm{in}\; C^{0,\alpha}_{\loc}\cap W^{1,2}_{\loc}$$
and the second derivatives of $\widetilde u_\infty$ vanish. Then $\Sigma_{s_j}^{\lambda_j}\rightarrow\Sigma_{-1}^\infty$ in the sense of Radon measures in $\R{n+m}$ and, moreover, $\Sigma_{-1}^\infty$ is the graph of an affine map. Therefore
$$\lim_{j\rightarrow+\infty}\int\rho_{0,0}d\mu_{s_j}^{\lambda_j}=\int\rho_{0,0}d\mu_{-1}^\infty=1,$$
which implies
$$\lim_{j\rightarrow+ \infty}\int\rho_{y_0,t_0}d\mu_{t_0+\frac{s_j}{\lambda_j^2}}=\lim_{j\rightarrow+\infty}
\int\rho_{y_0,t_0}d\mu_t=1.$$
Now White's theorem, theorem \ref{trmwhite}, gives the local estimates in the parabolic $C^{2,\alpha}_{\loc}$-space.  These, thanks to Ascoli-Arzel\`a's theorem, imply the convergence of a subsequence together with the time derivative and the space derivatives. On the other hand the limit of any such a subsequence is uniquely determined, so that $u_t\xrightarrow{C^2} u_{t_0}$ and $u$ solves the mean curvature flow \eqref{mcs} also in $t_0$.

\textbf{4.} The set of times for which the solution exists is open: set $t_0>0$ such that the solution exists up to $t_0$ meaning that the first order time derivative and the space derivatives of order less or equal to two have limit as $t\rightarrow t_0^-$. Then the small time existence theorem applies to the limit $u(t_0,\cdot)$ and we have a solution in $C^{2,\alpha}_{(\delta)}(\Omega_{t_0+\varepsilon})$ which is a solution also for the time $t_0$ (i.e. the solution before $t_0$ glues well with the solution in $(t_0,t_0+\varepsilon)$).

Putting together step 3 and step 4 the long time existence follows.
\end{proof}

\section{Convergence of the mean curvature flow}

The mean curvature flow decreases the area:
let $H_t=\Delta F_t$ be the mean curvature vector, where $F_t$ is the immersion of $\Omega$ in $\R{n+m}$ given by the graph of $u(-,t).$ Then from \eqref{eqmoto} used with $F$ instead of $\widetilde F$ we get
$$\frac{d}{dt}\M{H}^n(\M{G}_{u_t})=-\int_{\M{G}_{u_t}} H_t \cdot \frac{\partial F}{\partial t}=-\int_\Omega \abs{H_t}^2\sqrt{g_t}dx,$$
where we have used
$$\bra{0,\ldots,0,g^{ij}(Du)u_{ij}}^N=\Delta_\Sigma F=H.$$
The variation of the area in finite time is obtained integrating with respect to the time:
$$\int_0^{t_0} \int_\Omega\abs{H_t}^2\sqrt{g_t}dxdt=A(0)-A(t_0)\leq A(0),$$
and so the integral on the left is finite. This implies that there exists a sequence of times $t_i\rightarrow\infty$ such that
\begin{equation}\label{curvaturazero}
\int_\Omega \abs{H_{t_i}}^2\sqrt{g_{t_i}}dx\rightarrow 0.
\end{equation}

\begin{trm}\label{convergenza}
Consider a sequence of Lipschitz equibounded maps $u_j:\Omega\rightarrow\R{m}$ such that $|Du_j|\leq 1-\delta$, for some $\delta>0$. Assume that the first variations, defined in \eqref{varprimaeq}, satisfy  $\norm{\delta \M{G}_{u_j}}\rightarrow 0$. Then there exists a subsequence $u_{j'}$ converging uniformly and in the sense of varifolds to a Lipschitz function $u$ with $|Du|\leq 1-\delta$ whose graph is minimal in the sense of varifolds. Moreover
$$\V(\M{G}_{u_{j'}},1)\rightarrow\V(\M{G}_u,1)$$
in the sense of varifolds.
\end{trm}

For the elementary notions and definitions in the theory of varifolds, see the appendix.

\begin{proof}

\textbf{1.} By the theorem of Ascoli-Arzel\`a there exists a subsequence $u_{j'}\rightarrow u$ uniformly.
We want to prove that the convergence is also in the sense of varifolds.

By Allard's compactness theorem, theorem \ref{compallard}, there is a subsubsequence $u_{j''}$ such that $U_{j''}\rightarrow V$
in the sense of varifolds, where $V=\V(\Sigma,\theta)$ is an integer multiplicity rectifiable varifold, while $U_{j''}=\V(\M{G}_{u_{j''}},1)$ is the graph of $u_{j''}$ seen as an i.m. rectifiable varifold. We only need to prove that $V$ is the varifold induced by the graph of $u$, i.e. $\Sigma=\M{G}_u$ and $\theta=1$ up to a set $\M{H}^n$-negligible. In this case the whole sequence $U_{n'}$ would converge to $V$ in the sense of varifolds.

\textbf{2.} We prove that $V=\V(\M{G}_u,1)$. Clearly $\spt U\subset \M{G}_u$: let indeed $A$ be an open set non intersecting the graph of $u$. This latter is closed, thus for every continuous function compactly supported in $A$, we have
$\dist(\spt f,\M{G}_u)=\varepsilon>0$. Since the convergence is uniform, we may choose $j_0$ such that for $j\geq j_0$
we have
$\norm{u_j(x)-u(x)}_\infty<\varepsilon$, thus
$$\int_{\M{G}_u}f(x)d\M{H}^n(x)\rightarrow 0.$$
Then $V(f)=0$ for every $f$ supported in $A$ and by the arbitraryness of $A$ we have that the support of $V$ is included in $\M{G}_u.$

We now prove that for $\M{H}^n$-almost every
$p\in\M{G}_u$ we have $\theta(p)=1$. $U_{j''}\rightarrow V$ in the sense of varifolds implies that
\begin{equation}\label{proiezvar}
\pi_{\#}U_{j''}\rightarrow\pi_{\#}V
\end{equation}
in the sense of varifolds, where $\pi:\R{n+m}\rightarrow\R{n+m}$ is the orthogonal projection on $\R{n}\times\{0\}$.
To prove \eqref{proiezvar} we use \eqref{varimg}:
\begin{multline}
\pi_{\#}U_{j''}(f)=\\
=\int_{G_n}f(\pi(x),d\pi_x S)J\pi(x,S)dU_{j''}\rightarrow \int_{G_n}f(\pi(x),d\pi_x S)J\pi(x,S)dV_{j''}=\\
=\pi_{\#}V(f),
\end{multline}
where $G_n=G_n(\Omega\times\R{m})$ is the Grassmann bundle on $\Omega\times\R{m}$, as defined in \ref{defgrasm}.

The limit \eqref{proiezvar} may be rewritten as
$$\M{H}^n\res\Omega\times\{0\}\rightarrow \theta\M{H}^n\res\Omega\times\{0\},$$
whence $\theta=1$.

\textbf{3.} $V$ is minimal because for a given vector field
$$X\in C^1_0(\Omega\times\R{m};\R{n+m})$$ we have
$$|U_i(\diver X)|=\Big|\int_{\M{G}_{u_i}} X\cdot H_i\Big|\leq \sup|X|\norm{\delta U_i} \rightarrow 0.$$
And by the varifold convergence, we get
$$U_i(\diver X)\rightarrow V(\diver X)=0,$$
thus the limit graph is minimal in the sense of varifolds.
\end{proof}

\begin{oss} The same proof applies to minimal graphs defined on all of $\R{n+m}$ and to sequences defined on set invading $\R{n+m}$ as in the blow-up case, proposition \ref{blowup} and following.
\end{oss}

\chapter{Regularity in arbitrary codimension}\label{capitoloreg}

To study the regularity of minimal graphs we will use a blow-up procedure; we know that the blow-up of the graph of a smooth function converges to a plane. Allard's theorem says that, in the case of minimal graphs, the converse is true: if the blow-up in a point $p$ of a minimal graph converges to a plane, then the graph is smooth in a neighborhood of $p$. This reduces the regularity problem to the classification of the objects arising as blow-ups of minimal graphs. Since such objects are entire minimal graphs, the result we need is a Bernstein-type theorem: entire minimal graph, under suitable assumptions, are planes.

\section{Blow-ups and blow-downs: minimal cones}

\begin{prop}[Blow-up]\label{blowup}
Let $u:\Omega\rightarrow\R{m}$ be a Lipschitz map, $|Du|\leq K$, with $\M{G}_u$ minimal in the sense of varifolds.
Let $u_\lambda$ be defined by
$$u_\lambda (x)= \frac{1}{\lambda}(u(\lambda x)-u(x_0)), \quad x_0\in\Omega.$$
Then there exists a sequence $\lambda(i)\rightarrow 0$ such that
$u_{\lambda(i)}\rightarrow v$
uniformly on compact sets and in the sense of varifolds, where the graph of $v:\R{n}\rightarrow\R{m}$
is a cone minimal in the sense of varifolds.
\end{prop}

\begin{proof}
The convergence of a sequence $u_{\lambda(i)}$ to a Lipschitz minimal graph is an immediate consequence of theorem \ref{convergenza}. We may assume, without loss of generality, that $x_0=0$ and $u(0)=0$.

We prove that $v(\tau x)=\tau v(x)$ for every $\tau > 0$.

\begin{eqnarray}
\abs{v(x)-\tau v\bra{\frac{x}{\tau}}}\leq \abs{v(x)-\tau u_\lambda\bra{\frac{x}{\tau}}}
+\abs{\tau u_\lambda\bra{\frac{x}{\tau}}-\tau v\bra{\frac{x}{\tau}}}=\\
=\abs{v(x)- u_{\tau\lambda}(x)}+\tau\abs{u_\lambda\bra{\frac{x}{\tau}}-v\bra{\frac{x}{\tau}}}.
\end{eqnarray}

Thanks to the convergence of $u_\lambda$, the last two terms go to zero, therefore $v(\tau x)=\tau v(x)$.
\end{proof}

\begin{prop}[Blow-up of a cone]\label{blowup2} Let $u:\R{n}\rightarrow\R{m}$ be a Lipschitz cone\footnote{a cone $C\subset\R{n+m}$ with vertex in the origin is a set such that for every $\lambda>0$ we have $\lambda C=C.$ }, minimal as varifold, with $|Du|\leq K$. Let $x_0\in\R{n}$ and define
$$u_\lambda (x)= \frac{1}{\lambda}(u(x_0+\lambda x)-u(x_0)).$$
Then there exists a sequence $\lambda(i)\rightarrow 0$ such that $u_{\lambda(i)}\rightarrow v,$ uniformly on compact sets and in the sense of varifolds, where $\M{G}_v$
is a minimal cone in the sense of varifolds.
Moreover $\M{G}_v$ is a product of the form $\mathbb{R}\times C$, where $C$ is a minimal cone of dimension $n-1$ in $\R{n+m-1}$.
\end{prop}

Let $\widetilde x = (x_2,\ldots, x_n)$. The last assertion means that there exists an orthonormal system of coordinates $\R{n}$, a function $\widetilde v:\R{n-1}\rightarrow\R{m}$ and $\sigma\in\R{}$ such that
\begin{equation}\label{conooriz}
v(x_1,\ldots,x_n)=\sigma x_1 +\widetilde{v}(\widetilde x)
\end{equation}
and $\M{G}_{\widetilde{v}}$ is a minimal cone.

\medskip

\begin{proof} We may apply proposition \ref{blowup} to $u$ and obtain a Lipschitz map $v$, uniform limit of $u_\lambda$, with $\M{G}_v$ minimal varifold. We want to prove that
$$v(x+\tau x_0)=v(x)+\tau v(x_0)=v(x)+\tau u(x_0), \quad\forall \tau\in\mathbb{R}.$$

Using the convergence of $u_\lambda$ to $v$ we get

\begin{multline}
v(x+\tau x_0)=\lim_{\lambda\rightarrow 0^+}\frac{u(x_0+\lambda (x+\tau x_0))-u(x_0)}{\lambda}=\\
=\lim_{\lambda\rightarrow 0^+}\frac{u((1+\lambda\tau)x_0+\lambda x)-u(x_0)}{\lambda}=\\
=\lim_{\lambda\rightarrow 0^+}\frac{(1+\lambda\tau)u(x_0+\frac{\lambda}{1+\lambda\tau}x)-(1+\lambda\tau)u(x_0)
+\lambda\tau u(x_0)}{\lambda}\\
=\lim_{\lambda\rightarrow 0^+}\frac{(u(x_0+\frac{\lambda}{1+\lambda\tau}x)-u(x_0))}
{\frac{\lambda}{(1+\lambda\tau)}}+\tau u(x_0)=v(x)+\tau u(x_0).
\end{multline}

Observe that we used that for $|\lambda\tau|<1$ we have $1+\lambda\tau>0$ and, thus, 
$u((1+\lambda\tau)x_0+\lambda x)=(1+\lambda\tau)u(x_0+\frac{\lambda}{1+\lambda\tau}x)$ because $\M{G}_u$ is a cone.
Choosing a basis of $\R{n}$ of the form
$\left\{ \frac{x_0}{|x_0|}, v_2,\ldots,v_n \right\},$
where $v_2,\ldots,v_n$ is a completion of $\frac{x}{|x_0|}$ to an orthonormal basis, we have that $v$ satisfies
\eqref{conooriz}.

To see that the graph of $\widetilde{v}:=v\big{|}_{\{0\}\times\R{n-1}}$ is minimal, consider a vector field
$\widetilde{X}(\widetilde x,y_1,\ldots,y_m)$ in $C^1_0(\R{n-1+m};\R{n-1+m})$ and a function $\rho:\mathbb{R}\rightarrow\mathbb{R}_{\geq 0}$ compactly supported and non identically zero. Let
$$X(x_1,\ldots,x_n)=\rho(x_1)\widetilde{X}(\widetilde x).$$ Since $\M{G}_v$ is minimal and thanks to Fubini-Tonelli's theorem we get

\begin{multline}
0=\int_{\M{G}_u}\diver X=\int_{\R{n}}\diver X(x)\sqrt{g_u(x)} dx=\\
=\int_{\mathbb{R}}\frac{\partial \rho}{\partial x_1}\rho (x_1)dx_1\int_{\R{n-1}}\widetilde{X}(\widetilde{x})
\sqrt{g_u(\widetilde{x})}d\widetilde{x}+\\
+\int_{\mathbb{R}}\rho(x_1)dx_1\int_{\R{n-1}}\widetilde{\diver}\widetilde{X}(\widetilde{x})d(\widetilde{x})=
c\int_{\M{G}_{\widetilde{v}}} \widetilde{\diver}{\widetilde{X}}
\end{multline}
In the last inequality we used $\int_{\mathbb{R}}\frac{\partial\rho}{\partial x_1}(x_1)dx_1=0$, being the support of $\rho$ compact;
we also used that $\sqrt{g_u}$ depends only on $\widetilde{x}$ thanks to \eqref{conooriz}. We denoted by
$\widetilde{\diver}$ the divergence operator on $\{0\}\times\R{n-1}\times\R{m}$ and, finally,
$c=\int_{\mathbb{R}}\rho(x_1)dx_1.$
\end{proof}

\begin{prop}[Blow-downs]\label{blowdown}
Let $u:\R{n}\rightarrow\R{m}$ be a Lipschitz map, $|Du|\leq K$, with $\M{G}_u$ minimal as varifold.
Let $u_\lambda$ be defined by
$$u_\lambda (x)= \frac{1}{\lambda}u(\lambda x).$$
Then there exists a sequence $\lambda(i)\rightarrow \infty$ such that $u_{\lambda(i)}\rightarrow v$ uniformly on compact sets and in the sense of varifolds, where the graph of $v$ is a cone minimal in the sense of varifolds.
\end{prop}

\begin{proof}
As for proposition \ref{blowup}, with $\lambda\rightarrow\infty$ instead of $\lambda\rightarrow 0.$
\end{proof}

\section{A Bernstein-type theorem}

A Bernstein-type theorem is a rigidity theorem which, under suitable hypothesis, implies that an entire minimal graph is an affine subspace. The original statement is:

\begin{trm} Let $u:\R{2}\rightarrow\R{}$ be a $C^2$ function satisfying the minimal surface equation. Then $u$ is affine, i.e. $u(x,y)=y_0+\sigma_1 x+\sigma_2 y$, with  $\sigma_1,\sigma_2\in\R{}$.
\end{trm}

It comes from a memoir of Bernstein published in 1927, but several alternative proofs and several generalizations are now available. De Giorgi \cite{degiorgi2} proves that Bernstein's theorem holds true for 3 dimensional graphs in $\R{4}$, while Simons in \cite{simons} generalizes Bernstein's theorem in $\R{n+1}$ for $n\leq 7$. This result is sharp for what concerns the dimensions because in \cite{bdg} Bombieri, De Giorgi and Giusti show that there exists a non-affine function $u:\R{8}\rightarrow \R{}$ whose graph is minimal.

Some years before Moser had proved in \cite{moser} that the minimal graph of a scalar function whose gradient is bounded is an affine subspace. 

\medskip

In higher codimension, Lawson and Osserman \cite{lawson} have shown that the cone over Hopf's map \eqref{conohopf} is minimal, theorem \ref{controregol}. Moreover such a cone is the graph of a function with bounded gradient; this is in constrast with Moser's result for codimension 1.

The first Bernstein-type theorems in arbitrary codimension were proved in \cite{hildebrandt} by Hildebrandt, Jost and Widman who studied the Gauss map of a minimal graph. With a similar approach, Jost and Y. L. Xin in \cite{jost} improve the result in \cite{hildebrandt}, obtaining the following theorem.

\begin{trm}\label{bernjost} Let $u:\R{n}\rightarrow \R{m}$ be a smooth function satisfying the minimal surface system \eqref{mssnonpar}. Let $*\omega=(\det (I+Du^*Du))^{-\frac{1}{2}}$ and set $\beta_0>0$ such that
\begin{equation}\label{betazero}
\beta_0<\left\{
\begin{array}{ll}
2 & \mathrm{if}\;m\geq2\\
\infty &\mathrm{if}\;m=1.
\end{array}
\right.
\end{equation}
Then, if $*\omega\geq \frac{1}{\beta_0}$, $u$ is affine.
\end{trm}

Comparing this theorem with Moser's result, we remark that the hypothesis $*\omega\geq
\frac{1}{\beta_0}$ implies $Du\leq \beta_0^2-1$, while in codimension, though we require the gradient of $u$ to be bounded, we do not impose a specific constant to bound it. The theorem we will prove below, due to Mu-Tao Wang \cite{wangbern}, implies the result of Moser for codimension 1 and, as we will show, the result of Jost and Xin in arbitrary codimension. It is a natural extension of Moser's theorem because it only requires $|Du|$ to be bounded and area-decreasing. This last assumption is always met in codimension 1, thus the hypothesis of Wang's theorem, in codimension 1, reduces to the hypothesis of Moser's theorem.

\begin{trm}\label{teoremabern}
Let $u:\R{n}\rightarrow\R{m}$ be a $C^2$ \emph{area-decreasing} map satisfying the minimal surface system \eqref{mss}. Assume that $\abs{Du}\leq K$ for some $K >0$. Then $u$ is linear.
\end{trm}

\begin{proof} Let $\delta>0$ be such that $\lambda_i\lambda_j\leq 1-\delta$ for $i\neq j$, where the $\lambda_i$s are the singulare values of $Du$.

\textbf{1.} Denote by $\Delta_\Sigma$ the Laplacean on $\Sigma=\M{G}_u$, with respect to the parametrization $x\rightarrow (x,u(x))$.
\begin{equation}\label{bern}
\Delta_\Sigma(\ln *\omega)=\frac{*\omega\Delta_\Sigma*\omega-\abs{\nabla^\Sigma*\omega}^2}{\abs{*\omega}^2}.
\end{equation}
The covariant derivative of $*\omega$ may be computed using the singular value decomposition of $Du$ and equations
\eqref{pi} and \eqref{gradomega}:
\begin{equation}\label{gradomega2}
(*\omega)_k=-*\omega(\sum_{i,\alpha} \lambda_{\alpha i} h^\alpha_{ik})=-*\omega(\sum_i\lambda_{i} h^i_{ik}).
\end{equation}
Equation \eqref{evolomega3} may be rewritten as
\begin{equation}\label{bern2}
\Delta_\Sigma *\omega=- *\omega\Big(\sum_{\alpha,i,k} (h^\alpha_{ik})^2
+2\sum_{i,j, k\atop i< j}(-\lambda_j\lambda_i h^i_{ik}h^j_{jk}+\lambda_j\lambda_i h^i_{jk} h^j_{ik})\Big).
\end{equation}
This is easily seen by swapping $i$ with $j$ and $\alpha$ with $\beta$ when summing over these indices.
Inserting \eqref{bern2} and \eqref{gradomega2} into \eqref{bern} yields
\begin{multline}\label{lnsubarm}
\Delta_\Sigma(-\ln *\omega)=\\=*\omega\Big(\sum_{\alpha,i,k}(h^\alpha_{ik})^2
+2\sum_{i,j, k\atop i< j}\lambda_j\lambda_i h^i_{jk} h^j_{ik}-2\sum_{i,j, k\atop i< j}\lambda_j\lambda_i h^i_{ik}h^j_{jk}
+\sum_{k}\Big(\sum_{i}\lambda_i h^i_{ik} \Big)^2\Big)=\\
=*\omega\Big(\sum_{\alpha,i,k}(h^\alpha_{ik})^2+\sum_{i,k}\lambda_i^2(h^i_{ik})^2+2\sum_{k,i,j\atop i\neq j}
\lambda_i\lambda_jh^i_{jk}h^j_{ik}\Big)\leq\\
\geq *\omega\Big(\sum_{\alpha,i,k}(h^\alpha_{ik})^2+2\sum_{k,i,j\atop i\neq j}\lambda_i\lambda_jh^i_{jk}h^j_{ik}\Big)
\geq *\omega\Big(|A|^2-(1-\delta)|A|^2 \Big)\geq\delta|A|^2*\omega.
\end{multline}
We observe that, in order to prove that $-\ln*\omega$ is a subharmonic function, we only used the area-decreasing condition
$|\lambda_i\lambda_j|\leq 1$: thus we have shown that $-\ln*\omega$ is subharmonic on any area-decreasing graph.

The boundedness of $\abs{Du}$ implies that $*\omega\geq K_1>0$ for some $K_1$, whence
$\ln(*\omega)$ is bounded from below by $\ln K_1$.

\textbf{2.} We perform a blow-down of the graph of $u$, and by proposition \ref{blowdown} we get an equiLipschitz sequence
$$u_{\lambda(j)}(x)=\frac{1}{\lambda(j)}u(\lambda(j)x)$$
uniformly converging to a Lipschitz function $\Cl{u}$ for which the same condition on the singular values of the differential holds. In particular $\Cl u$
is area-decreasing. Moreover the convergence is also in the sense of varifolds and $\M{G}_{\Cl{u}}$ is a minimal cone with vertex in the origin. The differential of $\Cl{u}$ is positively homogeneous, that is
$$D\Cl{u}(tx)=D\Cl{u}(x),\quad t>0,\, x\in\R{n}\backslash \{0\}.$$
Observe that the cone in not necessarily regular in the origin, but we will assume that it is so in every other point. The general case is studied in step 3. The homogeneity of $D\Cl{u}$ implies that on every annulus with center in the origin $*\omega$ attains an interior minimum; on the other hand, the maximum principle applies to \eqref{bern2} and thus $\abs{A}=0$ in every annulus and so in all of
$\R{n}\backslash\{0\}.$ The vanishing of the second fundamental form implies that the cone is a linear subspace, so that $\Cl{u}$ is linear.
We now prove that $Du(x)=D\Cl{u}(0)$ for every $x\in \R{n}$ and, thus, $u$ is linear. Let $\delta$ and $\gamma$ be as in Allard's theorem \ref{regolallard}, $j_0$ and $\rho$ such that for every $j\geq j_0$ 
$$\frac{ \M{H}^n\res\M{G}_{u_j}(B(x_0,\rho))}{\rho^n \omega_n}\leq 1+\delta,$$
where this is possible because from the varifold convergence we get
$$\frac{ \M{H}^n\res\M{G}_{u_j}(B(x_0,\rho))}{\rho^n \omega_n}\rightarrow
\frac{ \M{H}^n\res\M{G}_{\Cl u}(B(x_0,\rho))}{\rho^n \omega_n}=1.$$
Then $u_j\in C^{1,\alpha}(\Cl{B}^n_{\gamma\rho}(0))$ and they are equibounded in $C^{1,\alpha}(\Cl{B}^n_{\gamma\rho}(0))$, thanks to \eqref{stimaallard2}, where $B^n_r(x)$ is the ball in $\R{n}$ centered in $x$ with radius $r$. By Ascoli-Arzel\`a's theorem, a subsequence, still denoted by $u_j$, $C^1$-converges in $B(0,\gamma\rho)$ to the linear map $\Cl{u}$.
For every $x\in \R{n}$ we have, for $j$ large enough, $\frac{x}{\lambda(j)}\in B(0,\gamma\rho)$ and
$$\abs{Du_j\bra{\frac{x}{\lambda(j)}}-D\Cl{u}\bra{\frac{x}{\lambda(j)}}}<\varepsilon.$$
As $\varepsilon$ goes to $0$, observing that $Du_j\bra{\frac{x}{\lambda(j)}}=Du(x)$, we get that $Du$ is constant and, by Lagrange's mean value theorem, we deduce that $u$ is linear.

\textbf{3.} If the blow-uo generates a cone which has at least a singularity in $x_0\neq 0$, we may perform a blow-up in 
$x_0$ and, by proposition \ref{blowup2} we obtain a minimal cone of dimension $n-1$ in $\R{n+m-1}$. If such a cone is $C^2$ but in the origin, we apply step 2 to prove that the singularity originating the cone couldn't exist, absurd. 
Otherwise we go on performing blow-ups, until we get a cone with at most a singularity in the origin. Such a cone must exist because in dimension one a cone is just a union of two half-lines. Again thanks to step 2 we obtain an absurd.
\end{proof}

\begin{oss} This theorem implies theorem \ref{bernjost} because the hypothesis \eqref{betazero} may be written as
$$1+\sum_{i=1}^n\lambda_i^4+2\sum_{1\leq i<j\leq n}\lambda_i^2\lambda_j^2
\leq\prod_{i=1}^n(1+\lambda_i^2)=\frac{1}{*\omega^2}\leq 4-\delta, \quad \delta>0.$$
If for every $1\leq i\leq n$ we have $\lambda_i^4<1-\frac{\delta}{2}$, then $\lambda_i\lambda_j\leq 1-\varepsilon$ for some $\varepsilon$ and thus $u$ is area-decreasing. If for some $i$ we have $\lambda_i^4\geq \frac{\delta}{2}$ we get
$$\sum_{1\leq i<j\leq n}\lambda_i^2\lambda_j^2\leq 1-\varepsilon$$
for some $\varepsilon,$ thus again the area-decreasing condition.
\end{oss}

\subsection{Remarks to Bernstein's theorem: the Gauss map}
Mu-Tao Wang's proof of theorem \ref{teoremabern} is based on inequality \eqref{lnsubarm} which says that $-\ln*\omega$ is a function \emph{subharmonic} on $\Sigma$ (with respect to the Riemannian metric of $\Sigma$). We explain this further.

\begin{defin}[The Gauss map]
Given a submanifold of dimension $n$ $\Sigma\subset\R{n+m}$, its Gauss map
$$\gamma:\Sigma\rightarrow G(n,m)$$
is the map associating to each $x\in\Sigma$ the tangent space $T_x\Sigma$, seen as an element of the Grassmannian of $n$-planes in $\R{n+m}$.
\end{defin}

The differentiable and Riemannian structure of $G(n,m)$ have been studied by Yung-Chow Wong in \cite{wong} and Jost and Xin in \cite{jost}. The fundamental theorem concerning the Gauss map of a minimal surface is due to Ruh and Vilms \cite{ruh}:

\begin{trm} The Gauss map $\gamma$ of a submanifold $\Sigma\subset\R{n+m}$ is harmonic if and only if the mean curvature $H$ of $\Sigma$ is parallel, i.e. $$\nabla^\Sigma H=0.$$
\end{trm}

In particular, if $\Sigma$ is minimal, i.e. $H=0$, it's Gauss map is harmonic. Jost and Xin observe that the condition $*\Omega\geq\frac{1}{\beta_0}$, determines a region in the Grassmannian in which
\begin{equation}\label{fconv}
f(L):=\ln \sqrt{\det(I+L^*L)}
\end{equation}
is convex\footnote{convex here means that, given a geodesic $\gamma\rightarrow\Xi$, where $\Xi\subset G(n,m)$ is the subset of the Grassmannian containing the graphs of area-decreasing linear maps, we have that
$$\frac{d^2}{dt^2}\Big(-\ln\sqrt{\det (I+L^*L)}\circ\gamma\Big)\geq 0.$$ This notion of convexity is different from the one used in codimension 1 when we say that $\sqrt{1+|Du|^2}$ is a convex function: in the latter case, indeed, the $1\times n$-matrix space where $Du$ lives is given the flat metric, which is different from the Riemannian metric on the Grassmannian.} (in \eqref{fconv} we identify a plane with the linear map
$L:\R{n}\rightarrow\R{m}$ of which it's the graph; we will only consider the region of the Grassmannian given by such planes). In \cite{wangmean}, Mu-Tao Wang proves that $f$ is convex on a larger region of the Grassmannian: the graphs of area-decreasing linear maps. Thus $*\omega=f\circ\gamma$ is subharmonic because it's the composition of a harmonic map and a convex function.

It is easily seen that a further extension of such a region of the Grassmannian on which $f$ is convex would yield an extension of theorem \ref{teoremabern} and of the consequent regularity theorem \ref{wangregol}.

\section{Regularity of area-decreasing minimal graphs}

The following regularity theorem, due to Mu-Tao Wang \cite{wangdir}, is consequence of Allard's regularity theorem and of the theorem of Bernstein \ref{teoremabern}. The hypothesis are similar to the one in codimension 1, since in the latter case every function is area-decreasing. We also remark that, due to the counterexample of Lawson and Osserman, theorem \ref{controregol}, an hypothesis on $Du$ is natural. 

\begin{trm}\label{wangregol}
Let be given a Lipschitz map $u:\Cl\Omega\rightarrow\R{m}$ satisfying the minimal surface system \eqref{mssnonpar} and assume that there exists $\varepsilon>0$ such that
$$\lambda_i \lambda_j\leq 1-\varepsilon, \quad i\neq j,$$
where the $\lambda_i$s are the singular values of $Du$. Then $u\in   C^\infty(\Omega;\R{m}).$
\end{trm}

\begin{proof}

\textbf{1.}
Let $x_0\in\Omega$. Up to translation, we may and do assume $x_0=0$ and $u(0)=0$. Performing a blow-up in $0$, see proposition \ref{blowup}, we get $u_i:=u_{\lambda(i)}\rightarrow v$ uniformly and in the sense of varifolds, where $\M{G}_v$ is a cone minimal as varifold. Moreover the convergence is uniform and preserves the area-decreasing condition as well as $\sup|Du|$.

If $v$ is $C^2$ but in $0$, then $v$ is affine by theorem \ref{teoremabern}. In particular $\M{H}^n\res\M{G}_v(B_1(0))=\omega_n$, where
$\omega_n=\M{L}^n(B^n_1(0))$.
From the varifold convergence (uniform convergence wouldn't be enough) we get
$$\frac{\M{H}^n\res\M{G}_u(B_\lambda(0))}{\omega_n \lambda^n}
=\frac{\M{H}^n\res\M{G}_{u_\lambda}(B_1(0))}{\omega_n}\rightarrow \frac{\M{H}^n\res\M{G}_v(B_1(0))}{\omega_n}=1$$
Set $\delta$ and $\gamma$ as in Allard's theorem \ref{regolallard}; let
$V=\V(\M{G}_u,1)$ and
$\rho>0$ such that $B_\rho(0)\subset\Omega\times\R{m}$ and
$$\frac{\mu_V(B_\rho(0))}{\omega_n\rho^n}\leq 1-\delta.$$
We apply Allard's regularity theorem, whence
$u\in C^{1,\alpha}(\Cl{B}_{\gamma\rho}(0))$.

\medskip

\textbf{2.} Now assume that the minimal cone generated by the blow-up in the above step 1 is not $C^2$ in all of $\R{n}\backslash\{0\}.$
As in the proof of theorem \ref{teoremabern}, assume that there exists a singularity in $x_0\neq 0$. We may generate another cone in $(x_0,v(x_0))$ with another blow-up. Thanks to proposition \ref{blowup2}, such a cone factorizes and we obtain an $(n-1)$-dimensional cone. If it's smooth but at most in the origin, applying step 1 we obtain that $v$ is smooth is $x_0$, absurd.

Then, by induction, we perform blow-ups and find cones with singularities until we find a cone of dimension 1, union of two straight lines. And this may not be singular, but in the origin.

\textbf{3.} The existence and smoothness of the higher order derivatives is consequence of Schauder estimates and is contained in the theorem which follows (which also says that $C^1$ solutions are $C^{1,\alpha}$, then smooth).
\end{proof}

\begin{oss} Thanks to the second Allard's theorem \ref{trmallard2}, the solutions $u$ of the Dirichlet problem are smooth up to the boundary if $\Omega$ is strictly convex.
\end{oss}

\begin{trm}[Morrey \cite{morrey2}]\label{morrey}
An application $u\in C^1(\Omega;\R{m})$ which is a solution to the minimal surface system in divergence form \eqref{mssnonpar} is analytic.
\end{trm}

\begin{cenni} For the analyticity see \cite{morrey}. We will only prove smoothness.

\textbf{1.} By the difference quotient method, a $W^{1,2}(\Omega;\R{m})$-weak solution to the minimal surface system \eqref{mssnonpar} is $W^{2,2}$.

\textbf{2.} Deriving the minimal surface system we obtain that the first derivatives $D_s u$ satisfy
\begin{equation}\label{eqderiv}
\frac{\partial}{\partial x^i}\Big( A^{ij}_{\alpha\beta}(Du(x))\frac{\partial}{\partial x^j}(D_s u^\beta(x))\Big)=0
\end{equation}
in the weak sense, where
$$A^{ij}_{\alpha\beta}(p):=\frac{\partial^2 F}{\partial p^\alpha_i\partial p^\beta_j}(p),$$ 
$$F(p):=\sqrt{\det\big(I+p^*p\big)},\forall p\in M^{m\times n}.$$
We see $p$ as an $m\times n$-matrix. $F$ is the area integrand \eqref{eqareagraf} and it's stricly polyconvex\footnote{a function $F$ defined on a space of $M^{m\times n}$-matrices is said polyconvex if there exists a convex function $g:\R{d}\rightarrow\R{}$ such that
$$F(p)=g(\Min(p)), \quad \forall p\in M^{m\times n},$$
where $\Min:M^{m\times n}\rightarrow \R{d}$ is the map associating to each matrix the set of its minors. Thanks to the  Cauchy-Binet's formula \cite{federer}, \cite{giaquinta},
$$F(p)=\sqrt{\sum_k M_k^2},$$
being $M_k^2$ the sum of the squares of the minors $k\times k$. The strict polyconvexity ($g$ is stricly convex) implies the Legendre-Hadamard condition.}; its derivatives satisfy the Legendre-Hadamard condition: for each $p\in M^{m\times n}$ there exists $\lambda>0$ such that
\begin{equation}\label{eqlh}
A^{ij}_{\alpha\beta}(p)\xi_i\xi_j\eta^\alpha\eta^\beta\geq\lambda|\xi|^2|\eta|^2.
\end{equation}

\textbf{3.} Since $u\in C^1(\Omega;\R{m})$, we have that $A^{ij}_{\alpha\beta}(Du)$ is continuous. From now on let $v:=D_s u$, $x_0\in\Omega$, $0<\rho<R<\dist(x_0,\partial\Omega)$ and let $\Cl v$ be the weak solution of
\begin{equation}\label{eqlin}
\left\{
\begin{array}{ll}
\displaystyle \frac{\partial}{\partial x^i}\Big( A^{ij}_{\alpha\beta}(Du(x_0))\frac{\partial}{\partial x^j}(\Cl v(x))\Big)=0 
&\mathrm{in}\; B_R(x_0);\\
\Cl v=v & \mathrm{su} \;\partial B_R(x_0).\rule{0cm}{.6cm}
\end{array}
\right.
\end{equation}
Such a solution exists because, being $A^{ij}_{\alpha\beta}(Du(x_0))$ constant, system \eqref{eqlin} is linear. The classical energy inequalities for $\Cl v$ are
\begin{equation}\label{eqclassica1}
\int_{B_\rho}|D\Cl v|^2\leq c\Big(\frac{\rho}{R}\Big)^n\int_{B_R}|D\Cl v|^2
\end{equation}
and, if we define the mean of a function $f_{x_0,\rho}:=\frac{1}{|B_\rho|}\int_{B_\rho(x_0)} f,$
\begin{equation}\label{eqclassica2}
\int_{B_\rho}|D\Cl v -(D\Cl v)_{x_0,\rho}|^2\leq c\Big(\frac{\rho}{R}\Big)^{n+2}\int_{B_R}|D\Cl v-(D\Cl v)_{x_0,\rho}|^2
\end{equation}

Thus $v= \Cl v+(v-\Cl v)$ satisfies
\begin{multline}\label{eqclassica3}
\int_{B_\rho}|Dv|^2\leq\\
\leq \int\limits_{B_\rho}|D\Cl v|^2+\int\limits_{B_\rho}|D\Cl v-Dv|^2\leq c\Big(\frac{\rho}{R}\Big)^n\int\limits_{B_R}|D\Cl v|^2+ \int\limits_{B_\rho}|D\Cl v-Dv|^2\leq \\ \leq c\Big(\frac{\rho}{R}\Big)^n\int_{B_R}|Dv|^2+c_1\int_{B_R}|D\Cl v-Dv|^2;
\end{multline}
\begin{multline}\label{eqclassica4}
\int_{B_\rho}|Dv-(Dv)_{x_0,\rho}|^2\leq \\ \leq\int_{B_\rho}|D\Cl v-(D\Cl v)_{x_0,\rho}|^2+\int_{B_\rho}|D(v-\Cl v)-(D\Cl v-v)_{x_0,\rho}|^2\leq\\
\leq c\Big(\frac{\rho}{R}\Big)^{n+2}\int_{B_R}|D\Cl v-(D\Cl v)_{x_0,\rho}|^2+\int_{B_\rho}|D\Cl v-Dv|^2\leq \\
\leq c\Big(\frac{\rho}{R}\Big)^{n+2}\int_{B_R}|Dv-(Dv)_{x_0,\rho}|^2+c_1\int_{B_\rho}|D\Cl v-Dv|^2.
\end{multline}

Putting together \eqref{eqderiv} and \eqref{eqlin}, omitting the indeces $i,j,\alpha,\beta$ and writing $A(Du(x))=A(x)$ we obtain
$$D \big(A(x_0)D (v-\Cl v)\big)=D\big([A(x_0)-A(x)]Dv \big)$$
Since $\Cl v-v\in W^{1,2}_0(B_R)$ we may take this as test function in the above equation; using ellipticity \eqref{eqlh} and integrating by parts we get

\begin{multline}
\int_{B_R}|D(v-\Cl v)|^2\leq\int_{B_R}A(x_0)D(v-\Cl v)D(v-\Cl v)\leq\\ \leq
\int_{B_R}[A(x_0)-A(x)]Dv(x)D(v(x)-\Cl v(x))dx.
\end{multline}
Applying $ab\leq \varepsilon a^2+\frac{b^2}{\varepsilon}$ with $a=|A(x)-A(x_0)||Dw|$ and $b=|D(v-\Cl v)|$ we obtain
\begin{equation}\label{eqdw}
\int_{B_R}|D(v-\Cl v)|^2\leq c\omega(R)^2\int_{B_R}|Dv|^2,
\end{equation}
being $c$ an absolute constant and $\omega(R):=\sup_{B_R}|A(x)-A(x_0)|.$

Estimate \eqref{eqdw} plugged into \eqref{eqclassica3} yields
$$\int_{B_\rho}|Dv|^2\leq c\Big[\Big(\frac{\rho}{R}\Big)^n+\omega(R)^2\Big]\int_{B_R}|Dv|^2;$$
choosing $R$ in such a way that $\omega(R)<\delta$ for some $\delta>0$ and applying an algebraic lemma
\footnote{\textbf{Lemma} Let a positive non decreasing function $\Phi$ and positive constants $A,B,\alpha,\beta$ be given with $\alpha>\beta$ an let $R_0$ be such that
$$\Phi(\rho)\leq A\Big[\Big(\frac{\rho}{R} \Big)^\alpha+\delta \Big]\Phi(R)+BR^\beta,\quad 0<\rho<R\leq R_0;$$
then there is a constant $c=c(\alpha,\beta,A,B,\delta)$ such that
$$\Phi(\rho)\leq c(\alpha,\beta,A,B,\delta)\Big[\Big(\frac{\rho}{R} \Big)^{\beta}\Phi(R)+B\rho^\beta\Big],\quad 0<\rho<R\leq R_0.$$}
we obtain
\begin{equation}\label{nepsilon}
\int_{B_\rho}|Dv|^2\leq c\rho^{n-\varepsilon},
\end{equation}
being $c$ a constant depending on the oscillation $\omega(R)$.
Estimate \eqref{nepsilon}, as $x_0$ ranges in an open set, implies that $Dv\in L^{2,n-\varepsilon}_{\loc}(\Omega)$.\footnote{$L^{p,\lambda}$ and $\mathcal{L}^{p,\lambda}$ are the spaces of Morrey and Campanato respectively.
$$L^{p,\lambda}(\Omega):=\Big\{u\in L^p(\Omega):\sup_{x_0\in\Omega\atop 0<\rho<\diam\Omega} \frac{1}{\rho^\lambda}\int\limits_{\Omega\cap B_\rho(x_0)}|u|^pdx<+\infty \Big\}$$
$$\mathcal{L}^{p,\lambda}(\Omega):=\Big\{u\in L^p(\Omega):\sup_{x_0\in\Omega\atop 0<\rho<\diam\Omega} \frac{1}{\rho^\lambda}\int\limits_{\Omega\cap B_\rho(x_0)}|u-u_{x_0,\rho}|^pdx<+\infty \Big\} $$
We use the following theorem of Campanato \cite{campanato}:
\begin{trm} Let $n<\lambda\leq n+p$; then
$$\mathcal{L}^{p,\lambda}(\Omega)\cong C^{0,\sigma}(\Omega),\quad \sigma=\frac{\lambda-n}{p},$$
while, for $\lambda> n+p$, the space $\mathcal{L}^{p,\lambda}$ contains only the constant functions.
\end{trm}}
Using Poincar\'e's inequality
$$\int_{B_\rho}|v-v_{x_0,\rho}|^2\leq c\rho^2\int_{B_\rho}|Dv|^2\leq c\rho^{n+2-\varepsilon},$$
which, given in an open set, implies $v\in \mathcal{L}^{2,n+2-\varepsilon}_{\loc}(\Omega)\cong C^{0,\sigma}(\Omega)$ for $\sigma=\frac{2-\varepsilon}{2},$ that is $u\in C^{1,\sigma}(\Omega)$.

\textbf{4.} Thanks to the preceding step $u\in C^{1,\sigma}(\Omega;\R{m})$ and consequently $A^{ij}_{\alpha\beta}\in C^{0,\sigma}(\Omega).$ We may, thus, estimate the modulus of continuity $A(x)$ having $\omega(R)\leq R^\alpha$.
Plugging \eqref{eqdw} into \eqref{eqclassica4} and applying \eqref{nepsilon} we get
\begin{multline}\label{eqdw2}
\int\limits_{B_\rho}|Dv-(Dv)_{x_0,\rho}|^2\leq c\Big(\frac{\rho}{R}\Big)^{n+2}\int\limits_{B_R}|Dv-(Dv)_{x_0,\rho}|^2+c\omega(R)^2\int\limits_{B_R}|Dv|^2\leq\\
\leq c\Big(\frac{\rho}{R}\Big)^{n+2}\int_{B_R}|Dv-(Dv)_{x_0,\rho}|^2+R^{n+2\sigma-\varepsilon}.
\end{multline}
Applying the algebraic lemma with $\beta=n+2\sigma-\varepsilon$ we conclude that
$$\int_{B_\rho}|Dv-(Dv)_{x_0,\rho}|^2\leq c\rho^{n+2\sigma-\varepsilon},$$
thus $Dv\in \mathcal L^{2,2+2\sigma-\varepsilon}_{\loc}(\Omega)\cong C^{0,\sigma-\frac{\varepsilon}{2}}(\Omega);$
now we know that $Dv$ is locally bounded and we may set $\varepsilon=0$ in \eqref{eqdw2}, concluding that $Dv\in C^{0,\sigma}(\Omega),$ i.e. $u\in C^{2,\sigma}(\Omega)$. 

\textbf{5.} Being $u\in C^{2,\sigma}(\Omega),$ the minimal surface system my be written in the non-variational form
$$\sum_{i,j=1}^n g^{ij}(Du)D_{ij}u^\alpha=0,\quad \alpha=1,\ldots,m.$$
Since $g^{ij}\in C^{1,\sigma}(\Omega),$ we may derive the system and obtain, again using the difference quotient method,
\begin{equation}\label{eqnonvar}
\sum_{i,j=1}^ng^{ij}(Du)D_{ij}(D_su)=-\sum_{i,j=1}^n D_sg^{ij}(Du)D_{ij}u=:h(Du).
\end{equation}
From the classical Schauder estimates we know that, if $a^{ij},f\in C^{k,\sigma}(\Omega)$, then the solution of
$$\sum_{i,j=1}^n a^{ij}(x)D_{ij}u(x)=f(x)$$
is in $C^{k+2,\sigma}(\Omega)$.

We prove inductively that the solutions $u$ to the minimal surface system are $C^{k,\sigma}$ for every $k\in\mathbb{N}$. For $k=1,2$ it's the result of steps 4 and 5. Assume inductively that $u\in C^{k,\alpha}(\Omega)$. Then $g^{ij}(Du(x)),h(Du(x))\in C^{k-1,\alpha}$ and thanks to the Schauder estimates applied to \eqref{eqnonvar}, $Du\in C^{k,\alpha}(\Omega)$. The induction is proved and, thus, $u\in C^{\infty}(\Omega;\R{m}).$
\end{cenni}

\appendix
\chapter{Geometry of Varifolds}\label{capitolovar}

\section{Rectifiable subsets of $\R{n+m}$}
Most of the definitions and propositions of chapter \ref{capitolosottov} may be applied to particular subsets of $\R{n+m}$ not necessarily having a $C^1$-submanifold structure. What we are seeking for is a class of measurable subsets of $\R{n+m}$ large enough to contain the graphs of Lipschitz functions, though containing only objects on which we may develop the standard notions of differential calculus.

The necessity to consider objects more general than smooth submanifolds may be appreciated in theorem \ref{convergenza}. Its proof uses the compactness theorem of Ascoli and Arzel\`a, giving the uniform convergence of equicontinuous and equibounded functions. On the other hand, the uniform limit of $C^1$ functions with equibounded gradients is not necessarily a $C^1$ function, but it is definitely a Lipschitz function.

This considerations suggest the interest of the notion of \emph{$n$-rectifiable} set:

\begin{defin}
A Borel subset $M\subset\R{n+m}$ is said to be countably $n$-rectifiable if
\begin{equation}\label{unionec1}
M\subset N_0\cup\Big(\bigcup_{j=1}^{\infty} N_j\Big),
\end{equation}
where $\M{H}^n(N_0)=0$ and, for $j\geq 1$, $N_j$ is a $C^1$ submanifold of $\R{n+m}$ of dimension $n$.
\end{defin}

The connection between rectifiable sets and Lipschitz functions is essentially a consequence of the theorems of Rademacher and Whitney; for their proofs see \cite{federer}, \cite{giaquinta} and \cite{simon}.

\begin{trm}[Rademacher]
Every Lipschitz function $f:\R{n}\rightarrow \R{}$ is $\M{L}^n$-almost everywhere differentiable, where $\M{L}^n$ is the Lebesgue measure in $\R{n}$. In particular its gradient is a.e. well defined
$$\nabla f:=\Big(\frac{\partial f}{\partial x^1},\cdots,\frac{\partial f}{\partial x^n}\Big)$$
and $\M{L}^n$-a.e. we have
$$\lim_{x\rightarrow x_0}\frac{f(x)-f(x_0)-\nabla f\cdot(x-x_0)}{|x-x_0|}=0. $$
\end{trm}

\begin{oss} $\nabla f$ is the a.e. limit of measurable functions (the difference quotients) and is thus measurable. Moreover, if $f$ is Lipschitz with Lipschitz constant $K$, it's clear that $|\nabla f|\leq K$, so that
$\nabla f\in L^{\infty}(\R{n};\R{n}).$
\end{oss}

The following theorem will be referred to as theorem of Whitney because it is an almost immediate consequence of a celebrated theorem of Whitney.

\begin{trm}[Whitney]
Let $f:\R{n}\rightarrow\R{}$ be a Lipschitz function. Then for every $\varepsilon>0$ there exists a function
$h:\R{n}\rightarrow \R{}$ of class $C^1$ such that
$$\M{L}^n\big(\{x\in\R{n}:f(x)\neq h(x)\}\cup\{x\in\R{n}:\nabla f(x)\neq \nabla h(x)\}\big)<\varepsilon.$$
\end{trm}
Thanks to Rademacher's theorem the right term in the union is well defined up to $\M{L}^n$-null sets.

\begin{prop}[Characterization of rectifiable sets]\label{carattrett}\ \\ A subset $M\subset\R{n+m}$ is countably $n$-rectifiable if and only if there exists a sequence of Lipschitz maps $F_j:\R{n}\rightarrow\R{n+m}$ and a set $M_0$ with $\M{H}^n(M_0)=0$ such that
\begin{equation}\label{unionelip}
M=M_0\cup\Big(\bigcup_{j=1}^n F_j(A_j)\Big),
\end{equation}
where $A_j\subset\R{n}$ is measurable for every $j$.
\end{prop}

\begin{proof} ($\Rightarrow$) Every $C^1$-submanifold $N_j$ in $\R{n+m}$ is locally the image of $C^1$-maps which we denote by
$h_{ij}:B^n\subset\R{n}\rightarrow\R{n+m}$. Therefore 
\begin{equation}\label{unione3}
N_j\subset E_j\cup\Big( \bigcup_{i=1}^{\infty}h_{ij}(B^n)\Big),\quad\M{H}^n(E_j)=0.
\end{equation}
If \eqref{unionec1} holds, choose $h_{ij}$ as said in such a way that \eqref{unione3} holds true. Let $A_{ij}:=g_{ij}^{-1}(M)$ and $N_0:=\bigcup_{j=1}^{n}E_j$. Then
$$M=N_0\cup\Big(\bigcup_{i,j=1}^{n}h_{ij}(A_{ij})\Big).$$
Since $A_{ij}$ is Borel (because inverse image of a Borel set) and since we may assume $g_{ij}$ to be Lipschitz, we get \eqref{unionelip}.

\medskip

($\Leftarrow$) Let $F_j$ be as in \eqref{unionelip}. By Whitney's theorem we may find a family
$h_{ij}:\R{n}\rightarrow\R{n+m}$
of $C^1$-maps such that
\begin{equation}\label{unione4}
F_j(A_j)\subset E_j\cup\Big(\bigcup_{i=1}^{\infty}h_{ij}(\R{n})\Big),\quad \forall j\geq 1.
\end{equation}
Indeed we may choose $h_{ij}$ as in the statement of Whitney's theorem with $\varepsilon =\frac{1}{i}$. If $D_{ij}$ is the set in which $h_{ij}$ or $\nabla h_{ij}$ are different from $F_j$ or $\nabla F_j$ and if $D_j:=\cap_{i} D_{ij}$ it's clear that
$\M{L}^n(D)=0$ and, by the area formula, $\M{H}^n(F_j(D))=0$. Then set $E_j:=F(D_j)$ and we have \eqref{unione4}.

Set $C_{ij}:=\{x\in\R{n}:\rank h_{ij}(x)<n  \}.$ Then $\M{H}^n(h_{ij}(C_{ij}))=0$ by Sard's lemma. Set
$$N_0:=\Big(\bigcup_{j=1}^{\infty}E_j\Big)\cup \Big(\bigcup_{i,j=1}^{\infty}C_{ij}\Big).$$
Then $\M{H}^n(N_0)=0$ and
$$M\subset N_0\cup \Big(N_{ij} \Big),$$
with $N_{ij}:=h_{ij}(\R{n}\backslash C_{ij})$ countable union of $C^1$-submanifold thanks to the rank-max theorem ($N_{ij}$ is a $C^1$-submanifold if $h_{ij}$ is injective, otherwise we use the local injectivity of $h_{ij}$ to write $N_{ij}$ as countable union of $C^1$-submanifolds and a null set).
\end{proof}

\begin{cor} The image of a Lipschitz map
$$F:\Omega\subset\R{n}\rightarrow\R{n+m}$$ is a countable $n$-rectifiable set. In particular the graph of a Lipschitz function $u:\Omega\rightarrow\R{m}$ is $n$-rectifiable.
\end{cor}

Since the only rectifiable sets $\Sigma$ we will use are the graphs of Lipschitz function, we may assume w.l.o.g that $\M{H}^n\res\Sigma$ is locally finite, that is, for every compact set $K\subset\R{n+m}$,
$\M{H}^n(\Sigma\cap K)<\infty$.

\begin{defin}[Tangent plane]\label{pianotan} Given a countably $n$-rectifiable set $\Sigma$ in $\R{n+m}$ we define the \emph{tangent plane to $\Sigma$ in $p$}, if it exists, to be the only $n$-dimensional subspace $P$ in $\R{n+m}$ such that 
$$\lim_{\lambda\rightarrow 0}\int_{\eta_{p,\lambda}}f(y)d\M{H}^n(y)=\int_Pf(y)d\M{H}^n(y),\quad\forall f\in C^0_c(\R{n+m}),$$
where $\eta_{p,\lambda(y)}:=\lambda^{-1}(y-p)$ for every $y\in\R{n+m}$. Such plane $P$ will be denoted by $T_p\Sigma$.
\end{defin}

Given $\Sigma$ $n$-retcifiable in $\R{n+m}$, for instance a Lipschitz submanifold, its tangent plane is well defined $\M{H}^n$-a.e. It's clear that if $\Sigma$ is of class $C^1$, then the tangent plane just defined is the same as the tangent plane defined for smooth submanifolds as the set of tangent vectors. Given $\Sigma$ $n$-rectifiable, thanks to proposition \ref{carattrett}, for $\M{H}^n$-a.e. $p\in\Sigma$ there exists $N_{j(p)}$ $C^1$-submanifold such that
$p\in N_{j(p)}$. It may be seen that $T_pM=T_p N_{j(p)}$ for $\M{H}^n$-a.e. $p\in\Sigma$; in particular $T_p N_{j(p)}$ doesn't depend on the choice of the manifolds $N_j$ covering $\Sigma$, nor on the choice of $j(p)$.

For these reasons, given $U\subset\R{n+m}$ open and given $f\in\Lip(U)$, it's $\M{H}^n$-a.e. well defined in $\Sigma\cap U$ the gradient
$\nabla^\Sigma f:=\nabla^{N_j}f$. The latter is $\M{H}^n\res N_j$-a.e. well defined thanks to Rademacher's theorem.

\section{Rectifiable varifolds}

A rectifiable $n$-varifold is, roughly speaking, an $n$-rectifiable subset $\Sigma$ endowed with a multiplicity function $\theta$. The importance of the multiplicity rests on the necessity of defining a concept of limit in the space of varifolds under which the "area" is continuous . Consider the following example:
$$\Sigma_j:=\{0\}\times(0,1)\cup\Big\{\frac{1}{j}\Big\}\times(0,1)\subset\R{2}.$$
Each $\Sigma_j$ is a $C^\infty$-submanifold of $\R{2}$, but the only reasonable limit in the category of submanifolds is $\Sigma:=\{0\}\times(0,1)$. Were it so we would have $\M{A}(\Sigma_j)\rightarrow 2>\M{A}(\Sigma)$.
The limit in the sense of varifolds instead is $2\Sigma$ and its \emph{mass} is $2$.

\begin{defin}[Rectifiable varifolds] A rectifiable $n$-varifold with support in $\Sigma$ and multiplicity $\theta$, $V=\V(\Sigma,\theta)$, is the Radon measure (Borel regular measure finite of compact sets)
$$V:=\theta\M{H}^n\res\Sigma,$$
i.e.
$$V(A):=\int_{A\cap\Sigma}\theta(y)d\M{H}^n(y),\quad\forall A\subset\R{n+m}\;\mathrm{Borel},$$
where $\Sigma\subset\R{n+m}$ is $n$-rectifiable and $\theta$ is positive and locally integrable on $\Sigma$.
\end{defin}

\begin{oss} Equivalently we may see a rectifiable varifold as an equivalence class of couples
$(\Sigma,\theta)$ under the relation
\begin{equation}\label{eqequiv} 
(\Sigma_1,\theta_1)\sim(\Sigma_2,\theta_2)\; \mathrm{if}\; \M{H}^n(\Sigma_1\backslash\Sigma_2\cup
\Sigma_2\backslash\Sigma_1)=0\; \mathrm{and}\; \theta_1=\theta_2,\; \M{H}^n-q.o.
\end{equation}
Indeed it's clear that if \eqref{eqequiv} holds, then $\V(\Sigma_1,\theta_1)=\V(\Sigma_2,\theta_2)$; conversly if $V=\V(\Sigma_1,\theta_1)=\V(\Sigma_2,\theta_2)$ the support $\Sigma$ of $V\subset\Sigma_1\cap\Sigma_2$ satisfies $\M{H}^n(\Sigma\backslash\Sigma_i)=0$ because $\theta_i>0$ on $\Sigma_i$. Finally, it's obvious that $\theta_1=\theta_2$ $\M{H}^n-$a.e.

In any case we see a varifold as a Radon measure which may be expressed in the form $\V(\Sigma,\theta)$. 
\end{oss}

\begin{oss}
A rectifiable subset $\Sigma\subset\R{n+m}$ such that $\M{H}^n\res\Sigma$ is locally finite is a rectifiable varifold (in this case we identify, without further comments $\Sigma$ and $\M{H}^n\res\Sigma$).
\end{oss}

\begin{defin}[Tangent plane and mass]
Given a rectifiable varifold $V=\V(\Sigma,\theta)$, the tangent plane of $V$ in $p\in\Sigma$ is defined as
$$T_pV:=T_p\Sigma,$$
the latter being defined as in \ref{pianotan}
The definition is well posed $\M{H}^n$-a.e. and doesn't depend on $\Sigma$ but for a $\M{H}^n$-null set.

The mass of $V$ is its total variation in the sense of measures and is denoted by
$\mathbf{M}(V)$. Clearly
$$\mathbf{M}(V)=V(\R{n+m})=\int_{\Sigma}\theta d\M{H}^n.$$
\end{defin}

The convergence we are going to define on the space of rectifiable varifolds, different from the convergence in the sense of varifolds which we will define for abstract varifolds, is the weak* convergence induced by the duality between Radon measures and compactly supported continuous functions:

\begin{defin}[Weak convergence] We will say that a sequence of varifolds $V_j$ converges weakly to $V$ (and we will write $V_j\rightharpoonup V$) if
$$\lim_{j\rightarrow\infty}\int_{\R{n+m}}fdV_j=\int_{\R{n+m}}fdV,$$
for every $f\in C^0_c(\R{n+m}).$
\end{defin}

\begin{prop} The mass is continuous with respect to the weak convergence in a compact set $K\subset\R{n+m}$, i.e. if $V_j\rightharpoonup V$,  $\spt V_j\subset K$ for every $j\geq 0$ and $\spt V\subset K,$ then $\mathbf{M}(V_j)\rightarrow\mathbf{M}(V)$.
\end{prop}
\begin{proof} Set $R>0$ such that $K\subset B_R(0)$ and $\varphi\in C^0_{c}(\R{n+m})$ such that
$\varphi=1$ on $B_R(0)$. Then
$$\mathbf{M}(V_j)=\int_{\R{n+m}}\varphi dV_j\rightarrow \int_{\R{n+m}} dV=\mathbf{M}(V).$$
\end{proof}

\subsection{First variation of a varifold}

The concept of first variation, which we defined for $n$-dimensional $C^1$-submanifolds of $\R{n+m}$ in \ref{defvar}, may be easily extended to a rectifiable varifold $V=\V(\Sigma,\theta)$ thanks to the following definition:

\begin{defin}[Image varifold]\label{varretimg} Given $f:\R{n+m}\rightarrow\R{n+m}$ Lipschitz and proper\footnote{for every compact $K\subset\R{n+m}$ $f^{-1}(K)$ is compact.} and an $n$-rectifiable varifold $V=\V(\Sigma,\theta)$, the image varifold of $V$ under $f$ is defined by
$$f_{\#} V:=\V(f(\Sigma),\widetilde\theta),$$
where
$$\widetilde \theta(y)=\sum_{x\in\Sigma\cap f^{-1}(y)}\theta(x).$$
\end{defin}

Thanks to proposition \ref{carattrett}, $f(\Sigma)$ is rectifiable and, since $f$ is proper, we have that $\widetilde\theta\M{H}^n\res f(\Sigma)$ is locally finite: indeed, given a compact set $K$, by the area formula we get
$$f_{\#} V(K)=\int_{K\cap f(\Sigma)}\widetilde\theta d\M{H}^n=\int_{f^{-1}(K)\cap\Sigma}
Jf\theta d\M{H}^n,$$
$Jf:=\sqrt{\det (dF^*dF)}.$ The last integral is finite because $Jf$ is bounded, $f^{-1}(K)$ is compact and $\theta\M{H}^n\res\Sigma$ is locally finite.

\begin{defin}[First variation]\label{varprimavar}
Let $\varphi:\R{n+m}\times(-1,1)$ be of class $C^2$ and such that
\begin{enumerate}
\item there exists a compact set $K\subset\R{n+m}$ such that $\varphi_t(x)=x$ for every $x\notin K$;
\item $\varphi_0(x)=x$ for every $x\in\R{n+m}$.
\end{enumerate}
Then the first variation of a varifold $V$ with respect to $\varphi$ is the first variation of the mass of the family of varifolds $V_t:=(\varphi_t)_{\#}V$, that is
$$\frac{d}{dt}\Big|_{t=0}\mathbf{M}(V_t).$$ 
\end{defin}

With the same proof of proposition \ref{vardiv} we get

\begin{prop}\label{vardivvar} Consider a family of diffeomorphisms $\varphi_t$ as in definition \ref{varprimavar} and an $n$-rectifiable varifold $V=\V(\Sigma,\theta)$. Let
$$X(x):=\frac{d}{dt}\Big|_{t=0}\varphi_t(x)$$
be the first variation field of $\varphi.$ Then
\begin{equation}\label{varprimaform}
\frac{d}{dt}\Big|_{t=0}\mathbf{M}(V_t)=\int_{\Sigma}  \diver X dV= V(\diver X).
\end{equation}
\end{prop}

\begin{defin}[Minimal varifold]\label{varmindef}
We will say that an $n$-rectifiable varifold $V=\V(\Sigma,\theta)$ is minimal if its first variation is zero for every choice of $\varphi$ in \ref{varprimavar} or, equivalently, if for every vector field $X\in C^1_0(\R{n+m};\R{n+m})$ we have
\begin{equation}\label{diverzero}
\int_{\Sigma}\diver^\Sigma X dV=0.
\end{equation}
\end{defin}

In the case of a varifold defined by the graph of a Lipschitz function $u:\Cl\Omega\rightarrow\R{m}$, we want a definition of minimal varifold we fixed boundary. In general there is not a satisfactory definition for the boundary of a varifold, but in the case of a graph, we will give the following definition:

\begin{defin} A varifold whose support is the graph of a Lipschitz function $u:\Cl\Omega\rightarrow\R{m}$, $V=\V(\M{G}_u,\theta)$ is said to be minimal if
$$\int_{\Sigma}\diver^\Sigma X dV=0$$
for every vector field $X\in C^1_0(\Omega\times\R{m};\R{n+m})$.
Similarly we may say that $V$ is minimal in $\Omega\times\R{m}$.
\end{defin}

What we are requiring is that the mass of $V(\M{G}_u,\theta)$ is stationary with respect to variations contained in $\Omega\times\R{m}$, thus leaving fixed "the boundary" of the graph.

\subsection{The generalized mean curvature}

For a submanifold $\Sigma$ smooth and, thus, having a mean curvature, and for a variation $\varphi$ with field $X$, we have seen that
$$\frac{d}{dt}\Big|_{t=0}\M{A}(\Sigma_t)=-\int_\Sigma H\cdot Xd\M{H}^n=\int_\Sigma \diver^\Sigma Xd\M{H}^n.$$
A natural generalization of the mean curvature to the class of varifolds may be obtained by last equality:
 
\begin{defin}[Generalized mean curvature]
Given a varifold $V=\V(\Sigma,\theta)$, we will say that $V$ that has generalized mean curvature $H$ if
$$\int_{\Sigma} H\cdot X dV=-\int_{\Sigma} \diver^\Sigma X dV$$
for every vector field $X\in C^1_0(\R{n+m};\R{n+m})$.
\end{defin}

Thus a varifold is minimal if and only if it has zero generalized mean curvature.

\subsection{The monotonicity formula}

\begin{prop}\label{propmonot} Let be given a rectifiable $n$-varifold ($H=0$) $V=\V(\Sigma,\theta)$ in $U\subset\R{n+m}$. Then, for every $x_0\in\R{n+m}$, the function defined by
$$\rho\rightarrow \frac{V(B_\rho(x_0))}{\rho^n},\quad 0<\rho<d(x_0,U^c)$$
is monotone increasing.
\end{prop}

\begin{proof} Fix $\rho>0$ and define a function $\gamma\in C^1(\R{})$ such that
\begin{enumerate}
\item $\dot\gamma(t)\leq 0$ for every $t\geq 0$;
\item $\gamma(t)=1$ for every $t\leq\frac{\rho}{2};$
\item $\gamma(t)=0$ for every $t\geq \rho$.
\end{enumerate}
Consider the vector field
$$X(x):=\gamma(r)(x-x_0),\quad r:=|x-x_0|.$$
Let $x\in\Sigma$ be such that $T_x\Sigma$ exists; then the divergence on $\Sigma$ of $X$ at $x$ is well defined:
$$\diver^\Sigma X(x)=\sum_{j=1}^{n+m}e_j\cdot(\nabla^\Sigma X^j)=\gamma(r)
\sum_{j=1}^{n+m}e^{jj}+r\dot\gamma(r)\sum_{j,l=1}^{n+m}\frac{x^j-x_0^j}{r}\frac{x^l-x_0^l}{r}
e^{jl}$$
where $e^{jl}$ is the $(n+m)\times(n+m)$-matrix projecting $\R{n+m}$ onto $T_x\Sigma$.
The trace of the projection is $\sum e^{jj}=n$; moreover
$$\sum_{j,l=1}^{n+m}\frac{x^j-x_0^j}{r}\frac{x^l-x_0^l}{r}e^{jl}=|(Dr)^T|^2=1-|(Dr)^N|,$$
being equal to the scalar product between the projection of $Dr$ on $T_x\Sigma$ and $Dr=\frac{x-x_0}{r}$ itself.
This implies
$$\diver^\Sigma X(x)=n\gamma(r)+r\dot\gamma(r)(1-|(\nabla r)^N|^2).$$

Apply $\eqref{diverzero}$ to $X$ and get
\begin{equation}\label{monotconto}
n\int_\Sigma \gamma(r)dV+\int_\Sigma r\dot\gamma(r)dV=
\int_{\Sigma} r\dot\gamma(r)|(\nabla r)^N|^2dV.
\end{equation}

Now consider a family of functions $\gamma$ arising from a rescaling of the function $\Phi\in C^1(\R{})$ and satisfying
\begin{enumerate}
\item $\dot\Phi(t)\leq 0$ for every $t\geq 0$;
\item $\Phi(t)=1$ for every $t\leq\frac{1}{2};$
\item $\Phi(t)=0$ for every $t\geq 1$.
\end{enumerate}
More precisely, let $\gamma(r):=\Phi\bra{\frac{r}{\rho}}$ for a fixed $\rho>0$. It's clear that
$$r\dot\gamma(r)=\frac{r}{\rho}\dot\Phi\bra{\frac{r}{\rho}}=-\rho\frac{d}{d\rho}
\bra{\Phi\bra{\frac{r}{\rho}}}.$$
It follows that, defining
$$I(\rho):=\int_\Sigma \Phi\bra{\frac{r}{\rho}}dV,\quad
J(\rho)=\int_\Sigma \Phi\bra{\frac{r}{\rho}}\abs{(\nabla r)^N}^2 dV,$$
we obtain
$$n I(\rho)-\rho \dot I(\rho)=-\dot J(\rho),$$
which may be rewritten multipling by $\rho^{-n-1}$ as
\begin{equation}\label{contomonot2}
\frac{d}{d\rho}\bra{\frac{I(\rho)}{\rho^n}}=\frac{\dot J(\rho)}{\rho^n}.
\end{equation}
Let $\Phi$ converge from below to the characteristic function of $(-\infty,1]$ and obtain
$$I(\rho)\rightarrow V(B_\rho(x_0)),\quad J(\rho)\rightarrow \int_{B_\rho(x_0)}
\abs{(Dr)^N}^2dV,$$
whence, in the sense of distributions, \eqref{contomonot2} becomes
$$\frac{d}{d\rho}\bra{\frac{\mu_V(B_\rho(x_0))}{\rho^n}}=
\frac{d}{d\rho}\int_{B_\rho(x_0)}\frac{\abs{(Dr)^N}^2}{r^n}dV.$$
The integrand on the right is positive, whence the monotonicity of the term on the left.
\end{proof}

\section{Abstract varifolds}
Rectifiable varifolds are Radon measures in $\R{n+m}$. A compactness theorem for measures assures that a sequence of varifolds with equibounded masses admits a subsequence converging in the sense of measures. The limit, though, is a Radon measure whose support is, in general, non rectifiable. We are, thus, motivated to introduce a stronger notion of convergence and, eventually, a class of objects larger than the class of rectifiable varifolds.

\begin{defin}\label{defgrasm} Given an open set $U\subset\R{n+m}$, its Grassmannian fiber bundle of $n$-planes is
$$G_n(U):=U\times G(n,m), \quad \pi:G_n(U)\rightarrow U$$
where $G(n,m)\cong\frac{O(n+m)}{O(n)\times O(n)}$ is the Grassmannian of $n$-planes in $\R{n+m}$ and $\pi(x,S)=x$ for every $x\in U$ and every $n$-plane $S$. 
We endow $G_n(U)$ with the product topology induced by $U$ and $G_n(m)$.
\end{defin}

\begin{defin} An $n$-varifold in $U\subset\R{n+m}$ is a Radon measure $V$ on the Grassmannian fiber bundle $G_n(U)$. Associated to $V$ there is a measure $\mu_V$ on $U$ defined by
$$\mu_V(A):=V(\pi^{-1}(A)),\quad \forall A\subset U.$$ Finally we define the mass of $V$, $$\mathbf{M}(V):=\mu_V(U).$$
\end{defin}

\begin{oss} To show that the class of abstract varifolds contains the class of rectifiable varifolds, we observe that to a rectifiable $n$-varifold $\V(\Sigma,\theta)$ corresponds an abstract varifold $V$ defined by
$$V(A)=\V(\Sigma,\theta)(\pi(A\cap T\Sigma)),$$
being $T\Sigma:=\{(x,T_x\Sigma):x\in\Sigma_*\}$ the tangent bundle of $\Sigma$ ($\Sigma_*$ is the set of point of $\Sigma$ where the approximate tangent plane is defined). Clearly
$\mu_V=\V(\Sigma,\theta)$ because
$$\mu_V(A)=V(\pi^{-1}(A))=\V(\Sigma,\theta)(\pi(\pi^{-1}(A)\cap T\Sigma))=\V(\Sigma,\theta)(A\cap\Sigma). $$
\end{oss}

We give the space of $n$-dimensional varifolds in $U$ the weak* topology of the Radon measures, so that $V_n\rightarrow V$ if and only if for every $f\in C^1_c(G_n(U))$ we have
$$\int_{G_n(U)}f(x,S)dV_n(x,S)\rightarrow\int_{G_n(U)}f(x,S)dV(x,S).$$

\begin{oss} The convergence just defined, which we call \emph{convergence in the sense of varifolds}, is stronger than the convergence defined for rectifiable varifolds. It requires that, in a certain sense, both the support and the tangent planes of the varifolds in the sequence converge.
\end{oss}

\subsection{Image of a varifold and first variation}

\begin{defin} Given a Lipschitz map $\varphi:U\subset\R{n+m}\rightarrow U$ and given an $n$-varifold $V$, the imege varifold of $V$ relative to $\varphi$ is defined by
\begin{equation}\label{varimg0}
\varphi_\#V(A):=\int_{F^{-1}(A)}J\varphi(x,S)dV(x,S),
\end{equation}
where $F:G_n(U)\rightarrow G_n(U)$ is given by
$$F(x,S):=(\varphi(x),d\varphi_xS)$$
while
$$J\varphi(x,S):=\sqrt{\det\big((d\varphi_x\big|_S)^*d\varphi_x\big|_S\big)}.$$
\end{defin}

\begin{oss} The image varifold $\varphi_\#V$ may be defined using the duality with continuous functions on $G_n(U)$:
\begin{equation}\label{varimg}
\varphi_\#V(f)=\int_{G_n(U)}fd\varphi_\# V=\int_{G_n(U)}f(\varphi(x),d\varphi_x S)J\varphi(x,S)
dV(x,S).
\end{equation}
To pass from \eqref{varimg0} to \eqref{varimg} we may use the characteristic functions of subsets $A\subset G_n(U)$ and then use an approximation process.
\end{oss}

We define the first variation of a varifold in a way similar to that used for rectifiable varifolds: let $\varphi_t$ be as in definition \ref{varprimavar}. Then the first variation of a varifold $V$ with respect to $\varphi_t$ is
\begin{equation}\label{varprimaeq0}
\delta V(X):=\frac{d}{dt}\Big|_{t=0}\mathbf{M}(\varphi_{t\#}V),
\end{equation}
with $X(x):=\frac{\partial \varphi_t(x)}{\partial t}(x,0).$

With the same computation of propositions \ref{propvarprima} and following, it may be proved that
$$\delta V(X)=\int_{G_n(U)}\diver_S X(x)dV(x,S),$$
being
$$\diver_S X(x):=\sum_{i=1}^n\dual{\tau_i,\nabla_{\tau_i}X},$$
for a choice of an orthonormal basis $\{\tau_1,\ldots,\tau_n\}$ of $S$.

\begin{defin}\label{varprima2}
Given a varifold $V$ in $U\subset\R{n+m}$, its first variation (not respect to a vector field) in $W\subset U$ is
\begin{equation}\label{varprimaeq}
\norm{\delta V}:=\sup_{X\in C^1_c(U;\R{n+m})\atop |X|\leq 1,\spt X\subset W}|\delta V(X)|,
\end{equation}
where $|\delta V(X)|$ is defined in \eqref{varprimaeq0}.
\end{defin}

\begin{oss} If $V$ is the abtract varifold induced by a rectifiable varifold $\V(\Sigma,\theta)$, then $\varphi_\#V$ is the abstract varifold corresponding to $\varphi_\#\V(\Sigma,\theta)$ defined in \ref{varretimg}. For this reason the first variation of a rectifiable varifold is the same as the first variation of the corresponding abstract varifold.
\end{oss}

\subsection{Allard's compactness theorem}

Allard's compactness theorem answers the following question: when does a sequence of rectifiable integer multiplicity varifolds admit a subsequence converging in the sense of varifolds (i.e. on the Grassmannian) to an \emph{integer multiplicity rectifiable} varifold?

\begin{esempio} Consider the sequence of functions $u_n:[0,1]\rightarrow \R{}$ defined by
$$u_n(x)=\frac{\{nx\}}{n},$$
where $\{x\}$ denotes $x$ minus it's integral part.\footnote{for instance $\{\pi\}=0,14159265\ldots$}
The graph of $u_n$ is an i.m. rectifiable $1$-varifold in $\R{2}$, and as $n\rightarrow +\infty$,
the limit of $\V(G_{u_n},1)$ as rectifiable varifolds is
$\sqrt{2}\M{H}^1\res([0,1]\times\{0\}),$ whose corresponding abstract varifold is
$$\sqrt{2}\M{H}^1\res([0,1]\times\{0\})\times\delta_{0},$$
identifying a line in $\R{2}$ with the angle it spans with the $x$ axis.
On the other hand, the limit in the sense of measures on the Grassmannian is:
$$\sqrt{2}\M{H}^1\res([0,1]\times\{0\})\times\delta_{\frac{\pi}{4}},$$
which is not rectifiable.
\end{esempio}

It's not hard to prove that in the preceding example $\norm{\delta\M{G}_{u_n}}\rightarrow+\infty.$ This is why Allards compactness theorem doesn't apply to this example.

\begin{trm}[Compactness]\label{compallard}
Let be given a sequence of i.m. rectifiable varifolds $V_j$ in $U$ whose masses and first variations, as defined in \ref{varprima2}, are locally equibounded, that is such that for every $W\subset\subset U$
$$\sup_{j\geq 1}\bra{\mathbf{M}(V_j\big|_{W})+\norm{\delta V_j}(W)}<+\infty.$$ Let also $U$ be bounded.
Then there exists a subsequence $V_{j'}$ converging in the sense of varifolds to an i.m. rectifiable varifold and we have
$$\norm{\delta V}(W)\leq\liminf_{j\rightarrow+\infty}\norm{\delta V_j}, \quad\forall W\subset\subset U.$$
\end{trm}

\chapter{Allard's regularity theorems}\label{capitoloallard}

\section{Interior regularity}

The following theorem is due to Allard, who published it in 1972 \cite{allard1}; it reduces the study of the regularity of a minimal varifold to the study of its tangent cones and, consequently, to the study of rigidity theorems as Bernstein's theorem, in order to prove that the density of a minimal varifold is close to 1.

\begin{defin}[Density]\label{defdens}
Given a rectifiable varifold $V=\V(\Sigma,\theta)$, its density in $p$, if it exists, is the following limit:
\begin{equation}\label{eqdens}
\Theta^n(V,p):=\lim_{r\rightarrow 0}\frac{V(B_r(p))}{|B_r(p)|}.
\end{equation}
\end{defin}

\begin{oss} The density of a minimal varifold is always well defined because, thanks to the monotonicity formula, la quantity $\frac{V(B_r(p))}{|B_r(p)|}$ is monotone and, thus, has limit.
\end{oss}

\begin{trm}\label{regolallard}
Consider $U\subset\R{n+m}$ and let $V=\V(\Sigma,\theta)$ be a rectifiable minimal varifold. Then there exist $\delta,$ $\gamma$ and $c$ depending on $m$ and $n$ such that if
\begin{equation}\label{ipotesiallard}
\left\{  \begin{array}{ll}
0\in\spt V,& B_\rho(0)\subset U\\
\theta\leq 1 & V-\mathrm{q.o.},\\
\displaystyle\frac{V(B_\rho(0))}{\omega_n\rho^n}\leq 1+\delta,
\end{array}
\right. 
\end{equation}
then $\forall\alpha\in(0,1)$ there exists a linear isometry $q:\R{n+m}\rightarrow\R{n+m}$ and $u\in C^{1,\alpha}(\Cl{B}^n_{\gamma\rho}(0))$ with
$u(0)=0$,
$$V\res B_{\gamma\rho}(0)=\theta\M{H}^n\res\big(q(\M{G}_u)\cap B_{\gamma\rho}(0)\big).$$
Moreover
\begin{equation}\label{stimaallard2}
\frac{1}{\rho}\sup_{B_{\gamma\rho}(0)}|u|+\sup_{B_{\gamma\rho}(0)}|Du|+\rho [Du]_{\alpha,B_{\gamma\rho}(0)}
\leq c\delta^{\frac{1}{4n}}.
\end{equation}
\end{trm}

For a proof see the book of Leon Simon \cite{simon}.

\section{Boundary regularity}

Also this theorem is due, substantially to W. Allard \cite{allard2}. For more details see \cite{lawson}, theorem 2.3.

\begin{trm}\label{trmallard2}
Let $u$ be a solution of the Dirichlet problem for the minimal surface system \eqref{dirnonpar} with boundary data $\psi\in C^{s,\alpha}(\Cl\Omega)$, $2\leq s\leq+\infty$ and assume that $\Omega$ is \emph{strictly} convex. Then there exists a neighborhood $V$ of $\partial\Omega$ such that $u\in C^{s,\alpha}(V)$. If $\psi$ is analytic, then also $u$ is in a neighborhood of $\partial\Omega$.
\end{trm}

\end{document}